\documentclass[a4paper,10pt]{article}

\usepackage{latexsym}

\usepackage{amssymb}

\def \zgran{\displaystyle}
\def \zpizq{\left(}
\def \zcizq{\left[}
\def \zpder{\right)}
\def \zcder{\right]}

\def \za{\alpha}
\def \zb{\beta}
\def \zg{\gamma}
\def \zd{\delta}

\def \zl{\lambda}
\def \zm{\mu}

\def \zp{\pi}

\def \zr{\rho}

\def \zt{\tau}

\def \zf{\varphi}
\def \zfi{\phi}

\def \zq{\psi}
\def \zw{\omega}

\def \zL{\Lambda}

\def \zW{\Omega}

\def \zlma{\ell}

\def \zsu{\sum}
\def \zpr{\prod}

\def \zin{\cap}

\def \zun{\cup}

\def \zex{\wedge}

\def \zdi{\oplus}
\def \zdig{\bigoplus}
\def \zte{\otimes}

\def \zpor{\times}
\def \zci{\circ}

\def \zmei{\leq}
\def \zmai{\geq}
\def \zco{\subset}

\def \zcco{\supset}

\def \zpe{\in}

\def \zeq{\equiv}

\def \znoi{\neq}

\def \znope{\not\in}

\def \zpar{\partial}
\def \zinf{\infty}

\def \zfl{\rightarrow}

\def \zim{\Rightarrow}

\def \zbv{\mid}

\def \z/{\over}
\begin{document}

{\bf ON THE LOCAL THEORY OF VERONESE WEBS}
\bigskip

{\it \noindent Francisco-Javier Turiel

\noindent Geometr{\'\i}a y Topolog{\'\i}a, Facultad de Ciencias,
Campus de Teatinos, 29071 M\'alaga, Spain

\noindent e-mail: turiel@agt.cie.uma.es}
\vskip 1truecm

{\bf Abstract.}
This work is an introduction to the local geometric theory of Veronese webs
developed in the last twenty years. Among the different possible approach, here
one has chosen the point of view of differential forms. Moreover, in order to make
its reading easier, this text is self-contained in which directly regards
Veronese webs.
\bigskip

{\bf Introduction.}

The aim of this work is to provide an introduction to the local
theory of Veronese webs from the geometric viewpoint. Although
the classical theory is only developed on real manifolds  there is
no difficulty for extending it to complex ones as well, so both
case will be considered here. In our approach differential forms
play a crucial role, which will allow us to benefit from the
advantages of Cartan exterior differential calculus.

The notion of Veronese web, due to Gelfand and Zakharevich for the
case of codimension one \cite{GEA,GEB,GEC} and some years later  extended to any
codimension by Panasyuk and Turiel \cite{PAA}, \cite{TUC} (see \cite{ZA} as well),
is a tool for the study of
generic bihamiltonian structures in odd dimension and more
generally of Kronecker bihamiltonian structures. As it is well
known bihamiltonian structures, introduced by Magri in \cite{MAG}, are
related to some differential equations many of them with a
physical meaning. Therefore it seems interesting to describe this
geometrical objects.

With respect to the local aspect of this subject here, among other
results, one shows that:

\noindent 1) giving a generic bihamiltonian structure in odd
dimension is like giving a codimension one Veronese web (theorem
3.2),

\noindent 2) in the analytic category Kronecker bihamiltonian
structures and Veronese webs are locally equivalent (theorem 3.2
again; to point out that in codimension one we may utilize the
theorem on symmetric hyperbolic systems, therefore on real
manifold the $C^{\zinf}$ class is enough, while in codimension two
or more the  Cauchy-Kowalewsky theorem  and the analyticity are
needed).

Moreover a completely classification of $1$-codimensional Veronese
webs is exhibited (theorem 6.1). In higher codimension no local
classification is known but, in the analytic category, one gives a
versal model for Veronese webs.

On the other hand a link between classical $3$-webs and Veronese
webs is established in the example at the end of section 2 (see
\cite{BO} by Bouetou-Dufour too).

For the global aspect of the question, still widely open, the
reader may consult the papers by Rigal \cite{RIA,RIB,RIC}.

The present text consists of six sections and, in order to make its
reading easier, it is largely self-contained in which directly
regards Veronese webs. The first paragraph is devoted to the
algebraic theory including the classification of pairs of
bivectors (proposition 1.4). In the second one the notion of
Veronese web, illustrated with different examples, and its main
properties are discussed.

In the third section one associates a Veronese web to every
Kronecker bihamiltonian structure and conversely; moreover the
local equivalence between Kronecker bihamiltonian structures and
Veronese webs is established. The fourth and fifth paragraphs,
rather technical, are aimed to solve some exterior differential
systems needed elsewhere. Finally the sixth section contains the
local classification of $1$-codimensional Veronese webs and the
versal models for higher codimension.
\bigskip

{\bf 1. Algebraic theory}

The first part of this section is devoted to the study of the algebraic properties of
Veronese webs; in particular one gives a method for constructing any Veronese
web by means of an endomorphism of the support vector space.
The second part contains the classification of pairs of bivectors.

{\it All vector spaces considered here are real or complex.}
\bigskip

{\bf 1.1. Algebraic Veronese webs.}

Given an endomorphism $J$  and a subset $A$ of a vector space $V$,
the vector subspace spanned by $(A,J)$ will mean that one spanned by
$A\zun J(A)\zun J^{2}(A)\zun...$. When $A$ itself is a vector subspace
and $(A,J)$ spans $V$ we will say that the couple $(A,J)$ is {\it admissible}.
\bigskip

{\bf Lemma 1.1.} {\it If $(W,J)$ is admissible and $1\zmei dimV<\zinf$ then there exist
$H\zpe End(V)$ and a basis $\{ e_{1},...,e_{r}\}$ of $W$ such that:

\noindent (a) $H$ is nilpotent and  $Im(H-J)\zco W$.

\noindent (b) $V=\zdig_{j=1}^{r}U_{j}$ where each $U_j$ is the vector subspace spanned
by $(e_{j},H)$.

\noindent (c) The number of vector subspaces $U_j$ of dimension $\zmai\zlma$ equals
$dim(W+JW+...+J^{\zlma-1}W) - dim(W+JW+...+J^{\zlma-2}W)$ if $\zlma\zmai 2$ and $r$ if
$\zlma=1$.

Therefore the family of natural numbers $\{dimU_{j}\}$, $j=1,...,r$, only depends,
up to permutation,  on $J$ and $W$.}\bigskip

{\bf Proof.} First remark that
$W+JW+...+J^{k}W=W+{\tilde J}W+...+{\tilde J}^{k}W$
when ${\tilde J}=J+{\tilde G}$ and $Im{\tilde G}\zco W$. Therefore it is enough to prove
lemma 1.1 for some ${\tilde J}$; moreover (c) directly follows from (a) and (b) because $H$
is a particular case of ${\tilde J}$.

We will prove (a) and (b) by induction on $r=dimW$. Let $\zlma$ be the first natural number such that
$\zgran  dim\zpizq{\frac {W+JW+...+J^{\zlma}W}  {W+JW+...+J^{\zlma-1}W}}\zpder<r$. Then
there exists $e\zpe W-\{ 0\}$ such that $J^{\zlma}e$ belongs to $W+JW+...+J^{\zlma-1}W$; that
is to say $J^{\zlma}e=v_{0}+...+v_{\zlma-1}$ where each $v_{k}\zpe J^{k}W$:

Given a basis $\{ d_{1},...,d_{r}\}$ of $W$ set ${\tilde G}=\zsu_{j=1}^{r}d_{j}\zte \za_{j}$
with $\za_{1},...,\za_{r}\zpe V^{*}$. Then
$(J+{\tilde G})^{\zlma}=J^{\zlma}+J^{\zlma-1}\zci {\tilde G}+A$ where $ImA\zco W+JW+...+J^{\zlma-2}W$.
Hence $(J+{\tilde G})^{\zlma}e=v_{\zlma-1}+\zsu_{j=1}^{r}\za_{j}(e)J^{\zlma-1}d_{j}+v'$
where $v'$ belongs to $ W+JW+...+J^{\zlma-2}W$, which allows us to choose $\za_{1}...,\za_{r}$
in such a way that $(J+{\tilde G})^{\zlma}e=v'$. So by considering $J+{\tilde G}$
instead of $J$ and calling it $J$, we can suppose that $J^{\zlma}e$ belongs to
$ W+JW+...+J^{\zlma-2}W$.

Starting the process again with another ${\tilde G}=\zsu_{j=1}^{r}d_{j}\zte \za_{j}$, where
this time $\za_{1}(W)=...=\za_{r}(W)=0$,
one has  $(J+{\tilde G})^{\zlma}e=v_{\zlma-2}+\zsu_{j=1}^{r}\za_{j}(Je)J^{\zlma-2}d_{j}+v''$
with $v''\zpe( W+JW+...+J^{\zlma-3}W)$ and we may suppose that $J^{\zlma}e$ belongs to
$ W+JW+...+J^{\zlma-3}W$. Then we choose $\za_{1},...,\za_{r}$ such that $\za_{j}(W)=
\za_{j}(JW)=0$, $j=1,...,r$, and so one. In short we can assume $J^{\zlma}e=0$ without
loss of generality.

Let $U$ denote the vector subspace spanned by $(e,J)$. By the choice of $e$ the set
$\{ e,Je,...,J^{\zlma-1}e\}$ is a basis of $U$ and $dim(W\zin U)=1$. Let
$\zp:V\zfl {\frac {V} U}$ be the canonical projection and ${\bar J}$ the endomorphism
of ${\frac {V} U}$ induced by $J$. By the induction hypothesis, applied to ${\frac {V} U}$,
${\frac {W} U}$ and ${\bar J}$, there exist vectors $e_{1},...,e_{r-1}\zpe W$ and an
endomorphism ${\bar G}=\zsu_{j=1}^{r-1}\zp(e_{j})\zte\zb_{j}$ of ${\frac {V} U}$ such that
$\{\zp(e_{1}),...,\zp(e_{r-1})\}$ and ${\bar H}={\bar J}+{\bar G}$ are as in lemma 1.1.
Let $\zlma_{j}$, $j=1,...r-1$, be the dimension of the vector subspace spanned by
$(\zp(e_{j}),{\bar H})$. Since
$J+{\tilde G}$ where ${\tilde G}=\zsu_{j=1}^{r-1}e_{j}\zte(\zb_{j}\zci\zp)$ projects into
${\bar H}$ and ${\tilde G}(U)=0$, by calling $J$ to $J+{\tilde G}$, one may directly assume
${\bar H}={\bar J}$. Thus each
${\bar J}^{\zlma_{j}}\zp(e_{j})=0$ whence
$J^{\zlma_{j}}e_{j}=\zsu_{k=0}^{\zlma-1}a_{kj}J^{k}e$.

Now suppose $\zlma_{j}<\zlma$ for some $j$. Let $m$ be the biggest $k>\zlma_{j}$, if any, such that
$a_{kj}\znoi 0$. Then $J^{m}e$ belongs to $W+JW+...+J^{m-1}W$, which contradicts the definition of
of $\zlma$; so $a_{kj}=0$ when $k>\zlma_{j}$. But in this case $J^{\zlma_{j}}(e_{j}-a_{\zlma_{j}j}e)$
belongs to $W+JW+...+J^{\zlma_{j}-1}W$ which again contradicts the definition of $\zlma$ because
$\{ e_{1},...,e_{r-1}, e\}$ is a basis of $W$
and $e_{j}-a_{\zlma_{j}j}e\znoi 0$. In short $\zlma\zmei\zlma_{j}$, $j=1,...,r-1$.

Let $V'_{j}$, $j=1,...,r-1$, the vector subspace spanned by $\{ e_{j},...,J^{\zlma_{j}-1}e_{j}\}$.
As $\zp:  V'_{j}\zfl {\bar U}_{j}$ is an isomorphism, $\{ e_{j},...,J^{\zlma_{j}-1}e_{j}\}$ is a basis
of $V'_{j}$ and $V=V'_{1}\zdi...\zdi V'_{r-1}\zdi U$. Set $G=e\zte\za$ with $\za(U)=0$. Then
$(J+G)^{\zlma_{j}}e_{j}=\zsu_{k=0}^{\zlma-1}(a_{kj}+\za(J^{\zlma_{j}-k-1}e_{j}))J^{k}e$, which allows us
to choose $\za$ in such a way that $(J+G)^{\zlma_{j}}e_{j}=0$. For finishing it suffices considering
the basis $\{ e_{1},...,e_{r-1},e\}$ of $W$ and the endomorphism $H=J+G$. $\square$
\bigskip

{\bf Lemma 1.2.} {\it If $(W,J)$ is admissible, $dimV=n\zmai 1$ and $\{ w_{1},...,w_{r}\}$ is a basis
of $W$ then:

\noindent (a) The curve $\zg(t)=\zf(t)((J+tI)^{-1}w_{1})\zex...\zex ((J+tI)^{-1}w_{r})$
in  $\zL^{r} V$, where $\zf(t)$ is the characteristic polynomial of $-J$, is polynomial of
degree $n-r$.

More precisely there exists a basis $\{e_{ij}\}$, $i=1,...,n_{j}$ and $j=1,...,r$, of $V$ such that
$\zg(t)=\zg_{1}(t)\zex...\zex \zg_{r}(t)$ where every $\zg_{j}(t)=\zsu_{i=1}^{n_{j}}t^{i-1}e_{ij}$
and $e_{n_{1}1}\zex...\zex e_{n_{r}1}= w_{1}\zex...\zex w_{r}$.

\noindent (b) Let $(W,{\tilde J})$ be a second admissible couple. If $Im({\tilde J}-J)\zco W$ then
${\tilde\zg}(t)=\zg(t)$ where
${\tilde\zg}(t)={\tilde\zf}(t)(({\tilde J}+tI)^{-1}w_{1})\zex...\zex (({\tilde J}+tI)^{-1}w_{r})$
and ${\tilde\zf}(t)$ is the characteristic polynomial of $-\tilde J$.}
\bigskip

{\bf Proof.} Consider $H\zpe End(V)$ and a basis $\{ e_{1},...,e_{r}\}$ of $W$ like in lemma 1.1.
Set $n_{j}=dimU_{j}$ were $U_j$ is the vector subspace spanned by $(e_{j},H)$. By multiplying $e_1$
by a suitable scalar one can suppose that $w_{1}\zex...\zex w_{r}=e_{1}\zex...\zex e_{r}$, so
$((J+tI)^{-1}e_{1})\zex...\zex ((J+tI)^{-1}e_{r})=((J+tI)^{-1}w_{1})\zex...\zex ((J+tI)^{-1}w_{r})$,
which allows us to work with $e_{1},...,e_{r}$ instead of $w_{1},...,w_{r}$.

Note that $\{e_{ij}=(-1)^{n_{j}-i}H^{n_{j}-i}e_{j}\}$, $i=1,...,n_j$, $j=1,...,r$,
is a basis of $V$. Set
$\zr(t)=t^{n}((H+tI)^{-1}e_{1})\zex...\zex ((H+tI)^{-1}e_{r})$. As
$t^{n_{j}}(H+tI)^{-1}= \zsu_{i=1}^{n_{j}}(-1)^{n_{j}-i}t^{i-1}H^{n_{j}-i}$ on $U_j$, then
$\zr(t)=\zr_{1}(t)\zex...\zex\zr_{r}(t)$ where every
$\zr_{j}(t)=\zsu_{i=1}^{n_{j}}t^{i-1}e_{ij}$.

Let us see that $\zg(t)=\zr(t)$. Since $Im(J-H)\zco W$ one has
$((J+tI)\zex...\zex (J+tI))\zr(t)=\zq(t)e_{1}\zex...\zex e_{r}$
while the action of $J+tI$ on $\zl=e_{11}\zex...\zex
e_{n_{1}-1,1}\zex...\zex e_{1r}\zex...\zex e_{n_{r}-1,r}$ equals
$t^{n-r}\zl+\zsu_{j=1}^{r}e_{j}\zex\zm_{j}$ where each
$\zm_{j}\zpe\zL^{n-r-1}V$. The $n$-vector
$\zr(t)\zex\zl=t^{n-r}e_{1}\zex...\zex e_{r}\zex\zl$ is
transformed in $det(J+tI)t^{n-r}e_{1}\zex...\zex e_{r}\zex\zl$ by
$J+tI$. But calculating its action on $\zr(t)$ and $\zl$
separately shows that $\zr(t)\zex\zl$ is transformed in
$\zq(t)t^{n-r}e_{1}\zex...\zex e_{r}\zex\zl$ as well; whence
$\zq(t)=det(J+tI)$, which is the characteristic polynomial of
$-J$. Thus $((J+tI)\zex...\zex (J+tI))\zr(t)=\zf(t)e_{1}\zex...\zex
e_{r} =((J+tI)\zex...\zex (J+tI))\zg(t)$ and $\zr(t)=\zg(t)$.

A similar argument shows that ${\tilde\zg}(t)=\zr(t)$. $\square$

A polynomial curve $\zg$ in $\zL^{r}V$, $r\zmai 1$, is named a
{\it Veronese curve} if there exists a basis $\{ e_{ij}\}$,
$i=1,...,n_{j}$, $j=1,...,r$, of $V$ such that
$\zg(t)=\zg_{1}(t)\zex...\zex\zg_{r}(t)$ where each
$\zg_{j}(t)=\zsu_{i=1}^{n_{j}}t^{i-1}e_{ij}$.
When $r=1$ one obtains the classical notion of Veronese curve.

For convenience one will set $\zg(\zinf)=
\zgran lim {\frac {\zg(t)} {t^{n-r}}} $, when $t\zfl \zinf$.

Lemma 1.2 provides us a method for constructing Veronese curve for
which $\zg(\zinf)=w_{1}\zex...\zex w_{r}$.
Conversely given a Veronese curve $\zg$ in $\zL^{r}V$ and a basis like in the definition,
let $H$  and $W$ be the nilpotent endomorphism of $V$ defined by
$He_{ij}=-e_{i-1,j}$, $i\zmai 2$,
$He_{1j}=0$, and the vector subspace of basis
$\{ w_{1}=e_{n_{1}1},..., w_{r}=e_{n_{r}r}\}$ respectively. Then $(W,H)$ is admissible,
$n_{1},...,n_{r}$ are the natural numbers associated to $(W,H)$ by lemma 1.1, and
$\{w_{1},...,w_{r}\}$, $H$ give rise to $\zg$. Thus any Veronese curve can be
constructed through lemma 1.2.

Every $\zg(t)\zpe \zL^{r}V$ is decomposable and defines a $r$-dimensional vector subspace of $V$.
The union of all these vector subspaces spans $V$ since each $\zg_{j}({\mathbb K})$ spans the vector
subspace of basis $\{e_{ij}\}$, $i=1,..., n_{j}$. Now assume that
$\zg(t)=\zf(t)((J+tI)^{-1}w_{1})\zex...\zex ((J+tI)^{-1}w_{r})={\tilde\zf}(t)
(({\tilde J}+tI)^{-1}{\tilde w}_{1})\zex...\zex (({\tilde J}+tI)^{-1}{\tilde w}_{r})$; then
$\zg(\zinf)=w_{1}\zex...\zex w_{r}={\tilde w}_{1}\zex...\zex {\tilde w}_{r}$.

On the other hand the action of ${\tilde J}-J=({\tilde
J}+tI)-(J+tI)$ on $\zg(t)$ equals
$({\tilde\zf}(t)-\zf(t))w_{1}\zex...\zex w_{r}$; so ${\tilde J}-J$
maps the vector subspace defined by $\zg(t)$ into $W$. Hence
$Im({\tilde J}-J)\zco W$.

Obviously if $w_{1}\zex...\zex w_{r}={\tilde w}_{1}\zex...\zex {\tilde w}_{r}$ and
$Im({\tilde J}-J)\zco W$ then $w_{1}\zex...\zex w_{r}$, $J$ and
${\tilde w}_{1}\zex...\zex {\tilde w}_{r}$, $\tilde J$ define the same Veronese curve.

Two admissible couples $(W,J)$ and $({\tilde W},{\tilde J})$ are named {\it equivalent} if
$W={\tilde W}$ and $Im({\tilde J}-J)\zco W$. Clearly the family of natural numbers given
by lemma 1.1 is the same for equivalent couples.
From all that said previously follows:
\bigskip

{\bf Proposition 1.1.} {\it (a) Giving a Veronese curve in
$\zL^{r}V$, $r\zmai 1$, is like giving a class of equivalent
admissible couples $(W,J)$, where $dim W=r$, and an element
$w_{1}\zex...\zex w_{r}\zpe\zL^{r}W-\{ 0\}$, by setting
$\zg(t)=\zf(t)((J+tI)^{-1}w_{1})\zex...\zex ((J+tI)^{-1}w_{r})$,
where $\zf(t)$ is the characteristic polynomial of $-J$.

\noindent (b) Consider a Veronese curve $\zg(t)=\zg_{1}(t)\zex...\zex\zg_{r}(t)$ in $\zL^{r}V$
 and a basis $\{e_{ij}\}$, $i=1,...,n_{j}$, $j=1,...,r$, of $V$ such that
$\zg_{j}(t)=\zsu_{i=1}^{n_{j}}t^{i-1}e_{ij}$, $j=1,...,r$. Then, up to permutation, the family
of natural numbers $\{ n_{1},...,n_{r}\}$ only depends on $\zg$ and corresponds to the family
$\{ dim U_{j}\}$, $j=1,...,r$, given by lemma 1.1 applied to $(W,J)$.

\noindent (c) Two Veronese curves in $\zL^{r}V$ are isomorphic (through an isomorphism of $V$)
if and only if they have the same family of natural numbers $\{ n_{1},...,n_{r}\}$ up to permutation.}
\bigskip

{\bf Remark.} Any vector subspace of $\zL^{r}V$ containing a Veronese curve $\zg$ is at least of
dimension $n-r+1$ since $\zg(0),\zg^{(1)}(0),...,\zg^{(n-r)}(0)$ are linearly independent. Indeed
if $n_{1}=...=n_{r}=1$ it is obvious; otherwise assume, for example, $n_{1}\zmai 2$ and consider
a linear combination $\zsu_{\zlma=0}^{n-r}a_{\zlma}\zg^{(\zlma)}(0)=0$.

Let $\bar\zg$ denote the projection of $\zg$ into $\zL^{r}V'$ where $V'$ is the quotient of $V$ by
the line spanned by $e_{n_{1}1}$. Then $\bar\zg$ is a Veronese curve in $\zL^{r}V'$ of degree
$n-r-1$. As $\bar\zg(0),\bar\zg^{(1)}(0),...,\bar\zg^{(n-r)}(0)$ are the projections of
$\zg(0),\zg^{(1)}(0),...,\zg^{(n-r)}(0)$ and $\bar\zg^{(n-r)}(0)=0$, the induction hypothesis
implies that $a_{0}=...=a_{n-r-1}=0$. So $a_{n-r}\zg^{(n-r)}(0)=0$ whence $a_{n-r}=0$.

Now we will introduce the notion of Veronese web on a $n$-dimensional vector space $V$ with
$n\zmai 1$. A family $w=\{w(t)\zbv t\zpe {\mathbb K}\}$
of $(n-r)$-planes is called {\it a Veronese
web of codimension $r$} if there exists a Veronese curve $\zg$ in $\zL^{r}V^{*}$ such that
$w(t)= Ker\zg(t)$, $t\zpe {\mathbb K}$. The curve $\zg$ will be named a {\it representative
of} $w$.

If $\tilde\zg$ is another representative of $w$ then ${\tilde\zg}(t)=f(t)\zg(t)$ for any
$t\zpe {\mathbb K}$. As $\zg$ and $\tilde\zg$ are polynomial curves of degree $n-r$ and never lie
into a $(n-r-1)$-plane of $\zL^{r}V^{*}$, $f$ is constant and ${\tilde\zg}(t)=a\zg(t)$,
$a\zpe{\mathbb K}-\{0\}$. This allows us to define $w(\zinf)=Ker\zg(\zinf)$, which does not
depend on the representative. Moreover if $\{ \zb_{ij}\}$, $i=1,...,n_{j}$, $j=1,...,r$, is a basis
of $V^*$  such that $\zg(t)=\zg_{1}(t)\zex...\zex \zg_{r}(t)$ where each
$\zg_{j}(t)=\zsu_{i=1}^{n_{j}}t^{i-1}\zb_{ij}$, then
$w(\zinf)=Ker(\zb_{n_{1}1}\zex...\zex \zb_{n_{r}r})$.

In view of lemma 1.2 and proposition 1.1 one has:
\bigskip

{\bf Proposition 1.2.} {\it Consider on a $n$-dimensional vector space $V$ and a natural number
$1\zmei r\zmei n$.

\noindent (a) Given a $r$-codimensional vector subspace $W$ and an endomorphism $J$ both two of $V$,
if $(W',J^{*})$ spans $V^*$ where $W'$ is the annihilator of $W$ in $V^*$ then
$\zg(t)=\zf(t)((J+tI)^{-1})^{*}\zb$, where $\zf$ is the characteristic polynomial of $-J$ and
$\zb$ a $r$-form such that $Ker\zb=W$, represents a Veronese web $w$ of codimension $r$.

Moreover $lim_{t\zfl \zinf}t^{r-n}\zg(t)= \zb$, $w(\zinf)=W$ and $(J+tI)w(\zinf)=w(t)$ for
any $t\zpe\mathbb K$.

\noindent (b) Any Veronese web on $V$ of codimension $r$ may be represented in this way.

\noindent (c) Assume that $\zg(t)=\zf(t)((J+tI)^{-1})^{*}\zb$ and
${\tilde\zg}(t)={\tilde\zf}(t)(({\tilde J}+tI)^{-1})^{*}{\tilde\zb}$ represent two Veronese webs
$w$ and $\tilde w$ respectively. Then $w=\tilde w$ if and only if ${\tilde\zb}=a\zb$,
$a\zpe {\mathbb K}-\{0\}$, and $Ker({\tilde J}-J)\zcco w(\zinf)= {\tilde w}(\zinf)$.

In this last case ${\tilde\zg}=\zg$ if and only if ${\tilde\zb}=\zb$.

\noindent (d) Up to permutation the family of natural numbers $\{n_{1},...,n_{r}\}$, associated to
a splitting of a representative of a Veronese web $w$, only depends on $w$.
This family characterizes the Veronese web up to isomorphism.

By definition $n_{1},...,n_{r}$ will be called the {\rm the characteristic numbers of}
$w$ and their maximum {\rm the height of} $w$.}
\bigskip

{\bf Remark.} Often hereafter we will write $\zl(G,...,G)$ or $\zl\zci G$ instead
of $G^{*}\zl$ when $G$ is a morphism and $\zl$ a form.

On the other hand, note that $(W',J^{*})$ spans $V^*$ if and only if $W$ does not contain any
non-zero $J$-invariant vector subspace.

By (c) of proposition 1.2 the restriction of $J$ to $w(\zinf)$
gives rise to a morphism $\zlma:w(\zinf)\zfl V$ with no
$\zlma$-invariant vector subspace different from zero (this notion
is meaningful since $:w(\zinf)\zco V$) and which only depends on
the Veronese web $w$. Moreover $(\zlma+tI)w(\zinf)=w(t)$, $t\zpe
\mathbb K$, that is to say $\zlma^{*}\za=-t\za_{\zbv w(\zinf)}$ for
any $\za\zpe V^{*}$ such that $\za(w(t))=0$ and any $t\zpe \mathbb
K$. This last property characterizes $\zlma$ completely because
the union of the annihilators of $w(t)$, $t\zpe \mathbb K$, spans
$V^*$.

Conversely given a morphism $\zlma:W\zfl V$ whose only
$\zlma$-invariant vector subspace is zero, we may construct a
Veronese web by considering an endomorphism $J$ of $V$ such that
$J_{\zbv W}=\zlma$ and applying (a) of proposition 1.2 to it. This
Veronese web only depends on $\zlma$. In fact $w(t)=(\zlma+tI)W$.
Thus:

{\it Giving a Veronese web of codimension $r\zmai 1$ is equivalent
to giving a morphism $\zlma:W\zfl V$, where $W$ is a
$r$-codimensional vector subspace, without non-zero
$\zlma$-invariant vector subspaces.}
\bigskip

{\bf Proposition 1.3.} {\it Consider a Veronese web $w$ of
codimension $r\zmai 1$, a basis $\{\za_{1},...,\za_{n}\}$ of $V^*$
and scalars $a_{1},...,a_{n}$. Assume that $\za_{j}(w(-a_{j}))=0$,
$j=1,...,n$. Then $w$ can be constructed through (a) of
proposition 1.2 by means of the endomorphism $J$ defined by
$J^{*}\za_{j}=a_{j}\za_{j}$, $j=1,...,n$.}
\bigskip

{\bf Proof.} As $\zlma^{*}\za_{j}=a_{j}{\za_{j}}_{\zbv W}$ then
$\zlma^{*}=(J_{\zbv W})^{*}$, so $J$ is an extension of $\zlma$.
$\square$
\bigskip

{\bf Lemma 1.3.} {\it Consider a Veronese web $w$ of codimension
$r\zmai 1$ on a $n$-dimensional vector space $V$ and
its characteristic numbers $n_{1}\zmai...\zmai n_{r}$.
Let $k_{j}$ be the number of $n_{\zlma}$ greater than
or equal to $j$. Then $r=k_{1}\zmai...\zmai k_{n_{1}}\zmai 1$,
$k_{j}=0$ if $j>n_{1}$, and
$k_{1}+...+k_{n_{1}}=n$. Moreover:

\noindent (1) Given non-equal scalars $b_{1},...,b_{n-k},b$, where
$1\zmei k\zmei r$, there exists a basis $\{\za_{1},...,\za_{n}\}$
of $V^*$ such that $\za_{j}(w(b_{j}))=0$, $j=1,...,n-k$,
$\za_{j}(w(b))=0$, $j=n-k+1,...,n$.

\noindent (2) Given, this time, non-equal scalars
$c_{1},...,c_{n_{1}}$ there exists a basis $\{\zb_{ij}\}$,
$i=1,...,k_{j}$, $j=1,...,n_{1}$, of $V^*$ such that
$\zb_{ij}(w(c_{j}))=0$, $i=1,...,k_{j}$, $j=1,...,n_{1}$.}
\bigskip

{\bf Proof.} First consider a basis $\{e_{ij}^{*}\}$,
$i=1,...,n_{j}$, $j=1,...,r$ and $n_{1}\zmai...\zmai n_{r}$, of
$V^*$ such that $\zg(t)=\zg_{1}(t)\zex...\zex\zg_{r}(t)$, where
each $\zg_{j}(t)=\zsu_{i=1}^{n_{j}}t^{i-1}e_{ij}^{*}$, is a
representative of $w$. Now if $\zf:\{1,...,n-k\}\zfl\{1,...,r\}$
is a map such that $\zf^{-1}(\zlma)$ has $n_{\zlma}-1$ elements
when $1\zmei\zlma\zmei k$ and $n_{\zlma}$ otherwise, it suffices
to set $\za_{j}=\zg_{\zf(j)}(b_{j})$, $j=1,..,n-k$, and
$\za_{j}=\zg_{j+k-n}(b)$, $j=n-k+1,..,n$, for proving (1).

With regard to (2) set $\zb_{ij}=\zg_{i}(c_{j})$, $i=1,...,k_{j}$,
$j=1,...,n_{1}$. $\square$
\bigskip

{\bf 1.2. Pairs of bivectors.}

In this paragraph we will give the classification of pairs of
bivectors, due to Gelfand and Zakharevich, by regarding them as
quotients of symplectic pairs. Consider, on a finite dimensional
vector space $W$, a pair of bivectors $(\zL,\zL_{1})$. One defines
{\it the rank of} $(\zL,\zL_{1})$  as the maximum of ranks of
$(1-t)\zL+t\zL_{1}$, $t\zpe\mathbb K$. Note that
$rank((1-t)\zL+t\zL_{1})=rank(\zL,\zL_{1})$ except for a finite
number of scalars $t$, which is $\zmei {\frac{dimW} {2}}$ (they
are given by the polynomial equation $((1-t)\zL+t\zL_{1})^{k}=0$
where $rank(\zL,\zL_{1})=2k$). We will say that $(\zL,\zL_{1})$ is
{\it maximal (or of maximal rank)} if
$rank(\zL)=rank(\zL_{1})=rank(\zL,\zL_{1})$. Obviously if
$(\zL,\zL_{1})$ is not maximal one may choose
$\zL'=(1-a)\zL+a\zL_{1}$, $\zL_{1}'=(1-a_{1})\zL+a_{1}\zL_{1}$,
with $a\znoi a_{1}$, which is maximal. Consequently it suffices
classifying maximal pairs.

Recall that to any symplectic form $\zw$ defined on a vector space
$V$ of dimension $2n$ one can associate a dual bivector
$\zL_{\zw}$ by means of the isomorphism $v\zpe V\zfl
\zw(v,\quad)\zpe V^*$ (or $v\zpe V\zfl \zw(\quad,v)\zpe V^*$; the
result is the same). Conversely any bivector whose rank equals
$2n$ can be defined in this way. More generally when $\zL$ is a
bivector on $W$, considered as a bivector on
$Im\zL=\zL(W^{*},\quad)$ it is the dual of a symplectic form. Thus
every bivector can be described by its image and a symplectic form
on it; that is to say by the annihilator of $Im\zL$, or one of its
basis, and a $2$-form whose restriction to $Im\zL$ is symplectic.

Let $V_{0}$, $\zp:V\zfl{\frac{V} {V_{0}}}$ and $\zL$ be a vector
subspace of $V$, the canonical projection and the bivector on
${\frac{V} {V_{0}}}$ image of $\zL_{\zw}$ by $\zp$ respectively.
\bigskip

{\bf Lemma 1.4.} {\it Consider a second vector subspace $V_{1}$
such that $V=V_{0}\zdi V_{1}$.
Assume isotropic $V_0$. Let $\zL'$ be the bivector on $ V_{1}$
pull-back of $\zL$ by the isomorphism $\zp:V_{1}\zfl{\frac{V}
{V_{0}}}$. Then $\zL'$ is defined by $\zw(V_{0},\quad)_{\zbv
V_{1}}$ and $\zw_{\zbv V_{1}}$.}
\bigskip

{\bf Proof.} Set $dimV_{0}=n-k$. There exists a basis
$\{e_{1},...,e_{2n}\}$ of $V$ such that
$\zw=\zsu_{j=1}^{n}e_{2j-1}^{*}\zex e_{2j-1}^{*}$ and
$\{e_{2j-1}\}$, $j=k+1,...,n$, is a basis of $V_0$. Then
$\zL=\zp(e_{1})\zex\zp(e_{2})+...+\zp(e_{2k-1})\zex\zp(e_{2k})$.

On the other hand, as $V=V_{0}\zdi V_{1}$ there exists a basis
${\mathcal B}=\{e_{1}+v_{1},...,e_{2k}+v_{2k},
\{e_{2j}+v_{2j}\}_{j=k+1,...,n}\}$ of $V_1$ where every $v_{i}\zpe
V_{0}$. Obviously $\zL'=(e_{1}+v_{1})\zex (e_{2}+v_{2})+...+
(e_{2k-1}+v_{2k-1})\zex (e_{2k}+v_{2k})$.

The restriction to $V_{1}$ of the family $\{e_{j}^{*}\}$,
$j=1,...,2k$ and $j=2(k+1),...,2n$, is the dual basis of $\mathcal
B$. So $\zL'$ will be defined by the restriction to $V_{1}$ of the
$2$-form $e_{1}^{*}\zex e_{2}^{*}+...+ e_{2k-1}^{*}\zex
e_{2k}^{*}$, which equals that of $\zw$, and by the basis
$\{{e_{2j}^{*}}_{\zbv V_{1}}={\zw(e_{2j-1},\quad)}_{\zbv V_{1}}\}$,
$j=k+1,...,n$ of the annihilator of $Im\zL'$  . $\square$

Warning lemma 1.4 can fail if $V_{0}$ is not isotropic. For
example on ${\mathbb K}^4$: $\zw=e_{1}^{*}\zex
e_{2}^{*}+e_{3}^{*}\zex e_{4}^{*}$, $V_{0}={\mathbb
K}\{e_{3},e_{4}\}$ and $V_{1}={\mathbb
K}\{e_{1}+e_{3},e_{2}+e_{4}\}$.

{\bf Remark.} On a finite dimensional vector space $E$ consider a
symplectic form $\zW$ and a $2$-form $\zW_{1}$. Let $K$ be the
endomorphism defined by $\zW_{1}=\zW(K,\quad)$, that is to say
$\zW_{1}(v,w)=\zW(Kv,w)$, $v,w\zpe E$. Then
$\zW(K,\quad)=\zW(\quad, K)$; thus every  $\zW(K^{k},\quad)$ is a
$2$-form on $E$. By definition the characteristic polynomial, the
minimal one and the elementary divisors of $(\zW,\zW_{1})$ will be
those of $K$.

Suppose that the characteristic polynomial of $(\zW,\zW_{1})$ is
the product $p_{1}p_{2}$ of two monic relatively prime
polynomials. Then $(\zW,\zW_{1}, E)$ can be identified to the
product of two similar structures $(\zW^{1},\zW_{1}^{1},
E_{1})\zpor (\zW^{2},\zW_{1}^{2}, E_{2})$ where $p_i$ is the
characteristic polynomial of $(\zW^{i},\zW_{1}^{i})$, $i=1,2$. In
this way classifying $(\zW,\zW_{1})$ reduces to the case where the
characteristic polynomial is a power of an irreducible polynomial.
It is not difficult to see that the model of $(\zW,\zW_{1})$ is
completely determined by the Jordan structure of $K$. Moreover
every elementary divisor occurs an even number of times, so $p$ is
the square of another polynomial, and the minimal polynomial
divides the square root of $p$.

Let us come back to the main question. Consider a second
symplectic form $\zw_1$ on $V$, the dual bivector $\zL_{\zw_{1}}$
and its image $\zL_{1}$ by $\zp$ on ${\frac{V} {V_{0}}}$. Let $J$
be the endomorphism (in fact the automorphism ) of $V$ defined by
$\zw_{1}=\zw(J,\quad)$.
\bigskip

{\bf Lemma 1.5.} {\it Assume that $V_0$ is isotropic for both
$\zw$ and $\zw_1$ and $(\zL,\zL_{1})$ is maximal. Then the vector
subspace spanned by $(V_{0},J)$ is $\zw$ and $\zw_1$ isotropic.}
\bigskip

{\bf Proof.} First note that $rank\zL=rank\zL_{1}=2r$, where
$dimV_{0}=n-r$, because $V_{0}$ is bi-isotropic. On the other hand
if $rank(\zL_{\zw}+t\zL_{\zw_{1}})=2n$ then
$\zw((I+tJ^{-1})^{-1},\quad)$ is its dual symplectic form (recall
that if $\zW_{1}=\zW(K,\quad)$ then
$\zL_{1}=\zL((K^{-1})^{*},\quad)$ when $\zL$ and $\zL_{1}$ are
regarded as $2$- forms on the dual space). Since
$rank(\zL+t\zL_{1})\zmei 2r$ this implies that $V_0$ is isotropic
for $\zw((I+tJ^{-1})^{-1},\quad)$, so $\zw((I+tJ^{-1})^{-1}v,w)=0$
for any $v,w\zpe V_{0}$.

Near $0\zpe \mathbb K$ one has $rank(\zL_{\zw}+t\zL_{\zw_{1}})=2n$
so deriving at $t=0$ successively yields, up multiplicative
constant, $\zw(J^{-k}v,w)=0$, $k\zmai 0$. Hence
$\zw(J^{-\zlma}v,J^{-s}w)=0$ for any $\zlma,s\zmai 0$ as
$\zw(J,\quad)=\zw(\quad,J)$. This implies that the vector subspace
spanned by $(V_{0},J^{-1})$ is $\zw$-isotropic; but this last one
equals the vector subspace spanned by $(V_{0},J)$ since $J$ is
invertible.

Finally as our vector subspace is $J$-invariant it has to be
$\zw_1$-isotropic. $\square$

For the remainder of this paragraph $(\zL,\zL_{1})$ will be a
maximal pair of bivectors defined on $m$-dimensional vector space
$W$. Set $r=corank(\zL,\zL_{1})$. Assume that $\zL$ is defined by
$\za_{1},...,\za_{r},{\tilde\zw}$, and $\zL_1$ by
$\zb_{1},...,\zb_{r},{\tilde\zw}_{1}$, where $\za_{1},...,\za_{r},
\zb_{1},...,\zb_{r}\zpe W^{*}$ and
${\tilde\zw},{\tilde\zw}_{1}\zpe\zL^{2} W^{*}$.

Let $V_0$ be a vector space of dimension $r$ and
$\{e_{1},...,e_{r}\}$ one of its basis. Let
$\{e_{1}^{*},...,e_{r}^{*}\}$ denote the extension of the dual
basis of $\{e_{1},...,e_{r}\}$ to $V=W\zdi V_{0}$ by setting
$e_{i}^{*}(W)=0$, $i=1,...,r$. On the other hand we will regard
$\za_{1},...,\za_{r},
\zb_{1},...,\zb_{r},{\tilde\zw},{\tilde\zw}_{1}$ as forms on $V$
such that $\za_{i}(V_{0})=\zb_{i}(V_{0})=0$, $i=1,...,r$,
${\tilde\zw}(V_{0},\quad)={\tilde\zw}_{1}(V_{0},\quad)=0$. Now on
$V$ one considers the symplectic forms $\zw={\tilde\zw}+
\za_{1}\zex e_{1}^{*}+...+\za_{r}\zex e_{r}^{*}$ and
$\zw_{1}={\tilde\zw}_{1}+\zb_{1}\zex e_{1}^{*}+...+\zb_{r}\zex
e_{r}^{*}$. If we identify $W$ to ${\frac{V} {V_{0}}}$ by means of
the canonical projection, by lemma 1.4 the pair $(\zL,\zL_{1})$ is
just the image of the dual pair $(\zL_{\zw},\zL_{\zw_{1}})$. Thus
any maximal pair is the quotient of a symplectic pair by a
bi-isotropic vector subspace.

By technical reasons we will deform $\zw_{1}$ for simplifying the
algebraic structure of the symplectic pair. Set
$\zw_{\zm}=\zw_{1}+\zb_{1}\zex\zm_{1}+...+\zb_{r}\zex\zm_{r}=
{\tilde\zw}_{1}+\zb_{1}\zex (e_{1}^{*}+\zm_{1})+...+\zb_{r}\zex
(e_{r}^{*}+\zm_{r})$ where $\zm_{1},...,\zm_{r}\zpe V^{*}$, which
is symplectic if and only if $\{(e_{1}^{*}+\zm_{1})_{\zbv
V_{0}},...,(e_{r}^{*}+\zm_{r})_{\zbv V_{0}}\}$ is still a basis of
$V_{0}^{*}$. In this last case $V_{0}$ is $\zw_{\zm}$ isotropic
and the dual bivector $\zL_{\zm}$ projects into $\zL_{1}$ as well
(apply lemma 1.4 again). Let $J$ and $J_{\zm}$ be the
endomorphisms defined by $\zw_{1}=\zw(J,\quad)$ and
$\zw_{\zm}=\zw(J_{\zm},\quad)$ respectively, and let ${\bar
e}_{j}$ be the vector defined by $\zw({\bar e}_{j},\quad)=\zm_{j}$,
$j=1,...,r$. Then $J_{\zm}=J+\zsu_{j=1}^{r}({\bar
e}_{j}\zte\zb_{j} +Je_{j}\zte\zm_{j})$.

Therefore, since ${J_{\zm}}_{\zbv V_{0}}=\zsu_{j=1}^{r}Je_{j}\zte
(e_{j}^{*}+\zm_{j})_{\zbv V_{0}}$, the form $\zw_{\zm}$ is
symplectic, that is to say $J_{\zm}$ is an isomorphism, if and
only if ${J_{\zm}}_{\zbv V_{0}}$ is a monomorphism.

Let $V_{1}$ denote the vector subspace spanned by $(V_{0},J)$, and
$V_{2}$ the $\zw$-orthogonal of $V_{1}$. As $JV_{1}=V_{1}$ the
vector subspace $V_{2}$ is the $\zw_1$-orthogonal of $V_{1}$ too.
From lemma 1.5 follows that $V_{1}$ is isotropic for $\zw$ and
$\zw_1$; thus $V_{1}\zco V_{2}$ and $\zb_{j}(V_{2})=0$,
$j=1,...,r$, since $\zb_{j}(V_{2})=-\zw_{1}(e_{j},V_{2})$. Hence
$J_{\zm}=J+\zsu_{j=1}^{r}Je_{j}\zte\zm_{j}$ on $V_{2}$ and
$(J_{\zm}-J)V_{2}\zco JV_{0}$.

{\it Hereafter assume $\zw_{\zm}$ symplectic}. Then
$J_{\zm}V_{0}=JV_{0}$. This implies that $(V_{0},J_{\zm})$ spans $V_1$ as well.
Again lemma 1.5, this time applied to $\zw,\zw_{\zm}$, shows that $V_1$
is $\zw_{\zm}$-isotropic; moreover $V_2$ is the $\zw_{\zm}$-orthogonal of
$V_1$ because $J_{\zm}V_{1}=V_{1}$. Obviously $JV_{2}=J_{\zm}V_{2}=V_{2}$ since
$JV_{1}=J_{\zm}V_{1}=V_{1}$ and $V_2$ is the orthogonal of $V_1$ for $\zw$,
$\zw_1$ and $\zw_{\zm}$.

The restricted forms $\zw_{\zbv V_{2}}$ and
${\zw_{1}}_{\zbv V_{2}}={\zw_{\zm}}_{\zbv V_{2}}$ (recall that $\zb_{j}(V_{2})=0$
so ${({\zw_{\zm}}-\zw_{1})}_{\zbv V_{2}}=0$) project into a pair
$({\bar\zw},{\bar\zw}_{1})$ of symplectic forms on ${\frac{V_{2}} {V_{1}}}$. As
$\zw_{1}=\zw(J,\quad)$ and $\zw_{\zm}=\zw(J_{\zm},\quad)$, the endomorphism
$\bar J$ of ${\frac{V_{2}} {V_{1}}}$ defined by
${\bar\zw}_{1}={\bar\zw}({\bar J},\quad)$ is just the projection of both
$J_{\zbv V_{2}}$ and ${J_{\zm}}_{\zbv V_{2}}$.

The next step will be to control the characteristic polynomial of $J_{\zm}$,
which is the product of three characteristic polynomials: that of the projection
of $J_{\zm}$ on ${\frac{V} {V_{2}}}$, that of ${J_{\zm}}_{\zbv V_{1}}$ and
that of the projection of $J_{\zm}$ on ${\frac{V_{2}} {V_{1}}}$. This last one
is the characteristic polynomial of $\bar J$, therefore it does not depend on
$\zm$; it will denote by $\zq(t)$.

As $J$ is an isomorphism $V_1$ is also the vector subspace spanned by
$(JV_{0},J)$. Now from lemma 1.1 applied to $V_1$ and $(JV_{0},J)$ follows the
existence of a nilpotent $H\zpe End(V_{1})$, such that
$Im(H-J_{\zbv V_{1}})\zco JV_{0}$, and a basis $\{d_{1},...,d_{r}\}$ of
$JV_0$ such that $V_{1}=\zdi_{j=1}^{r}U_{j}$ where each $U_j$ is the vector subspace
spanned by $(d_{j},H)$. Set $G=H+\zsu_{j=1}^{r}d_{j}\zte\zl_{j}$ where
$\zl_{1},...,\zl_{r}\zpe V_{1}^{*}$ and $\zl_{j}(U_{i})=0$ if $i\znoi j$. Then
we may choose $\zl_{1},...,\zl_{r}$ in such a way that $(d_{j},G)$ spans $U_{j}$,
$j=1,...,r$, $Im(G-J_{\zbv V_{1}})\zco JV_{0}$ and the characteristic polynomial
of $G_{\zbv U_{j}}$ is any monic polynomial whose degree equals the dimension of
$U_j$. Moreover if $G$ is invertible, so a monomorphism, there exist $\zw_{\zm}$
and $J_{\zm}$ such that  ${J_{\zm}}_{\zbv V_{1}}=G$ since
$J_{\zm}=J+\zsu_{j=1}^{r}Je_{j}\zte\zm_{j}$ on $V_2$.

Consider non-equal and non-zero scalars $a_{1},...,a_{k}$, where $k=dimV_{1}$, which
are not roots of $\zq(t)$. Then we can suppose, without loss of generality, that
$(d_{j},J_{\zm})$ spans $U_j$ and the characteristic polynomial of ${J_{\zm}}_{\zbv U_{j}}$
equals $\zpr_{i\zpe I_{j}}(t-a_{i})$ where $\{1,...,k\}$ is the disjoint union of
$I_{1}$,..., $I_{r}$. Thus the characteristic polynomial ${\zq}_{\zm}(t)$ of $J_{\zm}$
equals $\zq(t)\zr(t)\zpr_{i=1}^{k}(t-a_{i})$ where $\zr(t)$ is the characteristic
polynomial of the projection of $J_{\zm}$ on ${\frac{V} {V_{2}}}$.
But ${\zq}_{\zm}(t)$ has to be a square and $a_{1},...,a_{k}$ are not roots of $\zq(t)$,
so ${\zq}_{\zm}(t)=\zq(t)\zpr_{i=1}^{k}(t-a_{i})^{2}
=\zq(t)\zpr_{j=1}^{r}(\zpr_{i\zpe I_{j}}(t-a_{i})^{2})$.

Now we may identify $(\zw,\zw_{\zm},V)$ to a product
$\zpr_{j=0}^{r}(\zt_{j},\zt'_{j},L_{j})$ in such a way that $\zq(t)$ is the
characteristic polynomial of ${J_{\zm}}_{\zbv L_{0}}$ and
$\zpr_{i\zpe I_{j}}(t-a_{i})^{2}$ that of ${J_{\zm}}_{\zbv L_{j}}$,
$j=1,...,r$.
Then $V_{0}\zin L_{0}=\{0\}$, $dim(V_{0}\zin
L_{j})=1$, $j=1,...,r$, and $V_{0}=\zdi_{j=1}^{r}(V_{0}\zin
L_{j})$; indeed $J_{\zm}^{-1}d_{j}$ is a basis of $V_{0}\zin
U_{j}$, since $J_{\zm}V_{0}=JV_{0}$, and
$(J_{\zm}^{-1}d_{j},J_{\zm})$ spans $U_j$. Remark that
$\zpr_{i\zpe I_{j}}(t-a_{i})$ is the minimal polynomial of any
$v\zpe V_{0}\zin L_{j}-\{0\}$, $j=1,...,r$. Moreover
$(\zL,\zL_{1})$ is identified, in a natural way, to the product
of the dual pair $(\zL_{\zt_{0}},\zL_{\zt'_{0}})$, {\it called
symplectic}, times the projections of the dual pairs
$(\zL_{\zt_{j}},\zL_{\zt'_{j}})$ on $\frac{L_{j}} { V_{0}\zin
L_{j}}$, $j=1,...,r$, which will be called the {\it Kronecker
elementary pairs}. The case without symplectic factor and that
with no Kronecker elementary factor happen.

Let us describe the Kronecker elementary pair in dimension $2n-1$.
Consider, on a $2n$-dimensional vector space $E$, a pair of
symplectic forms $(\zW,\zW_{1})$ and the endomorphism $K$ defined
by $\zW_{1}=\zW(K,\quad)$. Suppose that
$\zpr_{i=1}^{n}(t-b_{i})^{2}$ is the characteristic polynomial of
$(\zW,\zW_{1})$, where all $b_{i}\znoi 0$ and $b_{i}\znoi b_{j}$
if $i\znoi j$. Let $E_0$ be a $1$- dimensional vector subspace of
$E$ such that the minimal polynomial of its non-zero elements is
$\zpr_{i=1}^{n}(t-b_{i})$. Then there exists a basis
$\{e_{1},...,e_{2n}\}$ of $E$ such that $\zW=e_{1}^{*}\zex
e_{2}^{*}+...+e_{2n-1}^{*}\zex e_{2n}^{*}$, $\zW_{1}=b_{1}e_{1}^{*}\zex
e_{2}^{*}+...+b_{n}e_{2n-1}^{*}\zex e_{2n}^{*}$ and
$e=-\zsu_{j=1}^{n}e_{2j}$ is a basis of $E_0$.

Denote by $E_1$, ${\tilde\zL}$ and ${\tilde\zL}_{1}$ the vector
subspace of basis $\{e_{1},...,e_{2n-1}\}$, and the images of
$\zL_{\zW}$ and $\zL_{{\zW}_{1}}$ on $\frac E {E_{0}}$
respectively. As $E=E_{0}\zdi E_{1}$ by lemma 1.4 the bivector
${\tilde\zL}$, considered on $E_1$ identified to $\frac E {E_{0}}$
in the natural way, is given by $\tilde\zw=
\zsu_{j=1}^{n-1}e_{2j-1}^{*}\zex e_{2j}^{*}$,
$\za=\zsu_{j=1}^{n}e_{2j-1}^{*}$ (obviously both of them
restricted to $E_1$) while ${\tilde\zL}_{1}$ is described by
${\tilde\zw}_{1}= \zsu_{j=1}^{n-1}b_{j}e_{2j-1}^{*}\zex
e_{2j}^{*}$, $\zb=\zsu_{j=1}^{n}b_{j}e_{2j-1}^{*}$. Moreover,
since $\zL_{\zW}+t\zL_{\zW_{1}}$ is the dual bivector of
$\zW((I+tK^{-1})^{-1},\quad)$
when $t\zpe {\mathbb K}-\{-b_{1},...,-b_{n}\}$,
the bivector ${\tilde\zL}+t{\tilde\zL}_{1}$ is given by $\zm_{t}=
\zsu_{j=1}^{n-1}b_{j}(t+b_{j})^{-1}e_{2j-1}^{*}\zex e_{2j}^{*}$
and $\za_{t}=(\zpr_{j=1}^{n}(t+b_{j}))\zW((I+tK^{-1})^{-1}e,\quad)
=\zsu_{j=1}^{n}b_{j}(\zpr_{i=1;i\znoi
j}^{n}(t+b_{i}))e_{2j-1}^{*}$.

But $\zm_{-b_{n}}$, $\za_{-b_{n}}$ still define a bivector on
$E_1$, which by continuity has to be equal to
${\tilde\zL}-b_{n}{\tilde\zL}_{1}$. Thus
$corank({\tilde\zL}+t{\tilde\zL}_{1})=1$ for
any $t\zpe{\mathbb K}-\{-b_{1},...,-b_{n-1}\}$.
Reasoning in the same way but with
other suitable direct summands of $E_0$ (for example for $-b_{1}$
the vector subspace spanned by $\{e_{2},...,e_{2n}\}$) finally
shows that $corank({\tilde\zL}+t{\tilde\zL}_{1})=1$,
$t\zpe{\mathbb K}$. Hence $Im({\tilde\zL}+t{\tilde\zL}_{1})=
Ker\za_{t}$, $t\zpe{\mathbb K}$.

Therefore $E'=\zin_{t\zpe\mathbb K}
Im({\tilde\zL}+t{\tilde\zL}_{1})$ is the $(n-1)$-dimensional
vector subspace of basis $\{e_{2j}\}$, $j=1,...,n-1$, and setting
$w(t)=\frac {Im({\tilde\zL}+t{\tilde\zL}_{1})} {E'}$ defines a
Veronese web $w$ of codimension one on $\frac {E_{1}} {E'}$.
Indeed, identify $\frac {E_{1}} {E'}$ to the vector subspace $E''$
spanned by  $\{e_{2j-1}\}$, $j=1,...,n$, and restrict $\za_t$ to
it (proposition 1.2 applied to $K_{\zbv E''}$ and
$(\zsu_{j=1}^{n}b_{j}e_{2j-1}^{*})_{\zbv E''}$
just yields ${\za_{t}}_{\zbv E''}$).

On the other hand $({\tilde\zL},{\tilde\zL}_{1})$ is a particular
case of $(\zL,\zL_{1})$ with $r=1$. So
$({\tilde\zL},{\tilde\zL}_{1})$ is isomorphic to a product of a
possible symplectic pair in dimension $2(n-k)$ and a Kronecker
elementary pair associated to scalars $a_{1},...,a_{k}$. As
$dim(Im({\tilde\zL}+t{\tilde\zL}_{1}))=2n-2$, $t\zpe{\mathbb K}$,
the characteristic polynomial of the symplectic factor has no
roots and in this case an elementary calculation yields $dimE'=
2n-k-1$. But $dimE'=n-1$ so $k=n$; that is to say there is no
symplectic factor. In other words our pair can be constructed from
any family of non-zero scalars $\{a_{1},...,a_{n}\}$ such that
$a_{i}\znoi a_{j}$ if $i\znoi j$, which shows that {\it the
Kronecker elementary pair $({\tilde\zL},{\tilde\zL}_{1})$ only
depends on the dimension $2n-1$ but not on $\{b_{1},...,b_{n}\}$.
Thus, up to isomorphism, in every odd dimension there exists just
one Kronecker elementary pair}.

Now we may state:
\bigskip

{\bf Proposition 1.4.} {\it Consider a maximal pair of bivectors
$(\zL,\zL_{1})$ on a finite dimensional vector space $W$. Set
$r=corank(\zL,\zL_{1})$. Let $L_0$ be the intersection of all the
vector subspaces $Im(\zL+t\zL_{1})$, ${t\zpe\mathbb K}$, of
codimension $r$. Denote by $L'_0$ its annihilator in $W^*$. One
has:

\noindent (a) $L_{0}\zco Im\zL_{1}$ and
$\zL(L'_{0},\quad)=\zL_{1}(L'_{0},\quad)$.

In what follows set $L_{1}=\zL(L'_{0},\quad)$.

\noindent (b) The restrictions to $L_0$ of the $2$-forms
associated to $\zL$ and $\zL_{1}$ respectively, which are unique
since $L_{0}\zco Im\zL\zin Im\zL_{1}$, have $L_1$ as kernel.

Therefore the projections on $\frac{L_{0}} {L_{1}}$ of these
restricted $2$-forms,denoted by $\bar\zw$ and $\bar\zw_{1}$
respectively, are symplectic.

\noindent (c) Setting $w(t)=\frac{Im(\zL+t\zL_{1})} {L_{0}}$, $t\zpe \mathbb K$,
defines a Veronese web on $\frac W {L_{0}}$.

\noindent (d) The elementary divisors of $({\bar\zw},{\bar\zw_{1}})$ and the
characteristic numbers $n_{1}\zmai ...\zmai n_{r}$of $w$
determine the algebraic structure of $(\zL,\zL_{1})$ completely. More precisely
$(\zL,\zL_{1},W)$ is isomorphic to a product
$\zpr_{\zlma=0}^{r}(\zL^{\zlma},\zL_{1}^{\zlma},W^{\zlma})$ where
$(\zL^{0},\zL_{1}^{0},W^{0})$ is isomorphic, in its turn, to the dual pair of
$({\bar\zw},{\bar\zw_{1}},{\frac {L_{0}} {L_{1}}})$ and every
$(\zL^{\zlma},\zL_{1}^{\zlma},W^{\zlma})$, $\zlma=1,...,r$, is the Kronecker
elementary pair in dimension $2n_{\zlma}-1$.

\noindent (e) $corank(\zL+a\zL_{1})>r$ if and only if $-a$ is a root of the
characteristic polynomial of $({\bar\zw},{\bar\zw_{1}})$.}
\bigskip

{\bf Remark.} Let $(\zw,\zw_{1})$ be a pair of symplectic forms on a $2n$-dimensional
vector space $V$ and let $V_0$ be a line in $V$. Denote by $V_1$ the vector subspace
spanned by $(V_{0},J)$ where $\zw_{1}=\zw(J,\quad)$. Then the dimension of the
symplectic factor, given by proposition 1.4, of the pair $(\zL,\zL_{1})$ induced by
$(\zL_{\zw},\zL_{\zw_{1}})$ on $\frac V {V_{0}}$ equals $2(n-dimV_{1})$.

Indeed, first note that
$rank\zL=rank\zL_{1}=rank(\zL,\zL_{1})=2n-2$ so $(\zL,\zL_{1})$ is
maximal. Let $e$ and $W$ be a basis and a direct summand of $V_0$
respectively. Then $\zw={\tilde\zw}+\za\zex e^{*}$,
$\zw_{1}={\tilde\zw}_{1}+\zb\zex e^{*}$ where $Ker\za$, $Ker\zb$,
$Ker{\tilde\zw}$ and $Ker{\tilde\zw}_{1}$ contain $V_0$,
$Kere^{*}=W$ and $e^{*}(e)=1$. Therefore, after identifying $W$
and $\frac V {V_{0}}$, bivectors $\zL$ and $\zL_{1}$ are given by
${\tilde\zw},\za$ and  ${\tilde\zw}_{1},\zb$ respectively. As
$V=W\zdi V_{0}$ we are just in the situation which allowed us
splitting any maximal pair. There it was showed that the dimension of
the symplectic factor equals that of $\frac {V_{2}} {V_{1}}$,
where $V_2$ was the orthogonal of $V_1$; in our case
$2n-2dimV_{1}$.
\bigskip

{\bf Proposition 1.5.} {\it Let $W$ be a $(2n-1)$-dimensional
vector space. The action of the linear group $GL(W)$ on
$(\zL^{2}W)\zpor(\zL^{2}W)$ possesses one dense open orbit whose
model is the elementary Kronecker pair in dimension $2n-1$.}
\bigskip

{\bf Proof.} First let us show that any pair $(\zL,\zL_{1})$ is
approachable in $(\zL^{2}W)\zpor(\zL^{2}W)$ by a Kronecker
elementary one. As bivectors of rank $2n-2$ are generic in
$\zL^{2}W$ one can suppose $rank\zL=rank\zL_{1}=2n-2$. Now assume
that the symplectic factor given by proposition 1.4 applied to
$(\zL,\zL_{1})$ has dimension $2k\zmai 2$ and minimal polynomial
$\zf$. Note that there is only one Kronecker elementary factor
since $corank(\zL,\zL_{1})=1$.
By constructing this Kronecker elementary factor with
scalars $\{a_{1},...,a_{n-k}\}$ which are not roots of $\zf$, the
pair $(\zL,\zL_{1})$ becomes the quotient by a line $V_0$ of a
dual symplectic pair $(\zL_{\zw},\zL_{\zw_{1}})$ defined on a
$2n$-dimensional vector space $V$, in such a way that the minimal
polynomial of $(\zw,\zw_{1})$ is $\zf\zpr_{j=1}^{n-k}(t-a_{j})$
and $\zpr_{j=1}^{n-k}(t-a_{j})$ that of each $e\zpe V_{0}-\{0\}$.
In particular $(e,J)$, where $\zw_{1}=\zw(J,\quad)$, spans a
$(n-k)$-dimensional vector subspace.

Set $V=W\zdi V_{0}$. By lemma 1.4 $(\zL,\zL_{1})$, regarded on
$W$, is given by $\zw_{\zbv W}$, $\zw(e,\quad)_{\zbv W}$,
${\zw_{1}}_{\zbv W}$ and ${\zw_{1}(e,\quad)}_{\zbv W}$. Now
consider a vector $e'$ near $e$ whose minimal polynomial is
$\zf\zpr_{j=1}^{n-k}(t-a_{j})$. Then $(e',J)$ spans a vector
subspace of $V$ of dimension $>n-k$ and the symplectic factor of
the quotient of $(\zL_{\zw},\zL_{\zw_{1}})$ by ${\mathbb K}\{e'\}$
has dimension $<2k$ (see the foregoing remark). Since
$V=W\zdi{\mathbb K}\{e'\}$ this last pair is given on $W$ by
$\zw_{\zbv W}$, $\zw(e',\quad)_{\zbv W}$, ${\zw_{1}}_{\zbv W}$ and
${\zw_{1}(e',\quad)}_{\zbv W}$; therefore we can choose it as
close to $(\zL,\zL_{1})$ as desired and, after a
finite number of steps, $(\zL,\zL_{1})$ will be
approached by a Kronecker pair.

On the other hand if $(\zL',\zL'_{1})$ is a Kronecker elementary
pair, consider scalars $\{a_{1},...,a_{n}\}$ all of them
different. Then $dim(Im(\zL'+a_{j}\zL'_{1}))=2n-2$, $j=1,...,n$,
and $dim(\zin_{j=1}^{n}Im(\zL'+a_{j}\zL'_{1}))=n-1$. Therefore
when $(\zL,\zL_{1})$ is close enough to $(\zL',\zL'_{1})$ one has
$dim(Im(\zL+a_{j}\zL_{1}))=2n-2$, $j=1,...,n$, and
$dim(\zin_{j=1}^{n}Im(\zL+a_{j}\zL_{1}))=n-1$. But by (d) of
proposition 1.4 this last dimension equals $n-k-1$ where $2k$ is
the dimension of the symplectic factor of $(\zL,\zL_{1})$; so
$k=0$ and $(\zL,\zL_{1})$ is Kronecker elementary too. $\square$
\bigskip

{\bf 2. Veronese webs on manifolds}

This section contains the basic theory of Veronese webs of any codimension.
The notion of Veronese web of codimension one was introduced by Gelfand and
Zakharevich for studying the generic bihamiltonian structures on odd
dimensional manifolds \cite{GEA,GEB,GEC}. Later on Panasyuk and Turiel dealt
with the case of higher codimension \cite{PAA}, \cite{TUC}; see \cite{ZA} as well.
The approach given from now on, different from that of Gelfand, Zakharevich and
Panasyuk, follows the Turiel's work \cite{TUA,TUB,TUC}.

{\it Hereafter all structures considered will be real $C^{\zinf}$ or complex
holomorphic unless another thing is stated.}

Let $N$ be a real or complex manifold of dimension $n$. A family
$w=\{w(t)\zbv t\zpe \mathbb K\}$ of involutive distributions (or foliations) on
$N$ of codimension $r\zmai 1$ is named a {\it Veronese web of codimension $r$},
if for any $p\zpe N$ there  exist an open neighborhood $A$ of this point and
a curve $\zg(t)$ in the module of sections of $\zL^{r}T^{*}A$ (that is to
say  $\zg(t)(q)\zpe\zL^{r}T_{q}^{*}A=\zL^{r}T_{q}^{*}N$ for every $q\zpe A$) such that:

\noindent 1) $w(t)=Ker\zg(t)$, $t\zpe {\mathbb K}$, on $A$

\noindent 2) for each $q\zpe A$, $\zg(t)(q)$ is a Veronese curve in
$\zL^{r}T^{*}N$.

The curve $\zg$ is called a {\it (local) representative of $w$}.

Although curves $\zg(t)(q)$ and $\zg(t)(q')$ could be not isomorphic
when $q\znoi q'$, $\zg(t)=\zsu_{i=0}^{n-r}t^{i}\zg_{i}$ where
$\zg_{0},...\zg_{n-r}$ are differentiable $r$-forms on $A$.
On the other hand $Ker\zg_{n-r}$ is an involutive distribution of
dimension $n-r$ since each $Ker\zg(t)$ was integrable and
$limt^{r-n}\zg(t)=\zg_{n-r}$, $t\zfl\zinf$. This allows us to define
$w(\zinf)=Ker\zg_{n-r}$, which does not depend on the representative
because if ${\tilde\zg}$ is another representative  then
${\tilde\zg}=f\zg$ on the common domain (see (c) of proposition 1.2).
In particular, there exists a global
representative if and only if $w(\zinf)$ is transversally orientable.
Obviously $w$ as map from ${\mathbb K}\zun\{\zinf\}\zeq{\mathbb K}P^{1}$
to the Grassmann manifold of $(n-r)$-plans of $TN$ is smooth.

{\bf Examples.} 1) On $S^3$ regarded as a Lie group consider three left invariant
contact forms $\zr_{1},\zr_{2},\zr_{3}$. Suppose that
$\zr_{1}\zex\zr_{2}\zex\zr_{3}\znoi 0$ and set $\zg(t)=(\zr_{1}+t\zr_{2})\zex\zr_{3}$.
Then $\zg$ defines a codimension two Veronese web which is not flat because
$Ker\zr_{3}=w(0)\zdi w(\zinf)$ is a contact structure.

\noindent 2) On ${\mathbb K}^4$ with coordinates $(x_{1},x_{2},y_{1},y_{2})$ set
$\zg(t)=(dx_{2}\zex dy_{2}+x_{2}dx_{1}\zex dx_{2})
+t(x_{2}dx_{2}\zex dy_{1}-dx_{1}\zex dy_{2})+t^{2}dy_{1}\zex dy_{2}$.
Then $\zg$ defines a Veronese web of codimension two since
$d\zg(t)=0$ and $\zg(t)=(-dx_{1}+x_{2}^{-1}dy_{2}+tdy_{1})\zex(-x_{2}dx_{2}+tdy_{2})$
when $x_{2}\znoi 0$, while $\zg(t)=(dx_{2}-tdx_{1}+t^{2}dy_{1})\zex dy_{2}$ if
$x_{2}=0$.

Note that $\zg(t)(q)$ and $\zg(t)(q')$ are not isomorphic as Veronese curves
when $q_{2}\znoi 0$ and $q'_{2}=0$.

\noindent 3) Let $V$ be the $3$-dimensional Lie algebra spanned by the vectors fields on
$\mathbb K$: $X_{1}=(\zpar/\zpar t)$, $X_{2}=t(\zpar/\zpar t)$
and $X_{3}=t^{2}(\zpar/\zpar t)$. Set ${\tilde w}(t)
=\{ v\zpe V\zbv v(t)=0\}$. As  ${\tilde w}(t)=Ker\{ e_{1}^{*}+te_{2}^{*}
+t^{2}e_{3}^{*}\}$ where $\{ e_{1}^{*},e_{2}^{*},e_{3}^{*}\}$ is the dual basis of
$\{ X_{1},X_{2},X_{3}\}$, ${\tilde w}=\{{\tilde w}(t)\zbv t\zpe \mathbb K\}$ is an
algebraic Veronese web on $V$. But $V$ is isomorphic to the Lie algebra of
$SL(2,{\mathbb K})$ and each ${\tilde w}(t)$ is a subalgebra of $V$; therefore ${\tilde w}$
gives rise to a Veronese web $w$  of codimension one on any $3$-dimensional homogeneous
space of $SL(2,{\mathbb K})$.

Now we shall give a local description of Veronese webs of codimension $r$ by means of
a $(1,1)$ tensor field and a $r$-form. Consider non-equal scalars $\{a_{1},...,a_{n-k},a\}$,
where $1\zmei k\zmei r$, and any point $p\zpe N$. By lemma 1.3 there exists a basis
$\{\zl_{1},...,\zl_{n}\}$ of $T_{p}^{*}N$ such that $Ker\zl_{j}\zcco w(-a_{j})(p)$,
$j=1,...,n-k$, and $Ker\zl_{j}\zcco w(-a)(p)$, $j=n-k+1,...,n$. Since every distribution
$w(t)$ is involutive, on some open neighbourhood $A$ of $p$ one may construct closed
$1$-forms $\zb_{1},...,\zb_{n}$, extensions of $\zl_{1},...,\zl_{n}$, such that
$Ker\zb_{j}\zcco w(-a_{j})$, $j=1,...,n-k$, $Ker\zb_{j}\zcco w(-a)$, $j=n-k+1,...,n$, and
$\zb_{1}\zex...\zex\zb_{n}$ is a volume form.

Let $J$ be the $(1,1)$ tensor field on $A$ defined by $\zb_{j}\zci J=a_{j}\zb_{j}$,
$j=1,...,n-k$, and $\zb_{j}\zci J=a\zb_{j}$, $j=n-k+1,...,n$. In coordinates
$(x_{1},...,x_{n})$ such that $\zb_{j}=dx_{j}$, $j=1,...,n$, one has
$J=\zsu_{j=1}^{n-k}a_{j}{\frac{\partial} {{\partial}x_{j}}}\zte dx_{j}
+a\zsu_{j=n-k+1}^{n}{\frac{\partial} {{\partial}x_{j}}}\zte dx_{j}$, so $J$ is
flat and diagonalizable.

Moreover, by propositions 1.2 and 1.3, if $\zb$ is a $r$-form on $A$ such that
$Ker\zb=w(\zinf)$ then
$\zg(t)=(\zpr_{j=1}^{n-k}(t+a_{j}))(t+a)^{k}((J+tI)^{-1})^{*}\zb$ is a
representative of $w$.

On the other hand if $n_{1}$ is the height of $w(p)$ and
$k_{1}\zmai...\zmai k_{n_{1}}$ are like in lemma 1.3, given non-equal scalars
${\tilde a}_{1},...,{\tilde a}_{n_{1}}$ a similar argument allows to construct closed
$1$-forms $\{{\tilde \zb}_{ij}\}$, $i=1,...,k_{j}$, $j=1,...,n_{1}$, linearly independent
everywhere and such that ${\tilde \zb}_{ij}(w(-{\tilde a}_{j}))=0$, $i=1,...,k_{j}$,
$j=1,...,n_{1}$. Then, by propositions 1.2 and 1.3,
$\zg(t)=(\zpr_{j=1}^{n_{1}}(t+{\tilde a}_{j})^{k_{j}})(({\tilde J}+tI)^{-1})^{*}\zb$ where
${\tilde J}$ is defined by ${\tilde \zb}_{ij}\zci{\tilde J}={\tilde a}_{j}{\tilde \zb}_{ij}$,
$i=1,...,k_{j}$, $j=1,...,n_{1}$.
\bigskip

{\bf Theorem 2.1.} {\it Let $N$ be a $n$-dimensional real or complex manifold.

\noindent (1) Consider a Veronese web $w$ on $N$ of codimension $r$ and non-equal scalars
$a_{1},...,a_{n-k},a$ where $1\zmei k\zmei r$. Then for each $p\zpe N$ there exist an open
set $p\zpe A$ and a $(1,1)$-tensor field $J$ on $A$ whit characteristic polynomial
$\zf(t)=(\zpr_{j=1}^{n-k}(t-a_{j}))(t-a)^{k}$, which is flat and diagonalizable, such that:

\noindent (I) $(Ker(J^{*}-a_{j}I))w(-a_{j})=0$, $j=1,...,n-k$, and
$(Ker(J^{*}-aI))w(-a)=0$.

\noindent (II) For any $q\zpe A$, $(w(\zinf)(q)',J^{*}(q))$ spans $T^{*}_{q}A$, that is to
say $w(\zinf)(q)$ contains no $J$-invariant vector subspace different from zero
(as before $\,'\,$ means the annihilator).

In particular, if $\zb$ is a $r$-form and $Ker\zb=w(\zinf)$ then
$\zg(t)=(\zpr_{j=1}^{n-k}(t+a_{j}))(t+a)^{k}((J+tI)^{-1})^{*}\zb$ represents $w$.

Moreover is $\zl$ is a closed $1$-form such that $Ker\zl\zcco w(\zinf)$ then
$d(\zl\zci J)_{\zbv w(\zinf)}=0$.

\noindent (2) Now consider non-equal scalars ${\tilde a}_{1},...,{\tilde a}_{n_{1}}$ instead
of $a_{1},...,a_{n-k},a$, where $n_{1}$ is the height of $w(p)$, and numbers
$k_{1}\zmai...\zmai k_{n_{1}}$ like in lemma 1.3. Then there exists a $(1,1)$-tensor field
${\tilde J}$ defined on an open neighbourhood $\tilde A$, which is flat and diagonalizable,
with characteristic polynomial ${\tilde\zf}=\zpr_{j=1}^{n_{1}}(t-{\tilde a}_{j})^{k_{j}}$
such that $(Ker({\tilde J}^{*}-{\tilde a}_{j}I))w(-{\tilde a}_{j})=0$,
$j=1,...,n_1$, $(w(\zinf)',{\tilde J}^{*})$ spans $T{\tilde A}^{*}$,
${\tilde\zg}(t)=\zpr_{j=1}^{n_{1}}(t+{\tilde a}_{j})^{k_{j}}(({\tilde J}+tI)^{-1})^{*}\zb$
represents $w$ and $\zg={\tilde\zg}$ on $A\zin \tilde A$.

Moreover $d(\zl\zci {\tilde J})_{\zbv w(\zinf)}=0$ for any $1$-closed form $\zl$ such that
$Ker\zl\zcco w(\zinf)$.

\noindent (3) Finally, on $N$ consider a foliation $\mathcal F$ of codimension $r\zmai 1$,
a $r$-form ${\bar\zb}$ such that $Ker{\bar\zb}=\mathcal F$ and $(1,1)$-tensor field
${\bar J}$ with characteristic polynomial ${\bar\zf}(t)$. Suppose that:

\noindent (I) $({\mathcal F}',{\bar J}^{*})$ spans $T^{*}N$, that is to say
${\mathcal F}$ does not contain any non-zero ${\bar J}$-invariant vector subspace.

\noindent (II) ${(N_{\bar J})}_{\zbv\mathcal F}=0$, where $N_{\bar J}$ is the Nijenhuis torsion
of ${\bar J}$, and $d(\zm\zci{\bar J})_{\zbv\mathcal F}=0$ for each
closed $1$-form $\zm$ such that $Ker\zm\zcco\mathcal F$ (note that if
${\mathcal F}=
Ker(\zl_{1}\zex...\zex\zl_{r})$ where each $\zl_j$ is a closed $1$-form, this last
condition is satisfied if and only if
$\zl_{1}\zex...\zex\zl_{r}\zex d(\zl_{j}\zci{\bar J})=0$, $j=1,...,r$).

Then ${\bar\zg}(t)=(-1)^{n}{\bar\zf}(-t)(({\bar J}+tI)^{-1})^{*}{\bar\zb}$ defines
a Veronese web $\bar w$ of codimension $r$ for which ${\bar w}(\zinf)={\mathcal F}$.
This Veronese web only depends on ${\mathcal F}$ and ${\bar J}$.}
\bigskip

In view of propositions 1.2 and 1.3 and all that said previously, for proving
theorem 2.1 it suffices to show that $d(\zl\zci J)_{\zbv w(\zinf)}
=d(\zl\zci {\tilde J})_{\zbv w(\zinf)}=0$ and that every ${\bar w}(t)$, $t\zpe\mathbb K$,
is involutive. For this purpose we need the following result:
\bigskip

{\bf Lemma 2.1.} {\it Given a $1$-form $\zr$ and a $(1,1)$-tensor field $G$ on manifold,
then $d(\zr\zci G)(G\quad,\quad)+d(\zr\zci G)(\quad,G\quad)=d\zr(G\quad,G\quad)
+d(\zr\zci G^{2})+\zr\zci N_{G}$.}
\noindent

{\bf Proof.} Consider two vector fields $X,Y$. One has:

$d(\zr\zci G)(GX,Y)=(GX)\zr(GY)-Y\zr(G^{2}X)-\zr(G[GX,Y])$

$d(\zr\zci G)(X,GY)=X\zr(G^{2}Y)-(GY)\zr(GX)-\zr(G[X,GY])$.

So $d(\zr\zci G)(GX,Y)+d(\zr\zci G)(X,GY)=d\zr(GX,GY)
+d(\zr\zci G^{2})(X,Y)+\zr(N_{G}(X,Y))$. $\square$

Let $\zl$ a closed $1$-form such that $Ker\zl\zcco w(\zinf)$. If
$t\zpe {\mathbb K}-\{-a_{1},...,-a_{n-k},-a\}$ lemma 2.1 applied to
$\zl\zci(J+tI)^{-1}$ and $(J+tI)$ yields
$d(\zl\zci J)=d(\zl\zci (J+tI))=-d(\zl\zci(J+tI)^{-1})((J+tI)\quad,
(J+tI)\quad)$.

But $Ker(\zl\zci(J+tI)^{-1})$ contains $w(t)$ which is involutive,
so $d(\zl\zci(J+tI)^{-1})_{\zbv w(t)}=0$ . Hence
$d(\zl\zci J)_{\zbv w(\zinf)}=-d(\zl\zci(J+tI)^{-1})((J+tI)\quad,
(J+tI)\quad)_{\zbv w(\zinf)}=0$ since $(J+tI)w(\zinf)=w(t)$.

The case of $\tilde J$ is similar.

Now we shall prove the involutivity of every ${\bar w}(t)$. Consider a point
$q\zpe N$ and $t\zpe\mathbb K$ such that ${\bar J}+tI$ is invertible around
$q$. If $\zm$ is a closed $1$-form and $Ker\zm\zcco\mathcal F$ then
$Ker(\zm\zci({\bar J}+tI)^{-1})\zcco{\bar w}(t)$ and by lemma 2.1

\noindent $d(\zm\zci({\bar J}+tI)^{-1})({\bar w}(t),{\bar w}(t))=
d(\zm\zci({\bar J}+tI)^{-1})(({\bar J}+tI){\mathcal F},({\bar J}+tI){\mathcal F})
=-d(\zm\zci({\bar J}+tI))({\mathcal F},{\mathcal F})
-\zm(N_{({\bar J}+tI)}({\mathcal F},{\mathcal F}))
=-d(\zm\zci {\bar J})({\mathcal F},{\mathcal F})
-\zm(N_{{\bar J}}({\mathcal F},{\mathcal F}))=0$.

That is to say $d(\zm\zci({\bar J}+tI)^{-1})_{\zbv{\bar w}(t)}=0$.

Around $q$ there exist closed $1$-forms $\zm_{1},...,\zm_{r}$ such that
$Ker(\zm_{1}\zex...\zex\zm_{r})=\mathcal F$; therefore
$\zm_{1}\zci({\bar J}+tI)^{-1},...,\zm_{r}\zci({\bar J}+tI)^{-1}$ define
${\bar w}(t)$. But $d(\zm_{j}\zci({\bar J}+tI)^{-1})_{\zbv{\bar w}(t)}=0$,
$j=1,...,r$, so ${\bar w}(t)$ is involutive near $q$. On the other hand if $A$
is an open neighbourhood of $q$, small enough, there exists a non-empty open set
$B\zco \mathbb K$ such that ${\bar J}+tI$ is invertible on $A$, and therefore
${\bar w}(t)$ involutive, for any $t\zpe \mathbb K$. As ${\bar\zg}(t)$ is polynomial
in $t$  this implies that every ${\bar w}(t)$ is involutive on $A$ (indeed if $X,Y$
are vector fields belonging to $\mathcal F$ then
${\bar\zg}(t)[({\bar J}+tI)X,({\bar J}+tI)Y]$ is polynomial in $t$), {\it which
proves theorem 2.1}.
\bigskip

{\bf Corollary 2.1.1.} {\it Consider a Veronese web $w$ on $N$ of codimension
$1\zmei r\zmei n-1$, an immersion $f:P\zfl N$ and a scalar $b$.

\noindent (1) If for every   $p\zpe P$ the characteristic numbers of $w(f(p))$ are
greater than or equal to 2 and $f_{*}(T_{p}P)\zcco w(b)(f(p))$, then the family
$\{{\tilde w}(t)=f_{*}^{-1}(w(t))\zbv t\zpe{\mathbb K}-\{b\}\}$ extends to a
Veronese web $\tilde w$ on $P$ of codimension $r$ by setting
${\tilde w}(b)=lim\,{\tilde w}(t)$, $t\zfl b$.

\noindent (2) Now assume that the characteristic numbers of $w(f(p))$ are
constant on $P$; let ${\tilde r}$ the number of them greater than or equal to 2. If
$f_{*}(T_{p}P)=w(b)(f(p))$ for any $p\zpe P$, then the family
$\{{\tilde w}(t)=f_{*}^{-1}(w(t))\zbv t\zpe{\mathbb K}-\{b\}\}$ extends to a
Veronese web $\tilde w$ on $P$ of codimension $\tilde r$ by setting
${\tilde w}(b)=lim\,{\tilde w}(t)$, $t\zfl b$.}
\bigskip

{\bf Proof.} As the problem is local we may suppose that $P$ is a regular (imbedded)
submanifold of $N$ of codimension $k$ and $f$ the canonical inclusion. Consider non-equal
scalars $a_{1},...,a_{n-k},a$ where $a=-b$. Then in the construction of $J$ we can take
$\zb_{n-k+1},...,\zb_{n}$ in such a way that $Ker(\zb_{n-k+1}\zex...\zex\zb_{n})(p)=T_{p}P$,
$p\zpe P$; even more one may suppose $P=\{x\zbv x_{n-k+1}=...=x_{n}=0\}$ when
$\zb_{j}=dx_{j}$, $j=1,...,n$. On the other hand the integrability  is clear since
${\tilde w}(t)=w(t)\zin TP$, $t\zpe {\mathbb K}-\{b\}$.

Now consider a $r$-form $\zb$ such that $Ker \zb=w(\zinf)$ and
$\zb=\zm_{1}\zex...\zm_{r}$, where $\zm_{1},...,\zm_{r}$ are $1$-forms, and set
${\bar J}=\zsu_{j=1}^{n-k}a_{j}{\frac{\partial} {{\partial}x_{j}}}\zte dx_{j}$
on $P$. As $(\zm_{1},...,\zm_{r},J^{*})$ spans $T^{*}N$ then
$({\zm_{1}}_{\zbv P},...,{\zm_{r}}_{\zbv P},{\bar J}^{*})$ spans $T^{*}P$.

In the first case of the corollary ${\zb}_{\zbv P}$ has no zeros and
${\tilde \zg}(t)=(\zpr_{j=1}^{n-k}(t+a_{j}))(({\bar J}+tI)^{-1})^{*}({\zb}_{\zbv P})$
is a representative of $\tilde w$. In the second one
$\{{\zm_{1}(p)}_{\zbv P},...,{\zm_{r}(p)}_{\zbv P}\}$ spans a $\tilde r$-dimensional
vector subspace of $T_{p}^{*}P$ at any $p\zpe P$, which allows us to assume, for example,
that ${(\zm_{1}\zex...\zm_{\tilde r})}_{\zbv P}$ never vanishes (our problem is
local); then ${\tilde \zg}(t)=\zpr_{j=1}^{n-k}(t+a_{j})(({\bar J}+tI)^{-1})^{*}
{(\zm_{1}\zex...\zex\zm_{\tilde r})}_{\zbv P}$ is a representative of $\tilde w$.
$\square$

A family $w$ of $r$-codimensional distributions which satisfies all the conditions
of Veronese web except, perhaps, for the involutivity of each $w(t)$ will be called
a {\it Veronese distribution}.
\bigskip

{\bf Corollary 2.1.2.} {\it Consider a Veronese distribution $w$, of codimension
$r\zmai 1$, on $N$ and a point $p$ of this manifold. Let $n_1$ be the height of $w(p)$.
Assume that $w(\zinf),w(b_{1}),...,w(b_{n_{1}+1})$ are integrable for $n_{1}+1$
non-equal scalars $b_{1},...,b_{n_{1}+1}$. Then $w$ is a Veronese web around $p$.}
\bigskip

Indeed, let $k_{1}\zmai...\zmai k_{n_{1}}$ like in lemma 1.3.
Set ${\tilde a}_{j}=-b_{j}$, $j=1,...,n_{1}$. Since
$w(b_{1}),...,w(b_{n_{1}})$ are involutive, reasoning as in the construction of
$\tilde J$ gives rise to a $(1,1)$-tensor field $H$ defined around $p$, flat and
diagonalizable, with characteristic polynomial
$\zpr_{j=1}^{n_{1}}(t+b_{j})^{k_{j}}$ such that
$(Ker(H^{*}+b_{j}I))w(b_{j})=0$, $j=1,...,n_{1}$, and
$\zr(t)=\zpr_{j=1}^{n_{1}}(t-b_{j})^{k_{j}}((H+tI)^{-1})^{*}\zb$, where
$Ker\zb=w(\zinf)$, represents $w$.

On the other hand ${d(\zl\zci H)}_{\zbv w(\zinf)}=0$ for any $1$-closed form $\zl$
such that $Ker\zl\zcco w(\zinf)$, because $w(b_{n_{1}+1})$ is involutive and
$H+b_{n_{1}+1}I$ invertible (do reason as in the first and the second paragraphs
after the proof of lemma 2.1). Now apply (3) of the foregoing theorem.

{\bf Remark.} Corollary 2.1.2, with another proof, is due to Panasyuk \cite{PAB} and it
was conjectured by Zakharevich \cite{ZA} (see \cite{BO} by Bouetou and Dufour as well).
Note that by means of a projective transformation of
${\mathbb K}P^{1}\zeq{\mathbb K}\zun\{\zinf\}$, one may replace the integrability of
$w(\zinf)$ by that of $w(b_{n_{1}+2})$ for some
$b_{n_{1}+2}\zpe {\mathbb K}-\{b_{1},...,b_{n_{1}+1}\}$; in other words it suffices
the involutivity of $w(t)$ for $n_{1}+2$ elements of ${\mathbb K}P^{1}$. Therefore if
$k$ is the maximum of the height of $w(q)$, $q\zpe N$, a Veronese distribution $w$ is
a Veronese web on $N$ if and only if $w(t)$ is involutive for $k+2$ values of
$t\zpe {\mathbb K}P^{1}$.

By a similar reason, corollary 2.1.1 still holds if $b=\zinf$.
\bigskip

{\bf Proposition 2.1.} {\it Let $w=\{ w(t)\zbv t\zpe {\mathbb K}\}$ a family of foliations
of codimension $r$ defined on a manifold $N$. Assume that each $w(p)$ is an algebraic Veronese
web on $T_{p}N$, which allows us to define a $r$-codimensional distribution ${\mathcal F}$ on $N$,
possibly not smooth, by setting ${\mathcal F}(p)=limw(t)(p)$, $t\zfl \zinf$. If ${\mathcal F}$
is smooth then $w$ is a Veronese web on $N$.}
\bigskip

{\bf Proof.} Note that the $(1,1)$ tensor field $J$ may be constructed, as before, around each point
of $N$ since every $w(t)$ is a foliation. On the other hand locally there exists a $r$-form $\zb$
such that $Ker\zb={\mathcal F}$ because ${\mathcal F}$ is smooth. So
$\zg(t)=(\zpr_{j=1}^{n-k}(t+a_{j}))(t+a)^{k}((J+tI)^{-1})^{*}\zb$ is a
representative of $w$. $\square$

{\bf Example.} On an open set $A$ of ${\mathbb K}^n$ consider a $(1,1)$-tensor field
$J= \zsu_{j=1}^{n}f_{j}(x_{j}){\frac{\partial} {{\partial}x_{j}}}\zte dx_{j}$ where
$f_{j}(x_{j})\znoi f_{k}(x_{k})$ whenever $x=(x_{1},...,x_{n})\zpe A$. Set
$\zb=\zsu_{j=1}^{n}dx_{j}$. As $N_{J}=0$, $(\zb, J^{*})$ spans $T^{*}A$ and
$d(\zb\zci J)=0$, by (3) of theorem 2.1 the curve
$\zg(t)= \zpr_{j=1}^{n}(t+f_{j})\zb\zci(J+tI)^{-1}
=\zsu_{j=1}^{n}(\zpr_{i=1;i\znoi j}^{n}(t+f_{i}))dx_{j}$ defines a Veronese web $w$ on
$A$ of codimension one, which generally is not flat.

Indeed, when $w$ is flat there exists a representative
${\tilde \zg}(t)=\zsu_{i=0}^{n-1}t^{i}{\tilde \zg}_{i}$ with each
${\tilde \zg}_{i}$ closed. Set $\zg(t)=\zsu_{i=0}^{n-1}t^{i}\zg_{i}$. As
$\zg=f\tilde\zg$ then $\zg_{i}\zex d\zg_{i}=0$, $i=0,...,n-1$. But
$\zg_{n-2}=\zsu_{j=1}^{n}(f_{1}+...+{\hat f}_{j}+...+f_{n})dx_{j}$; so the
coefficient of $dx_{i}\zex dx_{j}\zex dx_{k}$, where $i<j<k$, in the
expression of $\zg_{n-2}\zex d\zg_{n-2}=0$ equals
$f_{i}(f'_{k}-f'_{j})+f_{j}(f'_{i}-f'_{k})+f_{k}(f'_{j}-f'_{i})$ which almost never
vanishes.

For obtaining a $2$-codimensional Veronese web $\tilde w$, one may consider a second
$1$-form $\zb'=\zsu_{j=1}^{n}g_{j}(x_{j})dx_{j}$ such that $\zb\zex\zb'$ never vanishes
and set ${\tilde\zg}(t)=\zpr_{j=1}^{n}(t+f_{j})((J+tI)^{-1})^{*}(\zb\zex\zb')=
\zsu_{1\zmei j<k\zmei n}(\zpr_{i=1;i\znoi j,k}^{n}(t+f_{i}))(g_{k}-g_{j})dx_{j}\zex dx_{k}$.

Theorem 2.1 gives a method to construct all Veronese webs locally. Usually the scalars
$a_{1},...,a_{n-k},a$, respectively ${\tilde a}_{1},...,{\tilde a}_{n_{1}}$, do not determine
$J$, respectively $\tilde J$, which prevent us constructing them globally. Nevertheless if the
characteristic numbers are constant and equal, for example if $r=1$, then
$n_{1}=\frac n {r}$, $k_{1}=...=k_{n_{1}}=r$ and $\tilde J$ can be constructed
on all $N$ since, now, $Ker ({\tilde J}^{*}-{\tilde a}_{j}I)$ is the annihilator
of $w(-{\tilde a}_{j})$.

On the other hand, in view of proposition 1.2, the restriction of $J$ or $\tilde J$ to
$w(\zinf)$ gives rise to a morphism (of vector bundles) $\zlma:w(\zinf)\zfl TN$, which only
depends on the Veronese web,
without non-zero  $\zlma$-invariant vector subspace at any
point of $N$. Moreover $w(t)=(\zlma+tI)w(\zinf)$, $t\zpe \mathbb K$.

Now consider, on a manifold $M$, a foliation $\mathcal F$ and a morphism (of vector bundles)
$G:{\mathcal F}\zfl TM$. If $\za$ is a $s$-form defined on an open set $A$ of $M$, then
$G^{*}\za$ is a section on $A$ of $\zL^{s}{\mathcal F}^{*}$ and can regarded as a $s$-form on the
leaves of $\mathcal F$; thus we shall say that is {\it closed on $\mathcal F$} if it is closed
on its leaves. Besides, when ${\bar G}:TM\zfl TM$ is a prolongation of $G$, then
$d({\bar G}^{*}\za)_{\zbv{\mathcal F}}$ equals the exterior derivative of $G^{*}\za$
along the leaves of ${\mathcal F}$; thus $G^{*}\za$ is closed on $\mathcal F$ if and only if
$d({\bar G}^{*}\za)_{\zbv{\mathcal F}}=0$.
\bigskip

{\bf Lemma 2.2.} {\it Assume that $G^{*}\za$ is closed on $\mathcal F$ for every closed
$1$-form $\za$ such that $Ker\za\zcco\mathcal F$. Then the restriction of $N_{\bar G}$ to
$\mathcal F$, which will be named {\rm the Nijenhuis torsion of $G$} and denoted by $N_G$,
does not depend on the prolongation $\bar G$.}
\bigskip

{\bf Proof.} As the problem is local we may suppose that
${\mathcal F}=Ker(\za_{1}\zex...\zex\za_{k})$ where each $\za_j$ is a closed $1$-form and
$k=codim{\mathcal F}$. Since the difference between two prolongations equals
$\zsu_{j=1}^{k}Y_{j}\zte\za_{j}$, it suffices to consider the case
$H={\bar G}+Y\zte\za$ with $\za\zex\za_{1}\zex...\zex\za_{k}=0$ and $d\za=0$.
Now given $X\zpe\mathcal F$ one has:

\noindent $N_{H}(X,\quad)=L_{HX}H-HL_{X}H=L_{{\bar G}X}({\bar G}+Y\zte\za)-
{\bar G}L_{X}({\bar G}+Y\zte\za)-Y\zte\za(L_{X}{\bar G}+L_{X}(Y\zte\za))$

\noindent whence  $N_{H}(X,\quad)-N_{\bar G}(X,\quad)=Y\zte
( L_{{\bar G}X}\za-\za(L_{X}{\bar G}))+{\tilde Y}\zte\za$ because $L_{X}\za=d(\za(X))=0$.

On the other hand when $ Z\zpe \mathcal F$:

\noindent $( L_{{\bar G}X}\za-\za(L_{X}{\bar G}))(Z)=Z\za({\bar G}X)-\za([X,{\bar G}Z])
+\za({\bar G}[X,Z])=Z\za({\bar G}X)-X\za({\bar G}Z)+\za({\bar G}[X,Z])
=-d(\za\zci{\bar G})(X,Z)=0$

\noindent since $\za$ is closed and
$\za\zci{\bar G}$ is closed on $\mathcal F$. Therefore
$(N_{H})_{\zbv \mathcal F}=(N_{\bar G})_{\zbv \mathcal F}$.
$\square$

Note that the Nijenhuis torsion of $\zlma:w(\zinf)\zfl TN$ vanishes and $\zlma^{*}\za$ is
closed on $w(\zinf)$ for every closed $1$-form $\za$ such that $Ker\za\zcco w(\zinf)$ since $J$,
its local prolongation given by (1) of theorem 2.1, has zero Nijenhuis torsion and
$d(\za\zci J)_{\zbv w(\zinf)}=0$.

Conversely, given a foliation $\mathcal F$ on $N$ of codimension $1\zmei r\zmei n$ and a
morphism $\zlma:{\mathcal F}\zfl TN$ with the algebraic and differentiable properties stated
before, then $w(t)=(\zlma+tI){\mathcal F}$, $t\zpe \mathbb K$, defines a Veronese of
codimension $r$ for which $w(\zinf)={\mathcal F}$. Indeed apply (3) of theorem 2.1 to a
prolongation $\bar J$ of $\zlma$. Thus:

{\it Giving a Veronese web on $N$ of codimension $r\zmai 1$ is equivalent
to giving a morphism $\zlma:{\mathcal F}\zfl TN$, where ${\mathcal F}$ is a
$r$-codimensional foliation without non-vanishing $\zlma$-invariant vector subspace at
any point such that:

\noindent 1) whenever $\za$ is a closed $1$-form whose kernel contains $\mathcal F$,
restricted to the domain of $\za$, then $\zlma^{*}\za$ is closed on $\mathcal F$,

\noindent 2) $N_{\zlma}=0$.}

Note that if ${\mathcal F}=Ker(\za_{1}\zex...\zex\za_{r})$, where $d\za_{1}=...=d\za_{r}=0$, then
$\zlma^{*}\za$ is closed on $\mathcal F$ for any $1$-form $\za$ such that $d\za=0$ and
$Ker\za\zcco\mathcal F$, if and only if
$\zlma^{*}\za_{1},...,\zlma^{*}\za_{r}$ are closed on $\mathcal F$.

{\bf Example.} On an open set $A$ of ${\mathbb K}^{2m}$, endowed with coordinates
$(x,y)=(x_{1},...,x_{m},y_{1},...,y_{m})$, consider the foliation ${\mathcal F}$ defined by
$dy_{1}=...=dy_{m}=0$ and the morphism $\zlma:{\mathcal F}\zfl TA$ given by
$\zlma({\frac {\partial} {{\partial}x_{j}}})=
\zsu_{k=1}^{m}f_{jk}{\frac {\partial} {{\partial}y_{k}}}$, $j=1,...m$. Assume
$\zbv f_{jk}\zbv\znoi 0$ everywhere, which implies that $\zlma:{\mathcal F}\zfl TA$
defines a $m$-codimensional Veronese distribution $w$ on $A$ with characteristic numbers
$n_{1}=...=n_{m}=2$. Then $w$ is a Veronese web if and only if
$d(\zsu_{j=1}^{m}f_{jk}dx_{j})_{\zbv \mathcal F}=0$, $k=1,...,m$, and
$\zcizq\zsu_{k=1}^{m}f_{jk}{\frac {\partial} {{\partial}y_{k}}},
\zsu_{{\tilde k}=1}^{m}f_{{\tilde j}{\tilde k}}{\frac {\partial}
{{\partial}y_{\tilde k}}}\zcder=0$, $1\zmei j<{\tilde j}\zmei m$
(indeed consider the prolongation $J$ of $\zlma$ given by
$J({\frac {\partial} {{\partial}y_{k}}})=0$, $k=1,...,m$).

When $m=1$ there are no conditions at all. If $m=2$ one has a partial differential
system of order one with four equations and four functions; for $m\zmai 3$ the system is
over-determined.

More generally when $n=2m$, the $m$-dimensional Veronese webs on $N$, with characteristic
numbers $n_{1}=...=n_{m}=2$, are given by a morphism $\zlma:{\mathcal F}\zfl TN$ such that
$dim{\mathcal F}=m$ and $TN={\mathcal F}\zdi Im\zlma$. As $\zlma$ is determined by its image
and its graph, which may be identified to $w(1)=(\zlma+I){\mathcal F}$, from the algebraic
viewpoint giving a Veronese web $w$ with all its characteristic number equal to 2 is like giving
the $3$-web $\{{\mathcal F}=w(\zinf),w(0),w(1)\}$. Conversely, for any $3$-web
$D=\{{\mathcal D}_{1},{\mathcal D}_{2},{\mathcal D}_{3}\}$ on $N$ there exists just one
Veronese distribution $w_D$ such that $w_{D}(\zinf)={\mathcal D}_{1}$,
$w_{D}(0)={\mathcal D}_{2}$ and $w_{D}(1)={\mathcal D}_{3}$. It is easily seen that $w_D$ is
a Veronese web if and only if the torsion of the Chern connection of $D$ vanishes (the Chern connection
of $D$ is the only connection making ${\mathcal D}_{1},{\mathcal D}_{2},{\mathcal D}_{3}$ parallel
such that $T({\mathcal D}_{1},{\mathcal D}_{2})=0$, see \cite{NAG}).

For the link between $k$-webs, $k\zmai 4$, and Veronese webs see \cite{BO} (Bouetou-Dufour).
\bigskip

{\bf 3. Kronecker bihamiltonian structures}

Consider two Poisson structures $\zL,\zL_1$ defined on a real or complex manifold $M$
of dimension $m$. Following Magri \cite{MAG} we will say that $(\zL,\zL_{1})$ is a
{\it bihamiltonian structure} (or that $\zL,\zL_1$ are {\it compatible}) if
$\zL+\zL_1$ is still a Poisson structure, which is equivalent to say that their
Schouten bracket vanishes or that $\zL+b\zL_1$ is a Poisson structure
for some $b\zpe{\mathbb K}-\{0\}$. Recall
that if $\zL,\zL_1$ are compatible then $a\zL+a_{1}\zL_1$ is a Poisson structure for
all $a,a_{1}\zpe{\mathbb K}$.

A bihamiltonian structure $(\zL,\zL_{1})$ will be called {\it Kronecker} when there exists
$r\zpe{\mathbb N}-\{0\}$ such that each $(\zL(p),\zL_{1}(p))$, $p\zpe M$, is the product of
$r$ Kronecker elementary pairs. In this case from the algebraic model at each point
follows that
$m-r=rank(\zL,\zL_{1})=rank(\zL)=rank(\zL_{1})=rank(\zL+t\zL_{1})$ for any $t\zpe{\mathbb K}$;
moreover ${\mathcal D}=\zin Im(\zL+t\zL_{1})$, $t\zpe{\mathbb K}$, is a foliation of dimension
${\frac{m-r} 2}$ lagrangian for both $\zL$ and $\zL_1$, and ${\mathcal D}\zco Im\zL_{1}$.
This foliation will be named the {\it soul of} $(\zL,\zL_{1})$.

Let $N$ be the local quotient of $M$ by the foliation ${\mathcal D}$, which is a manifold
of dimension $n={\frac{m+r} 2}$, and let $\zp:M\zfl N$ be the canonical projection. Then
$w=\{w(t)=\zp_{*}(Im(\zL+t\zL_{1}))\zbv t\zpe\mathbb K\}$ is a family of foliation on $N$ of
codimension $r$, whose limit when $t\zfl\zinf$ equals $\zp_{*}(Im\zL_{1})$ since
$\zp_{*}(Im(\zL+t\zL_{1}))=\zp_{*}(Im(s\zL+\zL_{1}))$ where $s=t^{-1}$. Besides {\it $w$ is a
Veronese web of codimension $r$}.

Indeed, given $p\zpe N$ such that $\zp(q)=p$, proposition 1.4 applied to
$(\zL(q),\zL_{1}(q),T_{q}M)$ shows that $w(p)$ is an algebraic Veronese web. Now apply
proposition 2.1.

In short a Veronese web of codimension r is locally associated to
any Kronecker bihamiltonian structure with $r$ factors. Our next
goal is to study when this Veronese web locally determines the
Kronecker bihamiltonian structure.

Recall that a Poisson structure $\zL'$ on $M$ of constant rank
$m-r$ can be locally described by $r$ closed $1$-forms giving the
foliation $Im\zL'$ and a $2$-form whose restriction to $Im\zL'$ is
symplectic; this last one is only defined modulo the ideal spanned
by the $1$-forms. Consider non-equal and non-vanishing scalars
$a_{1},...,a_{n-r},a$, any point $p\zpe N$ and closed $1$-forms
$\za_{1},...,\za_{r}$, defined around $p$, such that
$Ker(\za_{1}\zex...\zex\za_{r})=w(\zinf)$. Let $J$ be a $(1,1)$
tensor field like in part (1) of theorem 2.1; then
$(\za_{1},...,\za_{r}, J^{*})$ spans the cotangent bundle near $p$
and $\za_{1}\zex...\zex\za_{r}\zex d(\za_{j}\zci J)=0$,
$j=1,...,r$. On the other hand one may choose coordinates
$(x_{1},...,x_{n-r},y_{1},...,y_{r})$, defined on an open
neighbourhood of $p\zeq 0$, such that $dx_{j}\zci J=a_{j}dx_{j}$,
$j=1,...,n-r$, and $w(0)=Ker(dy_{1}\zex...\zex dy_{r})$; indeed
the choice of $x_{1},...,x_{n-r}$ is obvious and
$dx_{1},...,dx_{n-r}$ restricted to $w(0)$ are linearly
independent everywhere since they are independent restricted to
$w(-a)$ and $w(-a)=(J-aI)J^{-1}w(0)$. As $\mathcal D$ is
$\zL$-lagrangian functions $x_{1}\zci\zp,...,x_{r}\zci\zp$ are in
$\zL$-involution, so around each $p'\zpe \zp^{-1}(p)$ there exist
functions $f_{1},...,f_{n-r}$, vanishing at $p'$, such that $\zL$
is given by $d(y_{1}\zci\zp),...,d(y_{r}\zci\zp)$ and
$d(x_{1}\zci\zp)\zex df_{1}+...+d(x_{n-r}\zci\zp)\zex df_{n-r}$.
Now by setting $z_{j}=f_{j}$ and writing $x_j$ and $y_k$ instead
of $x_{j}\zci\zp$ and $y_{k}\zci\zp$, for sake of simplicity, we
construct a system of coordinates
$(x,y,z)=(x_{1},...,x_{n-r},y_{1},...,y_{r},z_{1},...,z_{n-r})$
such that $p' \zeq 0$, $\zp(x,y,z)=(x,y)$ and $\zL$ is given by
$dy_{1},...,dy_{r}$, $\zsu_{j=1}^{n-r}dx_{j}\zex dz_{j}$.

But $\mathcal D$ is $\zL_1$-lagrangian too, so $x_{1},...,x_{n-r}$
are in $\zL_1$-involution. Moreover on $N$ forms $dx_{1},...,dx_{n-r}$
restricted to $w(\zinf)$ are linearly independent everywhere since
$w(-a)=(J-aI)w(\zinf)$; therefore around $p'$ there exist functions
$g_{1},...,g_{n-r}$ such that $\zL_1$ is given by $dx_{1}\zex
dg_{1}+...+dx_{n-r}\zex dg_{n-r}$ and $\za_{1},...,\za_{r}$ (more
exactly $\zp^{*}\za_{1},...,\zp^{*}\za_{r}$). On the other hand
$\zp_{*}^{-1}(w(-a_{j}))=Im(\zL-a_{j}\zL_{1})\zco Kerdx_{j}$ whence
$(\zpar/\zpar z_{j})
=\zL(dx_{j},\quad)=a_{j}\zL_{1}(dx_{j},\quad)$
and $(\zpar g_{k}/\zpar z_{j})
=\zd_{jk}a_{k}$.
So $\zL_1$ is given
by $\za_{1},...,\za_{r}$ and $\zsu_{j=1}^{n-r}a_{j}dx_{j}\zex
dz_{j}+\zw$ where

\centerline{$\zw=\zsu h_{ij}(x,y)dx_{i}\zex dx_{j}+\zsu
{\tilde h}_{ik}(x,y)dx_{i}\zex dy_{k}$ and $d\zw=0$.}

Thus $\zw$ may be regarded as a closed $2$-form on an open
neighbourhood of $p$ in $N$.

Given a $k$-form $\zt$, $k\zmai 1$, and a $(1,1)$ tensor field $H$
on a manifold, $\zt\zci H$ and $\zt_{H}$ will denote the $k$-forms
defined by $(\zt\zci H)(X_{1},...,X_{k})=\zt(HX_{1},...,HX_{k})$
and $\zt_{H}(X_{1},...,X_{k})=\zt(HX_{1},X_{2}...,X_{k})
+\zt(X_{1},HX_{2},...,X_{k})+...+\zt(X_{1},...,X_{k-1},HX_{k})$
respectively.

The next proposition, proved later on, characterizes the
compatibility of $\zL$ and $\zL_1$.
\bigskip

{\bf Proposition 3.1.} {\it The pair $(\zL,\zL_{1})$ is compatible
if and only if $\za_{1}\zex...\zex\za_{r}\zex d\zw_{J}=0$.}
\bigskip

The local determination of the bihamiltonian structure by the
Veronese web will be established if we are able to delete the term
$\zw$ in the expression of $\zL_1$, since $\za_{1},...,\za_{r}$
only depend on the web. Given a function $\zf(x,y)$ defined around
$p$ set $u_{j}=z_{j}-({\zpar \zf}/{\zpar x_{j}})$,
$j=1,...,n-r$. Then, in coordinates $(x,y,u)$,
$dy_{1},...,dy_{r}$, $\zsu_{j=1}^{n-r}dx_{j}\zex du_{j}$ define
$\zL$ (the other terms belong to the ideal spanned by
$dy_{1},...,dy_{r}$) while $\zL_1$ is given by
$\za_{1},...,\za_{r}$, $\zsu_{j=1}^{n-r}a_{j}dx_{j}\zex
du_{j}+(\zw-d(d\zf\zci J))$;
indeed each $(dy_{k}\zci J)\zex\za_{1}\zex...\zex\za_{r}=0$ since
$Jw(\zinf)=w(0)$, so
$(d\zf\zci J-\zsu_{j=1}^{n}a_{j}(\zpar \zf /\zpar x_{j})dx_{j})
\zex\za_{1}\zex...\zex\za_{r}=0$.
As the $2$-form expressing $\zL_1$ is
defined modulo the ideal spanned by $\za_{1},...,\za_{r}$, it
suffices to find a function $\zf$ such that
$\za_{1}\zex...\zex\za_{r}\zex d(d\zf\zci
J)=\za_{1}\zex...\zex\za_{r}\zex\zw$ for deleting $\zw$. To remark
that if a such function $\zf$ exists, by adding a suitable linear
function of $(x,y)$ we may suppose $d\zf(p)=0$ and $u_{j}(p')=0$,
$j=1,...,n-r$.
\bigskip

{\bf Theorem 3.1.} {\it On a manifold $N$ consider closed
$1$-forms $\za_{1},...,\za_{r}$, $r\zmai 1$, linearly independent
everywhere and a $(1,1)$ tensor field $J$, which is flat and
diagonalizable with characteristic polynomial
$(t-a)^{r}\zpr_{j=1}^{n-r}(t-a_{j})$ where $a_{1},...,a_{r},a$ are
non-equal scalars. Assume that $(\za_{1},...,\za_{r},J^{*})$ spans
$T^{*}N$ and $\za_{1}\zex...\zex\za_{r}\zex d(\za_{j}\zci J)=0$,
$j=1,...,r$.

Given a closed $2$-form $\zw$ on $N$ if $d\zw_{J}=0$ then, around
each point of $N$, there exists a function $\zf$ such that
$\za_{1}\zex...\zex\za_{r}\zex d(d\zf\zci
J)=\za_{1}\zex...\zex\za_{r}\zex\zw$  at least in the following
three cases:

\noindent (1) on complex manifold,

\noindent (2) in the real analytic category,

\noindent (3) in the $C^{\zinf}$ category when $r=1$.}
\bigskip

This theorem will be proved in the next section.
\bigskip

{\bf Theorem 3.2.} {\it From the local viewpoint the Veronese web
completely determines the Kronecker bihamiltonian structure, at
least, in the following four cases: complex manifold, real
analytic category, $C^{\zinf}$ category when $r=1$, and flat
Veronese web.}
\bigskip

Theorem 3.2 is an obvious consequence of theorem 3.1 except for real
flat webs. In this last case in some coordinates
$(v_{1},...,v_{n})$ the expression of $w(t)$ does not depend on
the point considered, which allows us to choose
$\za_{1},...,\za_{r}$ and $J$ with constant coefficients. Thus in
these coordinates the partial differential equation
$\za_{1}\zex...\zex\za_{r}\zex d(d\zf\zci
J)=\za_{1}\zex...\zex\za_{r}\zex\zw$ is homogeneous of of order
two with constant coefficients and $C^{\zinf}$ independent term.
By the Ehrenpreis-Malgrange theorem (see \cite{MAL}) there exist local
solutions provided that it has formal solutions.

Let $\zw_k$ be the $k$th term of the Taylor expansion of $\zw$,
always in coordinates $(v_{1},...,v_{r})$, at point $q$. Then
$d\zw_{k}=0$ and $\za_{1}\zex...\zex\za_{r}\zex d({(\zw_{k})}_{J})=0$,
so by theorem 3.1 the equation
$\za_{1}\zex...\zex\za_{r}\zex d(d\zf\zci
J)=\za_{1}\zex...\zex\za_{r}\zex\zw_{k}$ has a solution ${\tilde
f}$ around $q$. Note that the $(k+2)$th term $f_{k+2}$ of the Taylor
expansion of ${\tilde f}$ at $q$ is a solution of this equation
too. Thus if $f$ is a polynomial of degree $\zlma\zmai 2$ such
that $\za_{1}\zex...\zex\za_{r}\zex d(df\zci
J)=\za_{1}\zex...\zex\za_{r}\zex(\zw_{0}+...+\zw_{\zlma-2})$ then
$\za_{1}\zex...\zex\za_{r}\zex d(d(f+f_{\zlma+1})\zci
J)=\za_{1}\zex...\zex\za_{r}\zex(\zw_{0}+...+\zw_{\zlma-1})$.
Therefore the equation $\za_{1}\zex...\zex\za_{r}\zex d(d\zf\zci
J)=\za_{1}\zex...\zex\za_{r}\zex\zw$ is formally integrable and
there exist local solutions of it around each point.

Theorem 3.2 was proved by Gelfand and Zakharevich \cite{GEA,GEB} for
analytic Veronese web of codimension 1; the flat case, the $C^{\zinf}$
case of codimension 1 and the analytic one of any codimension are
due to Turiel \cite{TUA,TUC}.

Now we will prove proposition 3.1
\bigskip

{\bf Lemma 3.1.}{ \it If $t\znope\{ -a_{1},...,-a_{n-r},-a\}$ then
$\zL+t\zL_{1}$ is defined by
$\za_{1}\zci(J+tI)^{-1},...,\za_{r}\zci(J+tI)^{-1}$ and
$\zsu_{j=1}^{n-r}a_{j}(t+a_{j})^{-1}dx_{j}\zex dz_{j} +t\zw\zci
(J+tI)^{-1}$.}
\bigskip

{\bf Proof.} First we replace coordinates $(y_{1},...,y_{r})$ by
coordinates $(u_{1},...,u_{r})$ such that $du_{k}\zci J=adu_{k}$,
thus $J=\zsu_{j=1}^{n-r}a_{j}{\frac {\zpar }{\zpar x_{j}}\zte
dx_{j}+\zsu_{k=1}^{r}a\frac {\zpar }{\zpar u_{k}}}\zte du_{k}$
in coordinates $(x_{1},...,x_{n-r},u_{1},...,u_{r})$.
Let $V$ be a $r$-dimensional vector space and let
$\{e_{1},...,e_{r}\}$ be a basis of $V$. It will be enough to
prove the result for each  point $q$. On $T_{q}M\zdi V$ set
$\zW=\zsu_{j=1}^{n-r}dx_{j}\zex dz_{j}+\zsu_{k=1}^{r}du_{k}\zex
e_{k}^{*}$, $\zW_{1}=\zsu_{j=1}^{n-r}a_{j}dx_{j}\zex
dz_{j}+\zsu_{k=1}^{r}adu_{k}\zex e_{k}^{*}+\zw$ where
$dx_{j},dz_{j},dy_{k},du_{k},e_{k}^{*}$ and $\zw$ are extended to
$T_{q}M\zdi V$ in the obvious way and the point $q$ is omitted in
the notation.

Let $G$ and $H$ be the endomorphisms of $T_{q}M\zdi V$ defined by
$\zW(G,\quad)=\zW_{1}-\zw$ and $\zW(H,\quad)=\zw$ respectively.
Note that $G=\zsu_{j=1}^{n-r}a_{j}({\frac {\zpar } {\zpar
x_{j}}}\zte dx_{j} +{\frac {\zpar } {\zpar z_{j}}}\zte dz_{j})
+\zsu_{k=1}^{r}a({\frac {\zpar } {\zpar u_{k}}}\zte du_{k}
+e_{k}\zte e_{k}^{*})$, $dx_{j}\zci G=dx_{j}\zci J$, $dy_{k}\zci
G=dy_{k}\zci J$, $du_{k}\zci G=du_{k}\zci J$, $j=1,...,n-r$,
$k=1,...,r$, and $ImH\zco U\zco KerH$, so $H^2=0$, where $U$ is
the vector space spanned by
$(\zpar/\zpar z_{1}),...,(\zpar/\zpar z_{n-r}),e_{1},...,e_{r}$.

Let $W$ be the $r$-dimensional vector subspace of $T_{q}M\zdi V$ whose
image by $\zW$ is the space spanned by $dy_{1},...,dy_{r}$
(note that this last space is
the annihilator of $w(0)\zdi U$). Obviously $W\zco U$ so
$W$ is $\zW$-isotropic; moreover $W$ is a direct factor of $T_{q}M$
since $dx_{1},...,dx_{n-r},dy_{1},...,dy_{r}$ are linearly
independent. On the other hand $\zW_{1}(W,\quad)=\zW(GW,\quad)$ is
spanned by $dy_{1}\zci G=dy_{1}\zci J,...,dy_{r}\zci G=dy_{r}\zci
J$. As $Jw(\zinf)=w(0)$, $\zW_{1}(W,\quad)$ is spanned by
$\za_{1},...,\za_{r}$ too; that is to say $\zW_{1}(W,\quad)$ is
the annihilator of $w(\zinf)\zdi U$ and $W$ is $\zL_{1}$-isotropic too.

By lemma 1.4 bivectors $\zL,\zL_{1}$ are the projection on ${\frac
{T_{q}M\zdi V} {W}}\zeq T_{q}M$ of the dual bivectors $\zL_{\zW}$
and $\zL_{\zW_{1}}$. Therefore $\zL+t\zL_{1}$ is the projection of
$\zL_{\zW}+t\zL_{\zW_{1}}$, which is the dual bivector of
$\zW((I+t(G+H)^{-1})^{-1},\quad)$.

By lemma 1.5 the space $W$ is isotropic for this last symplectic form, so
$\zL+t\zL_{1}$ will be given by the restriction to $T_{q}M$ of
$\zW((I+t(G+H)^{-1})^{-1},\quad)$ and $\zW((I+t(G+H)^{-1})^{-1}W,\quad)$.

Recall that if $A$ is an automorphism and $B$ an endomorphism such
that $B^{2}=0$ and $A^{-1}(ImB)\zco Ker B$, then
$(A+B)^{-1}=A^{-1}-A^{-1}BA^{-1}$. So
$(G+H)^{-1}=G^{-1}-G^{-1}HG^{-1}$ and
$(I+t(G+H)^{-1})^{-1}=((I+tG^{-1})-tG^{-1}HG^{-1})^{-1}
=(I+tG^{-1})^{-1}+t(G+tI)^{-1}H(G+tI)^{-1}$.

Hence $\zW((I+t(G+H)^{-1})^{-1},\quad)
=\zsu_{j=1}^{n-r}a_{j}(t+a_{j})^{-1}dx_{j}\zex dz_{j}
+\zsu_{k=1}^{r}a(t+a)^{-1}du_{k}\zex e^{*}_{k}+t\zw\zci (J+tI)^{-1}$ and
$\zW((I+t(G+H)^{-1})^{-1}W,\quad)=\zW((I+tG^{-1})^{-1}W,\quad)$ equals the
vector space spanned by $dy_{1}\zci(I+tG^{-1})^{-1}$,...,
$dy_{r}\zci(I+tG^{-1})^{-1}$, that is to say by
$\za_{1}\zci(J+tI)^{-1}$,..., $\za_{1}\zci(J+tI)^{-1}$, since
$dy_{k}\zci(I+tG^{-1})^{-1}=dy_{k}\zci(I+tJ^{-1})^{-1}
=(dy_{k}\zci J)\zci(J+tI)^{-1}$ and $Jw(\zinf)=w(0)$. $\square$
\bigskip

{\bf Lemma 3.2.} {\it Consider a $k$-form $\zt$, $k\zmai 1$, and a (1, 1) tensor
field $G$ on a manifold. Suppose that the Nijenhuis torsion of $G$ vanishes. Then
$(d(\zt\zci G))_{G} = d((\zt\zci G)_{G}) + (d\zt)\zci G$.}
\bigskip

{\bf Proof.} By induction on $k$. The case $k=1$ follows from
lemma 2.1; on the other hand if $k\zmai 2$ it suffices proving the
lemma when $\zb=\zb_{1}\zex\zb_{2}$ and $\zb_1$ is a $1$-form.
Then
\medskip

\noindent $(d(\zb\zci G))_{G}=(d(\zb_{1}\zci G)\zex(\zb_{2}\zci
G))_{G} -((\zb_{1}\zci G)\zex d(\zb_{2}\zci G))_{G}= (d(\zb_{1}\zci
G))_{G}\zex(\zb_{2}\zci G) +d(\zb_{1}\zci G)\zex(\zb_{2}\zci
G)_{G}- (\zb_{1}\zci G)_{G}\zex d(\zb_{2}\zci G)- (\zb_{1}\zci
G)\zex (d(\zb_{2}\zci G))_{G}$
\medskip

\noindent $d((\zb\zci G)_{G})=d((\zb_{1}\zci
G)_{G}\zex(\zb_{2}\zci G))+ d((\zb_{1}\zci G)\zex(\zb_{2}\zci
G)_{G})= d((\zb_{1}\zci G)_{G})\zex(\zb_{2}\zci G)+ d((\zb_{1}\zci
G))\zex(\zb_{2}\zci G)_{G}- (\zb_{1}\zci G)_{G}\zex d(\zb_{2}\zci
G)- (\zb_{1}\zci G)\zex d((\zb_{2}\zci G)_{G})$
\medskip

\noindent $(d\zb)\zci G=((d\zb_{1})\zci G)\zex(\zb_{2}\zci G)-
(\zb_{1}\zci G)\zex (d\zb_{2})\zci G$.
\medskip

Now take into account that the formula is true for $\zb_1$ (lemma
2.1) and $\zb_2$ (induction hypothesis), and remark that the
second and third terms of the expansion of $(d(\zb\zci G))_{G}$
equal the second and third ones of $d((\zb\zci G)_{G})$. $\square$
\bigskip

By lemma 3.1, $\,\zL$ and $\zL_1$ are compatible if and only if
$(\za_{1}\zci(J+tI)^{-1})\zex...\zex(\za_{r}\zci(J+tI)^{-1})\zex
d(\zw\zci(J+tI)^{-1})=0$ for some
$t\znope\{-a_{1},...,-a_{n-r},-a\}$, that is to say when
$\za_{1}\zex...\zex\za_{r}\zex
(d(\zw\zci(J+tI)^{-1})\zci(J+tI))=0$. Lemma 3.2 applied to
$\zw\zci(J+tI)^{-1}$ and $J+tI$ yields
$d(\zw\zci(J+tI)^{-1})\zci(J+tI)=-d(\zw_{(J+tI)})= -d\zw_{J}$.
Therefore $\zL,\zL_1$ are compatible if and only if
$\za_{1}\zex...\zex\za_{r}\zex d\zw_{J}=0$, {\it which proves
proposition 3.1}.
\bigskip

Consider a foliation ${\mathcal F}$ of codimension $s$ defined on a
$k$-manifold $P$. Let ${\mathcal F}'$ be the foliation, on the
cotangent bundle $T^{*}{\mathcal F}$ of the first foliation,
pull-back of ${\mathcal F}$ by the canonical projection
$\zp:T^{*}{\mathcal F}\zfl P$; that is to say ${\mathcal F}'(\zb)=
(\zp_{*}(\zb)^{-1})({\mathcal F}(\zp(\zb)))$
(until the end of this section one will write $T^{*}{\mathcal F}$ instead
of ${\mathcal F}^*$ for pointing out that $T^{*}{\mathcal F}$ is regarded
as a manifold itself).
On the leaves of
${\mathcal F}'$ one defines the Liouville  $1$-form $\zr$ by
setting $\zr(\zb)(X)=\zb(\zp_{*}(X))$ for any $X\zpe{\mathcal
F}'(\zb)\zco T_{\zb}(T^{*}{\mathcal F})$ and any $\zb\zpe
T^{*}{\mathcal F}$, and the Liouville $2$-form ${\tilde{\zw}}
=-d\zr$; then ${\tilde{\zw}}$ is symplectic on the leaves of
${\mathcal F}'$ and, by duality, gives rise to a Poisson structure
$\zL_L$ such that $Im\zL_{L}={\mathcal F}'$, which  will be named
the {\it Liouville-Poisson structure of $T^{*}{\mathcal F}$}. In
coordinates $({\tilde x},{\tilde y})= ({\tilde x}_{1},...,{\tilde
x}_{k},{\tilde y}_{1},...,{\tilde y}_{k-s})$, associated to
coordinates ${\tilde x}= ({\tilde x}_{1},...,{\tilde x}_{k})$ on
$P$ such that ${\mathcal F}$ were defined by $d{\tilde
x}_{k-s+1}=...=d{\tilde x}_{k}=0$, $\zL_L$ is given by $d{\tilde
x}_{k-s+1},...,d{\tilde x}_{k}$, $\zsu_{j=1}^{k-s}d{\tilde
x}_{j}\zex d{\tilde y}_{j}$; so

\centerline{$\zgran\zL_{L}=\zsu_{j=1}^{k-s}{\frac
{\zpar} {\zpar{\tilde x}_{j}}}\zex{\frac {\zpar} {\zpar{\tilde
y}_{j}}}$.}
\bigskip

{\bf Proposition 3.2.} {\it Consider on a $n$-manifold $N$ a
Veronese web $w$ of codimension $r$. Let $\zL$ and $\zL'$ be the
Liouville-Poisson structures of $T^{*}w(0)$ and $T^{*}w(\zinf)$
respectively, and let $\zf_{\zlma}:T^{*}w(0)\zfl T^{*}w(\zinf)$ be
the vector bundle isomorphism  defined by
$\zf_{\zlma}(\zb)=\zb\zci{\zlma}$ where $\zlma:w(\zinf)\zfl w(0)$
is the canonical isomorphism attached to $w$. Note $\zL_{1}$ the
pull-back of $\zL'$ by $\zf_{\zlma}$ (regarded as a
diffeomorphism).

Then $(\zL,\zL_{1})$ is a Kronecker bihamiltonian structure on
$T^{*}w(0)$ with $r$ factors, whose soul $\mathcal{D}$ is given by
the fibres of the canonical fibration $T^{*}w(0)\zfl N$; therefore
the quotient manifold ${\frac{M}{\mathcal{D}}}=N$. Moreover $w$ is
the Veronese web induced by $(\zL,\zL_{1})$ on $N$.}
\bigskip

{\bf Proof.} Let $\zp:T^{*}w(0)\zfl N$ and $\zp':T^{*}w(\zinf)\zfl
N$ be the canonical projections. Choose non-equal and
non-vanishing scalars $\{a_{1},...,a_{n-r},a\}$. On an open
neighbourhood $A$ of a generic point consider a $(1,1)$ tensor
field $J$ like in part (1) of theorem 2.1, coordinates
$(x,y)=(x_{1},...,x_{n-r},y_{1},...,y_{r})$ such that $dx_{j}\zci
J=a_{j}dx_{j}$, $j=1,...,n-r$, and $Ker(dy_{1}\zex...\zex
dy_{r})=w(0)$, and closed $1$-forms $\za_{1},...,\za_{r}$ such
that $Ker(\za_{1}\zex...\zex\za_{r})=w(\zinf)$.

The restriction to $w(0)$ of $dx_{1},...,dx_{n-r}$ is a basis on
$A$ of $T^{*}w(0)$; so on $\zp^{-1}(A)\zeq A\zpor{\mathbb
K}^{n-r}$ one has coordinates $(x,y,u)$, $u=(u_{1},...,u_{n-r})$,
where $(x,y)(\zb)$ are the coordinates of $\zp(\zb)$ and
$\zb=\zsu_{j=1}^{n-r}u_{j}(\zb)dx_{j}$ for each
$\zb\zpe\zp^{-1}(A)$. In the same way one constructs coordinates
$(x,y,u')$, $u'=(u'_{1},...,u'_{n-r})$, on $(\zp')^{-1}(A)$.

In this kind of coordinates, $\zL$ is given by
$dy_{1},...,dy_{r}$, $\zsu_{j=1}^{n-r}dx_{j}\zex du_{j}$ while
$\za_{1},...,\za_{r}$, $\zsu_{j=1}^{n-r}dx_{j}\zex du'_{j}$ define
$\zL'$. On the other hand

\centerline{$\zf_{\zlma}(x,y,u)=(x,y,a_{1}u_{1},...,a_{n-r}u_{n-r})$}

\noindent since $J$
is an extension of $\zlma$ and each $dx_{j}\zci J=a_{j}dx_{j}$.
Therefore $\za_{1},...,\za_{r}$, $\zsu_{j=1}^{n-r}a_{j}dx_{j}\zex
du_{j}$ define $\zL_{1}$. By lemma 3.1 (here $\zw=0$)
$\zL+t\zL_{1}$, $t\znope\{-a_{1},...,-a_{n-r},-a\}$, is given by
$\za_{1}\zci(J+tI)^{-1},...,\za_{r}\zci(J+tI)^{-1}$ and the closed
$2$-form  $\zsu_{j=1}^{n-r}a_{j}(t+a_{j})^{-1}dx_{j}\zex du_{j}$,
which shows the compatibility of $\zL$ and $\zL_{1}$.

The remainder statements are obvious. $\square$
\bigskip

{\bf 4. The equation $d(df\zci J) = \zw$ modulo $I(E)$}

By technical reasons for studying the equation above we shall need parameters that will be
regarded as transverse variables to a $n$-foliation ${\mathcal F}$ defined on a
$m$-dimensional manifold $M$. Let $E$ be an involutive vector subbundle of ${\mathcal F}$ of
dimension $n-r$ where $r\geq 1$.
Consider along ${\mathcal F}$ a diagonalizable (1, 1) tensor field $J$ with characteristic polynomial
$(t-a)^{r}\zpr_{j=1}^{n-r}(t-a_{j})$ where $a_{1},..., a_{n-r}, a$ are non-equal
scalars. Suppose that its Nijenhuis torsion $N_J$ vanishes.

    Let $E^c$ and $I(E)$ be the annihilator of $E$ on ${\mathcal F}^{*}$ and
the differential ideal spanned by the sections of $E^c$ respectively. Assume that $(E^{c}, J^{*})$ spans
${\mathcal F}^{*}$ and that for all closed 1-form $\za$ belonging to $I(E)$ the 2-form
$d(\za\zci J)$ belongs to $I(E)$ as well, where $d$ is the exterior derivative along ${\mathcal F}$.

    As $N_J =0$, distributions $Im(J-a_{j}I)$, $j=1,...,n-r$, and
$Im(J-aI)$ are involutive.
Therefore around every point $p\zpe M$ there exist
functions $x_{1},..., x_{n-r}, y_{1}, ... ,y_{r}$ such that
$dx_{1}\zex ...\zex dx_{n-r}\zex dy_{1} \zex ... \zex dy_{r}$ is a volume form on ${\mathcal F}$,
$dx_{j}\zci J = a_{j}dx_{j}$, $j = 1,..., n-r$, and $dy_{1}=...=dy_{r}=0$ defines $E$.
Indeed, since $(E^{c}, J^{*})$ spans ${\mathcal F}^{*}$ one has $E\zin Ker(J-aI)=\{0\}$, so $E$
is a direct factor of $ Ker(J-aI)$ in ${\mathcal F}^{*}$.

On the other hand $dy_{k}\zci J = ady_{k} + \zsu_{j=1}^{n-r}f_{kj}dx_{j}$, $k=1,...,r$.
As $(E^{c},J^{*})$ spans ${\mathcal F}^{*}$, by linearly recombining functions
$y_{1},..., y_{r}$ and considering $b_{j}x_{j}$
instead $x_j$ for a suitable $b_{j}\zpe{\mathbb K}-\{0\}$, {\it from now on we may
assume that every $f_{1j}(p)$, $j = 1,..., n-r$, is a positive real number}.

    Set $\tilde\za_{k} = \zsu_{j=1}^{n-r}f_{kj}dx_{j}$, $k=1,...,r$. Since $d(dy_{k}\zci J)$
belongs to $I(E)$ one has $dy_{1}\zex ...\zex dy_{r}\zex  d\tilde\za_{k} =0$.
On the other hand vector fields

\centerline{${\partial} /{\partial}x_{1},...,
{\partial}/{\partial}x_{n-r}, {\partial}/{\partial}y_{1} , ...,
{\partial}/{\partial}y_{r}$}

\noindent are defined as the dual basis of $dx_{1},..., dx_{n-r},
dy_{1},..., dy_{r}$.

    In the domain of functions $x_{1},..., x_{n-r}, y_{1}, ..., y_{r}$,
we consider the submanifold $S$ defined by
$x_j - x_{n-r}= x_{j}(p) - x_{n-r}(p)$, $j = 1,..., n-r-1$ ($S =M$ if $n = r, r+1$). Denote
by ${\mathcal F}\zin S$ the $(r+1)$-foliation induced by ${\mathcal F}$ on $S$.

Given a 1-form $\zb$ along ${\mathcal F}$ defined on a open set $M'\zco M$, we denote by $\zb'$
its restriction to $S\zin M'$ {\it as a section of ${\mathcal F}^{*}$}. That is to say $\zb'$
is a section of ${\mathcal F}^{*}$ over $S\zin M'$ and $\zb \zfl\zb'$ is a linear map.
Recall that if $\zm$ is a section of $\zL^{k}{\mathcal F}^{*}$ on $S\zin M'$, its restriction
$\zm_{\zbv {\mathcal F}\zin(S\zin M')}$ can be considered as a $k$-form
on ${\mathcal F}\zin(S\zin M')$. In our particular case when $\zb$ is closed,
$\zb'_{\zbv {\mathcal F}\zin (S\zin M')}$ is closed as well.

Hereafter {\it the standard case} will mean that the structures considered are complex, real
analytic, or $C^{\zinf}$ with $r=1$ in this last case.

    Let $\za_0$ be a 1-form on ${\mathcal F}$.
\bigskip

    {\bf Theorem 4.1.} {\it Suppose that each $f_{1j}(p)$, $j = 1,..., n-r$,
is a positive real number. Then in the standard case the linear map
$\zb\zfl\zb'$ defines an injective correspondence between germs, at $p$, of closed 1-forms
$\zb$ on ${\mathcal F}$ such that

\centerline{$(d(\zb\zci J) + \zb\zex\za_{0})\zex dy_{1} \zex ... \zex dy_{r} =0$}

\noindent and germs, at $p$ on $S$, of sections $\zb'$ of ${\mathcal F}^{*}$ whose
restriction to ${\mathcal F}\zin S$ are closed.

    When $\za_0 =0$ this correspondence becomes bijective.}
\bigskip

    We shall prove this theorem by induction on $n$. For $n = r, r+1$ the result is
obvious since $S =M$. Now assume that the theorem holds up to $n-1$ (whichever $m$ and
$a_{1},..., a_{n-r}, a$ are).

By sake of convenience {\it we will suppose} $a_1 =0$ by replacing $J$ by $J-a_{1}I$
(the equation of theorem 1 does not change
because $d(\zb\zci I)=d\zb=0$).
Set $\za_0 = \zsu_{j=1}^{n-r}h_{j}dx_{j} + \zsu_{k=1}^{r}h_{n+k-r}dy_{k}$ and
$\zb = \zsu_{j=1}^{n-r}\zfi_{j}dx_{j} + \zsu_{k=1}^{r}\zfi_{n+k-r}dy_{k}$. Since
$d\zb=0$ and $dy_{1}\zex...\zex dy_{r}\zex d{\tilde\za}_{k}=0$ we have:
\medskip

\noindent  $\zgran (d(\zb\zci J) + \zb\zex\za_{0})\zex dy_{1}\zex ...\zex dy_{r}=$

\noindent $\zgran  dx_{1}\zex\zsu_{j=2}^{n-r}\zpizq a_{j}{\frac{{\partial}\zfi_{j}} {{\partial}x_{1}}}+
\zsu_{k=1}^{r}(f_{kj}{{\partial}\zfi_{1}\over {\partial}y_{k}}
- f_{k1}{{\partial}\zfi_{j}\over {\partial}y_{k}})  + h_{j}\zfi_{1} - h_{1}\zfi_{j}\zpder dx_{j}
\zex dy_{1}\zex ...\zex dy_{r}$

\hskip 4truecm $\zgran + \zsu_{2\leq i<j\leq n-r}{\tilde h}_{ij}dx_{i}\zex
dx_{j}\zex dy_{1}\zex ...\zex dy_{r}$.
\medskip

    Therefore the part of $d(\zb\zci J) + \zb\zex\za_{0}$ which is divisible by
$dx_1$ modulo $dy_{1},..., dy_{r}$ vanishes if and only if the following system holds:

$\,$

\noindent (1) \hskip 1truecm $\zgran a_{j}{{\partial}\zfi_{j}\over {\partial}x_{1}}
+ \zsu_{k=1}^{r}(f_{kj}{{\partial}\zfi_{1}\over {\partial}y_{k}}
- f_{k1}{{\partial}\zfi_{j}\over {\partial}y_{k}})
+ h_{j}\zfi_{1} - h_{1}\zfi_{j} =0$, $j = 2,..., n-r$.

$\,$

    Let $S'$ be the submanifold defined by $x_j - x_{n-r}= x_{j}(p) - x_{n-r}(p)$, $j = 2,...,
n-r-1$ ($S' =M$ if $n = r+2$). By construction $S$ is a 1-codimension submanifold
of $S'$ and the induced foliation ${\mathcal F}\zin S'$ has dimension $r+2$.

    Set $z_1 =x_1$, $z_2 = x_{n-r}$, $z_3 = y_{1}$,..., $z_{r+2} = y_{r}$. Let
${{\partial}/ {\partial}z_{1}},..., {{\partial}/
{\partial}z_{r+2}}$ be the dual basis of the restriction of $dz_{1},...,
dz_{r+2}$ to ${\mathcal F}\zin S'$. Vector fields
${{\partial}/ {\partial}x_{1}}$, ${{\partial}/ {\partial}x_{2}} +...+ {{\partial}/
{\partial}x_{n-r}}$, ${{\partial}/ {\partial}y_{k}}$, $k=1,...,r$, are tangent to
${\mathcal F}\zin S'$; even more ${{\partial}/ {\partial}z_{1}} = {{\partial}/
{\partial}x_{1}}$, ${{\partial}/ {\partial}z_{2}} = {{\partial}/ {\partial}x_{2}}
+...+ {{\partial}/ {\partial}x_{n-r}}$ and ${{\partial}/ {\partial}z_{k+2}} =
{{\partial}/ {\partial}y_{k}}$, $k=1,...,r$, on $S'$. Besides $dx_1 = dz_1$, $dy_{k} =
dz_{k+2}$, $k=1,...,r$, and the restriction to ${\mathcal F}\zin S'$ of each $dx_j$, $j =
2,..., n-r$, equals that of $dz_2$.

    On $S'$ system (1) becomes:

$\,$

\noindent (2) \hskip 1truecm $\zgran a_{j}{{\partial}\zfi_{j}\over {\partial}z_{1}}
+  \zsu_{k=1}^{r}(f_{kj}{{\partial}\zfi_{1}\over {\partial}z_{k+2}}
- f_{k1}{{\partial}\zfi_{j}\over {\partial}z_{k+2}}) + h_{j}\zfi_{1} - h_{1}\zfi_{j}
=0$, $j = 2,..., n-r$.

$\,$

    The restriction of $\zb$ to ${\mathcal F}\zin S'$, whose expression is

\centerline {$\zgran{\zfi_{1}dz_{1}}_{\zbv {\mathcal F}\zin S'} +
(\zsu_{j=2}^{n-r}\zfi_{j}){dz_{2}}_{\zbv {\mathcal F}\zin S'} +
\zsu_{k=1}^{r}\zfi_{n+k-r}{dz_{k+2}}_{\zbv {\mathcal F}\zin S'}$}

\noindent is a closed 1-form. Hence

\centerline {$\zgran{{\partial}\zfi_{1}\over {\partial}z_{2}}
- \zsu_{j=2}^{n-r}{{\partial}\zfi_{j}\over {\partial}z_{1}} =0$.}

    Now on $S'$ we can consider the system:

$\,$

\noindent (3) \hskip 1truecm $\cases{\zgran{{\partial}\zfi_{1}\over {\partial}z_{2}}
- \zsu_{j=2}^{n-r}{{\partial}\zfi_{j}\over {\partial}z_{1}} =0 \cr\zgran
a_{j}{{\partial}\zfi_{j}\over {\partial}z_{1}}
+  \zsu_{k=1}^{r}(f_{kj}{{\partial}\zfi_{1}\over {\partial}z_{k+2}}
- f_{k1}{{\partial}\zfi_{j}\over {\partial}z_{k+2}}) + h_{j}\zfi_{1} - h_{1}\zfi_{j}
=0\,\, ; j = 2,..., n-r. \cr}$

$\,$
\bigskip

    {\bf Lemma 4.1.} {\it  In the standard case, given a germ at $p$
on $S$ of functions $({\hat\zfi}_{1},...,
{\hat\zfi}_{n-r})$ there exists one and only one germ, at $p$ on $S'$, of functions
$(\zfi_{1},..., \zfi_{n-r})$ which is a solution to (3) and such that
${\zfi_{j}}_{\zbv S} = {\hat\zfi}_{j}$, $j = 1,..., n-r$.}
\bigskip

{\bf Proof.} Consider functions $u_{1},..., u_{m-n}$, on a neighbourhood of $p$ on $S'$, which are
basic for ${\mathcal F}\zin S'$ and such that $(z_{1},..., z_{r+2}, u_{1},..., u_{m-n})$ is a
system of coordinates. Since $ u_{1},..., u_{m-n}$ are basic for ${\mathcal F}\zin S'$ vector
fields ${\partial}/{\partial}z_{1},..., {\partial}/{\partial}z_{r+2}$
defined above equal to partial derivative vector fields, with the same name,
which are associated to coordinates $(z_{1},..., z_{r+2}, u_{1},..., u_{m-n})$.

    Therefore (3) can be regarded like a system on an open set of ${\mathbb K}^{m+r+2-n}$, with
coordinates $(z_{1},..., z_{r+2}, u_{1},..., u_{m-n})$, while $S$ is identify to the
hypersurface defined by $z_1 - z_2 = z_{1}(p) - z_{2}(p)$. In particular
${\partial}/{\partial}z_{1} - {\partial}/{\partial}z_{2}$ is normal to $S$.

In this system ${\partial}/{\partial}z_{1} - {\partial}/{\partial}z_{2}$ is represented by
an invertible triangular matrix with entries on the diagonal $-1, a_{2},..., a_{n-r}$.
Therefore in the complex case or in the real analytic one, lemma
4.1 follows from the Cauchy-Kowalewsky theorem.

Now one will proves the result in the $C^{\zinf}$ case when $r=1$.

Set $f_{j}=f_{1j}$. By adding up to the first equation the second one
multiplied by $a_{2}^{-1}$, the third
one multiplied by $a_{3}^{-1}$, etc..., we obtain the system:

$\,$

\noindent (4) \hskip 1truecm $\cases{\zgran {{\partial}\zfi_{1}\over {\partial}z_{2}}
+(\zsu_{j=2}^{n-1}a_{j}^{-1}f_{j}){{\partial}\zfi_{1}\over {\partial}z_{3}}
- \zsu_{j=2}^{n-1}a_{j}^{-1}f_{1}{{\partial}\zfi_{j}\over {\partial}z_{3}}
+\zsu_{j=2}^{n-1}a_{j}^{-1}(h_{j}\zfi_{1} - h_{1}\zfi_{j}) =0 \cr\zgran
a_{j}{{\partial}\zfi_{j}\over {\partial}z_{1}}
+ f_{j}{{\partial}\zfi_{1}\over {\partial}z_{3}}
- f_{1}{{\partial}\zfi_{j}\over {\partial}z_{3}} + h_{j}\zfi_{1} - h_{1}\zfi_{j}
=0\,\, ; j = 2,..., n-1. \cr}$

$\,$

    In this system ${\partial}/{\partial}z_{1}$ and ${\partial}/{\partial}z_{2}$ are
represented by diagonal matrices with entries on the diagonal $0, a_{2},..., a_{n-1}$ and
$1, 0, ..., 0$ respectively.

    On the other hand ${\partial}/{\partial}z_{3}$ is represented by the matrix:

$$\pmatrix{\zgran \zsu_{j=2}^{n-1}a_{j}^{-1}f_{j} & -a_{2}^{-1}f_{1} & -a_{3}^{-1}f_{1} &
. & . & . & -a_{n-1}^{-1}f_{1} \cr\zgran
f_{2} & -f_{1} & \quad & \quad & \quad & \quad & \quad  \cr\zgran
f_{3} & \quad & -f_{1} &  \quad & \quad & \quad & \quad  \cr\zgran
. &  \quad & \quad & . \quad & \quad & \quad  \cr\zgran
. &  \quad & \quad & \quad & . & \quad & \quad  \cr\zgran
. &  \quad & \quad & \quad &\quad & . & \quad  \cr\zgran
f_{n-1} & \quad & \quad & \quad & \quad & \quad & -f_{1} \cr
}$$

$\,$

    Obviously each ${\partial}/{\partial}u_{i}$ is represented by the zero matrix.

If one multiplies the jth equation, $j = 2,..., n-1$, by $-a_{j}^{-1}f_{1}f_{j}^{-1}$, we
obtain  a linear symmetric system. In this new system
${\partial}/{\partial}z_{1} - {\partial}/{\partial}z_{2}$ is represented by a diagonal
matrix with entries on the diagonal $-1, -f_{1}f_{2}^{-1},..., -f_{1}f_{n-1}^{-1}$. This
matrix is negative definite around $p$, then the new system is symmetric hyperbolic and $S$ is
space-like.

Therefore this case of lemma 4.1 follows from the classical results on the Cauchy problem
\cite{CO}, \cite{TA}. $\square$

    Let us come back to the proof of theorem 4.1.
\medskip

{\bf Uniqueness.} Let $\zb = \zsu_{j=1}^{n-r}\zfi_{j}dx_{j} +
\zsu_{k=1}^{r}\zfi_{n+k-r}dy_{k}$ and
$\zg = \zsu_{j=1}^{n-r}\zf_{j}dx_{j} + \zsu_{k=1}^{r}\zf_{n+k-r}dy_{k}$ be two solutions to
the equation of theorem 4.1, such that $\zb' = \zg'$. On $S'$ functions
$\zfi_{1},...,\zfi_{n-r}$ and
$\zf_{1},..., \zf_{n-r}$ are solutions to (3), which agree on $S$, then by lemma 4.1 we have
$\zfi_{j} = \zf_{j}$, $j = 1,..., n-r$, as germs at $p$ on $S'$.

    The restriction of $\zb -\zg$ to $S'$, which equals
$\zsu_{k=1}^{r}(\zfi_{n+k-r}-\zf_{n+k-r})dy_{k\zbv {\mathcal F}\zin S'}$, is closed. Therefore
each $\zfi_{n+k-r}-\zf_{n+k-r}$, $k=1,...,r$, is constant on the leaves of the foliation
defined by $Ker(dy_{1}\zex ...\zex dy_{r})_{\zbv {\mathcal F}\zin S'}= E\zin S'$.
But $S$ is transverse to this foliation and
$(\zfi_{n+k-r}-\zf_{n+k-r})_{\zbv S} = 0$ then $\zfi_{n+k-r}=\zf_{n+k-r}$, $k=1,...,r$,
on $S'$. In other words $\zb$ and $\zg$ agree on $S'$ as
sections of ${\mathcal F}^{*}$.

    The next step will be to regard $x_1$ like a new parameter. By shrinking $M$ we may
suppose that function $x_1$ is defined on the whole $M$.

Set ${\mathcal F'}=Ker dx_{1} \zco {\mathcal F}$, which is a $(n-1)$-foliation,
and let $d'$ be the the exterior derivative along it.
Denote by $J'$ and $\za'_{0}$ the restriction
to ${\mathcal F'}$ of $J$ and $\za_0$ respectively (recall that $dx_{1}\zci J =0$).
Set $E'=E\zin {\mathcal F'}$. Let $E'^{c}$ and $I(E')$ be the annihilator of $E'$ on
$({\mathcal F'})^{*}$  and the differential ideal spanned by the sections of $E'^{c}$ respectively.
Then $(E'^{c}, J'^{*})$ spans $({\mathcal F'})^{*}$ and, for any closed 1-form $\zt$ belonging to
$I(E')$, the 2-form $d'(\zt\zci J')$ belongs to $I(E')$ as well. On the other hand
$d'x_{j}\zci J' = a_{j}d'x_{j}$, $j = 2,..., n-r$,
$d'y_{k}\zci J'=\zsu_{j=2}^{n-r}f_{kj}d'x_{j} + ad'y_{k}$, $k=1,...,r$,
and $d'y_{1}=...=d'y_{r}=0$ defines $E'$.

    Since $S'$ plays the same role with respect to $(x_{2},..., x_{n-r}, y_{1},..., y_{r})$ as
$S$ does with respect to $(x_{1},..., x_{n-r}, y_{1},..., y_{r})$, $\zb_{\zbv {\mathcal F'}}$
and $\zg_{\zbv {\mathcal F'}}$ satisfy to the equation of theorem 4.1 for ${\mathcal F'}$, $J'$,
$E'$ and $\za'_{0}$, and
$\zb_{\zbv {\mathcal F'}}=\zg_{\zbv {\mathcal F'}}$ on $S'$, from the induction hypothesis
follows that $\zb_{\zbv {\mathcal F'}}=\zg_{\zbv {\mathcal F'}}$
like germs at $p$ on $M$, i.e. $\zfi_j = \zf_j$, $j = 2,..., n$.

    Finally, as $\zb -\zg = (\zfi_{1}-\zf_{1})dx_1$ is closed, function $\zfi_{1}-\zf_{1}$ is
constant along the leaves of ${\mathcal F'}$. But $S$ is transverse to ${\mathcal F'}$ and
$(\zfi_{1}-\zf_{1})_{\zbv S}=0$ then $\zfi_{1}=\zf_{1}$ and $\zb=\zg$ as germs at $p$
on $M$.
\medskip

    {\bf Existence.}  Now $\za_0 =0$, i.e. $h_{1} = ...= h_n =0$. Given functions
$\zfi_{1},..., \zfi_{n}$ on $S$ such that the restriction of
$\zb' = \zsu_{j=1}^{n-r}\zfi_{j}dx_{j} + \zsu_{k=1}^{r}\zfi_{n+k-r}dy_{k}$
to ${\mathcal F}\zin S$ is closed, by means of system (3) we extend functions $\zfi_{1},...,
\zfi_{n-r}$ to $S'$ (around $p$).

    Since $\zfi_{1}{dz_{1}}_{\zbv {\mathcal F}\zin S'} + (\zsu_{j=2}^{n-r}\zfi_{j}){dz_{2}}_{\zbv
{\mathcal F}\zin S'}$ is closed modulo ${dz_{k+2}}_{\zbv {\mathcal F}\zin S'}$, $k=1,...,r$,
(first equation of (3)), there exist functions ${\hat\zfi}_{n+1-r},..., {\hat\zfi}_{n}$ on
$S'$ such that the restriction to ${\mathcal F}\zin S'$ of
$\zfi_{1}dz_{1} + (\zsu_{j=2}^{n-r}\zfi_{j})dz_{2} + \zsu_{k=1}^{r}{\hat\zfi}_{n+k-r}
dz_{k+2}$ is closed. Consequently its restriction to ${\mathcal F}\zin S$ is closed as well.
On the other hand, by hypothesis, the restriction to ${\mathcal F}\zin S$ of
$\zfi_{1}dz_{1} + (\zsu_{j=2}^{n-r}\zfi_{j})dz_{2} + \zsu_{k=1}^{r}\zfi_{n+k-r} dz_{k+2}$ is
closed. Therefore
${\zsu_{k=1}^{r}({\hat\zfi}_{n+k-r} - \zfi_{n+k-r})dz_{k+2}}_{\zbv {\mathcal F}\zin S}$ is
closed.

    In coordinates $(z_{1},..., z_{r+2}, u_{1},..., u_{m-n})$ like in the proof of lemma 4.1,
this implies the existence, on $S'$, of a function  $h(z_{3},..., z_{r+2}, u_{1},...,
u_{m-n})$ such that $dh_{\zbv {\mathcal F}\zin S} = {\zsu_{k=1}^{r}({\hat\zfi}_{n+k-r} -
\zfi_{n+k-r})dz_{k+2}}_{\zbv {\mathcal F}\zin S}$ on $S$.
Obviously functions ${\hat\zfi}_{n+k-r} - {\partial}h/{\partial}z_{k+2}$, $k=1,...,r$,
have the same property as functions ${\hat\zfi}_{n+k-r}$, $k=1,...,r$. Then by
replacing each ${\hat\zfi}_{n+k-r}$ by
${\hat\zfi}_{n+k-r} - {\partial}h/{\partial}z_{k+2}$, we can suppose that
${\hat\zfi}_{n+k-r}$ is an extension of $\zfi_{n+k-r}$ and call it $\zfi_{n+k-r}$ from now
on.

If we consider ${\mathcal F'}$, $J'$, $E'$ and the section of $({\mathcal F'})^{*}$ over $S'$:

\noindent $\zsu_{j=2}^{n-r}\zfi_{j}d'x_{j} + \zsu_{k=1}^{r}\zfi_{n+k-r}d'y_{k}$, whose
restriction to  ${\mathcal F'}\zin S'$ is closed, the induction hypothesis allows us to
extend functions $\zfi_{2},..., \zfi_{n}$ to an open set of $M$ containing $p$, in such a way
that ${\bar\zb}= \zsu_{j=2}^{n-r}\zfi_{j}d'x_{j} + \zsu_{k=1}^{r}\zfi_{n+k-r}d'y_{k}$  is a
closed 1-form along ${\mathcal F'}$ and $d'({\bar\zb}\zci J') \zex d'y_{1}\zex ...\zex d'y_{r}
=0$.

Since $d'{\bar\zb} =0$ there exists a function $\zf$ such that

\noindent $\zr = \zf dx_{1}+
\zsu_{j=2}^{n-r}\zfi_{2}dx_{j} + \zsu_{k=1}^{r}\zfi_{n+k-r}dy_{k}$ is a closed form along
${\mathcal F}$. On the other hand

\centerline {$\zr_{\zbv {\mathcal F}\zin S'} -
(\zfi_{1}dz_{1}+ {(\zsu_{j=2}^{n-r}\zfi_{j})dz_{2} +
\zsu_{k=1}^{r}\zfi_{n+k-r}dz_{k+2})}_{\zbv {\mathcal F}\zin S'} = (\zf
-\zfi_{1}){dz_{1}}_{\zbv {\mathcal F}\zin S'}$}

\noindent is closed; i.e. $\zf -\zfi_{1}$
is constant on the leaves of the foliation associated
to $Ker{dz_{1}}_{\zbv {\mathcal F}\zin S'}$.

Around $p$ on $M$ consider coordinates $(x_{1},..., x_{n-r}, y_{1},..., y_{r}, v_{1},...,
v_{m-n})$ where
$v_{1},..., v_{m-n}$ are basic functions for ${\mathcal F}$. Then  as $x_1 = z_1$
there exists a function ${\bar h}(x_{1}, v_{1},..., v_{m-n})$, around $p$ on $M$, such
that $\zf-\zfi_{1}={\bar h}$ on $S'$ and, by replacing $\zf$ by $\zf-{\bar h}$, we
may suppose that $\zf$ extends $\zfi_1$ and call $\zfi_1$ this extension too.

    In short we have constructed a closed 1-form, along ${\mathcal F}$,

 \noindent $\zb = \zsu_{j=1}^{n-r}\zfi_{j}dx_{j} + \zsu_{k=1}^{r}\zfi_{n+k-r}dy_{k}$
which extends $\zb'$ and such that $d(\zb\zci J)\zex dx_{1}\zex dy_{1}\zex ...\zex dy_{r}=0$
(this is another way for writing
$d'({\bar\zb}\zci J')\zex d'y_{1}\zex ...\zex d'y_{r}=0$). Therefore there exist closed
1-forms $\zg_{0},..., \zg_{r}$ along ${\mathcal F}$ such that

\centerline{$d(\zb\zci J)= dx_{1}\zex \zg_{0} + \zsu_{k=1}^{r}\zg_{k}\zex dy_{k}$.}

    Set $\zg_0 = \zsu_{j=1}^{n-r}g_{j}dx_{j} + \zsu_{k=1}^{r}g_{n+k-r}dy_{k}$. Then

$\,$

\centerline {$\zgran g_{j} = a_{j}{{\partial}\zfi_{j}\over {\partial}x_{1}}
+ \zsu_{k=1}^{r}(f_{kj}{{\partial}\zfi_{1}\over {\partial}y_{k}}
- f_{k1}{{\partial}\zfi_{j}\over {\partial}y_{k}})\,$; $j = 2,..., n-r$}

$\,$

\noindent (recall the construction of system (1)). Therefore each $g_j$, $j = 2,..., n-r$,
vanishes on $S'$ because $\zfi_{1},..., \zfi_{n-r}$ satisfy to system (3).

    On the other hand $(d(\zb\zci J))_{J}$ is closed (apply lemma 2.1 along the leaves
of ${\mathcal F}$). Then

\centerline {$-dx_{1}\zex d(\zg_{0}\zci J) +  \zsu_{k=1}^{r}(d(\zg_{k}\zci J) \zex dy_{k} -
\zg_{k}\zex d({\tilde\za}_{k} + ady_{k}))=0\,\,$,}

\noindent whence $dx_{1} \zex d(\zg_{0}\zci J)\zex dy_{1}\zex ...\zex dy_{r}=0$.
That is to say $d'({\bar\zg}_{0}\zci J')\zex d'y_{1}\zex ...\zex d'y_{r}=0$ where
${\bar\zg}_{0} = \zsu_{j=2}^{n-r}g_{j}d'x_{j} + \zsu_{k=1}^{r}g_{n+k-r}d'y_{k}$.

    On $S'$, ${\bar\zg}_{0}$ is a combination of $d'y_{1},..., d'y_{r}$. Since the restriction
of ${\bar\zg}_{0}$ to ${\mathcal F'}\zin S'$ is closed there exists
a function $\zlma (x_{1}, y_{1},..., y_{r}, v_{1},..., v_{m-n})$,
defined near $p$ on $M$, such that ${\bar\zg}_{0}=
d'\zlma $ on $S'$.

    But $d'\zlma$ is a closed 1-form along ${\mathcal F'}$ defined on an
open set of $M$ and

\noindent $d'(d'\zlma \zci J')\zex d'y_{1}\zex ...\zex d'y_{r}=0$. Therefore the uniqueness
in dimension $n-1$ implies that ${\bar\zg}_{0}=d'\zlma$.
In other words $\zg_0$ is a combination of $dx_{1}, dy_{1},..., dy_{r}$. Then $d(\zb\zci
J)\zex dy_{1}\zex ...\zex dy_{r}=0$ and {\it the proof of theorem 4.1 is finished}.
\medskip

The following result will be needed in the next section.
\bigskip

{\bf Lemma 4.2.} {\it Suppose that each $f_{1j}(p)$,
$j=1,...,n-r$,
is a positive real number. Consider $1$-forms $\zr_{\zlma q}$,
$\zlma,q=1,...,s$. In the standard case, given two families of $s$
closed $1$-forms, which are solution to the system

\centerline{$(d(\zb_{q}\zci J) + \sum_{\zlma=1}^{s}\zb_{\zlma}
\zex\zr_{\zlma q})\zex dy_{1} \zex ... \zex dy_{r} =0$, $q=1,...,s$,}

\noindent if they agree around $p$ on $S$ then they agree around
$p$ on $M$.}
\bigskip

{\bf Proof.} Just adapt the proof of the uniqueness of theorem 4.1
(in fact the case $s=1$ is the first assertion of this theorem).
Now system $(3)$ is replaced by a system ${\mathcal
S}(\zb_{1},...,\zb_{s})$ with $s$ boxes corresponding each of them
to a $\zb_{q}$. Note that the symbol of every box, which only
depends on $\zb_{q}$, is similar to the symbol of system $(3)$.
Therefore lemma 4.1 extends to ${\mathcal
S}(\zb_{1},...,\zb_{s})$. Finally if $\zb_{q}=\zsu_{j=1}^{n-r}\zfi_{qj}dx_{j}
+ \zsu_{k=1}^{r}\zfi_{qn+k-r}dy_{k}$ and $\zg_{q}=\zsu_{j=1}^{n-r}\zf_{qj}dx_{j}
+ \zsu_{k=1}^{r}\zf_{qn+k-r}dy_{k}$, $q=1,...,s$, are two
solutions to the system of lemma 4.2 such that
$\zb'_{q}=\zg'_{q}$, $q=1,...,s$, reasoning as in the proof of the
uniqueness of theorem 4.1 shows that $\zb_{q}=\zg_{q}$,
$q=1,...,s$. $\square$
\bigskip

{\bf Theorem 4.2.} {\it Suppose that every $f_{1j}(p)$, $j = 1,..., n-r$, is a
positive real number. In the standard case given, on an open
neighbourhood of $p$ on $M$, a closed 1-form $\zg$ along ${\mathcal F}$ such that
$d(\zg\zci J)\zex dx_{1}\zex dy_{1}\zex ...\zex dy_{r}=0$, then around $p$ there exists a
closed 1-form $\zb$ along ${\mathcal F}$ such that $d(\zb\zci J)\zex dy_{1}\zex ...\zex
dy_{r}=dx_{1}\zex \zg \zex dy_{1}\zex ...\zex dy_{r}$.}
\bigskip

{\bf Proof.} As above we shall suppose that $a_{1}=0$ by  replacing, if necessary, $J$
by $J-a_{1}I$. Set $\zg = \zsu_{j=1}^{n-r}\zf_{j}dx_{j} + \zsu_{k=1}^{r}\zf_{n+k-r}dy_{k}$.

    On $S'$ we consider the following system:

$\,$

\noindent (3') \hskip 1truecm $\cases{\zgran{{\partial}\zfi_{1}\over {\partial}z_{2}}
- \zsu_{j=2}^{n-r}{{\partial}\zfi_{j}\over {\partial}z_{1}} =0 \cr\zgran
a_{j}{{\partial}\zfi_{j}\over {\partial}z_{1}}
+  \zsu_{k=1}^{r}(f_{kj}{{\partial}\zfi_{1}\over {\partial}z_{k+2}}
- f_{k1}{{\partial}\zfi_{j}\over {\partial}z_{k+2}}) = \zf_{j} \,\, ; j = 2,..., n-r. \cr}$

$\,$

    This system has some solution around $p$ because its symbol is the same as that of system
(3). Let $\zfi_{1},..., \zfi_{n-r}$ be a solution to (3'). The first equation of (3') allows
us to find functions $\zfi_{n+1-r},..., \zfi_n$, on a neighbourhood of $p$ on $S'$, such
that ${\zpizq\zfi_{1}dz_{1} + (\zsu_{j=2}^{n-r}\zfi_{j})dz_{2} + \zsu_{k=1}^{r}\zfi_{n+k-r}
dz_{k+2}\zpder}_{\zbv {\mathcal F}\zin S'}$ is closed. Obviously the restriction of this form to
${\mathcal F'}\zin S'$ is closed too.

    Now we apply theorem 1 to ${\mathcal F'}$, $J'$ and $E'$ for extending functions
$\zfi_{2},..., \zfi_{n}$ to an open set of $M$ containing $p$, in such a way that
$d'{\bar\zb}=0$ and $d'({\bar\zb}\zci J')\zex d'y_{1}\zex ...\zex d'y_{r}=0$
where ${\bar\zb}= \zsu_{j=2}^{n-r}\zfi_{j}d'x_{j} + \zsu_{k=1}^{r}\zfi_{n+k-r}d'y_{k}$.

    The rest of the proof is very similar to that of the existence in theorem 1. First we
extend function $\zfi_1$ to a neighbourhood of $p$ on $M$ in such a way that
$\zb = \zsu_{j=1}^{n-r}\zfi_{j}dx_{j} + \zsu_{k=1}^{r}\zfi_{n+k-r}dy_{k}$ is closed. Since
$d'({\bar\zb}\zci J')\zex d'y_{1}\zex ...\zex d'y_{r}=0$ we get $d(\zb\zci J)\zex dx_{1}\zex
dy_{1}\zex ...\zex dy_{r} =0$. Therefore
$d(\zb\zci J)= dx_{1}\zex \zg_{0} + \zsu_{k=1}^{r}\zg_{k}\zex dy_{k}$
where $\zg_{0},..., \zg_{r}$ are closed 1-forms along ${\mathcal F}$.

    Set $\zg_0 = \zsu_{j=1}^{n-r}g_{j}dx_{j} + \zsu_{k=1}^{r}g_{n+k-r}dy_{k}$ and
${\bar\zg}_{0}= \zsu_{j=2}^{n-r}g_{j}d'x_{j} + \zsu_{k=1}^{r}g_{n+k-r}d'y_{k}$. Then

$\,$

\centerline {$\zgran g_{j} = a_{j}{{\partial}\zfi_{j}\over {\partial}x_{1}}
+ \zsu_{k=1}^{r}(f_{kj}{{\partial}\zfi_{1}\over {\partial}y_{k}}
- f_{k1}{{\partial}\zfi_{j}\over {\partial}y_{k}})\,$; $j = 2,..., n-r$}

$\,$

    Besides $d'({\bar\zg}_{0}\zci J')\zex d'y_{1}\zex ...\zex d'y_{r}=0$ because $(d(\zb\zci
J))_{J}$ is closed (lemma 2.1).

    By hypothesis $d'({\bar\zg}\zci J')\zex d'y_{1}\zex ...\zex d'y_{r}=0$ where

\noindent ${\bar\zg} = \zsu_{j=2}^{n-r}\zf_{j}d'x_{j} + \zsu_{k=1}^{r}\zf_{n+k-r}d'y_{k}$.
On the other hand ${\bar\zg} - {\bar\zg}_{0}$ is a closed 1-form along ${\mathcal F'}$ which is a
combination of $d'y_{1},..., d'y_{r}$ on $S'$ since $(\zfi_{1},..., \zfi_{n-r})$ is a
solution to (3'). This fact implies the existence, on an open neighbourhood of $p$ on $M$, of a
function $\zlma (x_{1}, y_{1},..., y_{r}, v_{1},..., v_{m-n})$ such that ${\bar\zg} -
{\bar\zg}_{0} = d'\zlma$ on $S'$.

    Obviously $d'(d'\zlma\zci J')\zex d'y_{1}\zex ...\zex d'y_{r}=0$.
Now  from theorem 4.1 applied to ${\mathcal F'}$, $J'$ and $E'$
follows that ${\bar\zg} - {\bar\zg}_{0} =
d'\zlma$ around $p$ on $M$.
Hence $({\zg} - {\zg}_{0})\zex dx_{1}\zex dy_{1} \zex ...\zex dy_{r} =0$ and
$d(\zb\zci J)\zex dy_{1} \zex ...\zex dy_{r} = dx_{1}\zex\zg_{0}\zex dy_{1} \zex ...\zex
dy_{r} = dx_{1}\zex\zg\zex dy_{1} \zex ...\zex dy_{r}$. $\square$
\bigskip

    {\bf Theorem 4.3.} {\it Let ${\mathcal F}$ be a $n$-foliation defined on a $m$-manifold $M$ and
let $E \zco {\mathcal F}$ be a second foliation of dimension $n-r$ where $r\geq 1$. On
${\mathcal F}$ we consider a diagonalizable (1, 1) tensor field $J$ with characteristic
polynomial $(t-a)^{r}\zpr_{j=1}^{n-r}(t-a_{j})$ where $a_{1},..., a_{n-r}, a$
are non-equal scalars. Suppose $N_J =0$.

    Let $E^c$ and $I(E)$ be the annihilator of $E$ on ${\mathcal F}^{*}$ and
the differential ideal spanned by the sections of $E^c$ respectively. Assume that
$(E^{c}, J^{*})$ spans ${\mathcal F}^{*}$ and that for all closed 1-form $\za$
belonging to $I(E)$ the 2-form $d(\za\zci J)$ belongs to $I(E)$ as well, where
$d$ is the exterior derivative along ${\mathcal F}$.

    In the standard case, given a closed 2-form $\zw$ on ${\mathcal F}$, the
following statements are equivalents:

\noindent (a) Around each  point $p\zpe M$ there exists a function $f$ such that
$d(df\zci J) =\zw$ modulo $I(E)$.

\noindent (b) $d\zw_{J}$ belongs to $I(E)$.}
\bigskip

{\bf Proof.} $(a) \zim (b)$ In this case locally $\zw = d(df\zci J) +
\zsu_{k=1}^{r}\zm_{k}\zex\za_{k}$ where

\noindent $\zm_{1},..., \zm_{r}, \za_{1},..., \za_{r}$ are
closed 1-forms and $\za_{1},..., \za_{r}$ belong to $I(E)$. Since $(d(df\zci J))_{J}$ is
closed (lemma 2.1) and each 2-form $d(\za_{k}\zci J)$ belongs to $I(E)$, it follows that
$d\zw_{J}$ belongs to $I(E)$.

    $(b) \zim (a)$ As the problem is local we will use the concepts and notations of
the proofs of theorems 4.1 and 4.2. The implication will be proved by induction on $n$. For $n
= r, r+1$ the results is obvious. Now, assume that it holds up to $n-1$ (whichever $m$ is).

    Let ${\bar\zw}$ be the restriction of $\zw$ to ${\mathcal F'}$. Then $d'({\bar\zw}_{J'})\zex
d'y_{1}\zex ...\zex d'y_{r}=0$. By the induction hypothesis there exists a function ${\bar
f}$ around $p$ such that $d'(d'{\bar f}\zci J')\zex d'y_{1}\zex ...\zex d'y_{r} =
{\bar\zw}\zex d'y_{1}\zex ...\zex d'y_{r}$. Hence
$d(d{\bar f}\zci J)\zex dx_{1}\zex dy_{1}\zex ...\zex dy_{r} = \zw\zex dx_{1}\zex
dy_{1}\zex ...\zex dy_{r}$.

    Therefore $\zw - d(d{\bar f}\zci J) = dx_{1}\zex \zg_{0} + \zsu_{k=1}^{r}\zg_{k}\zex
dy_{k}$ where $\zg_{0},..., \zg_{r}$ are 1-closed forms along ${\mathcal F}$.
As $(d\zw_{J})\zex dy_{1}\zex ...\zex dy_{r} =0$ and $(d(d{\bar f}\zci J))_{J}$ is closed
(lemma 2.1) we obtain $d(\zg_{0}\zci J)\zex dx_{1}\zex dy_{1}\zex ...\zex dy_{r}=0$.

    By theorem 4.2,  around $p$ there exists a closed 1-form $\zb$, along ${\mathcal F}$, such that
$d(\zb\zci J) \zex dy_{1}\zex ...\zex dy_{r}= dx_{1}\zex\zg_{0}\zex dy_{1}\zex ...\zex
dy_{r}$. Now it is enough to set $f = h + {\bar f}$ where $h$ is a primitive of $\zb$.
$\square$
\medskip

Finally remark that theorem 3.1 is just the implication $(b) \zim
(a)$ of the foregoing theorem when $n=m$ and
$E=Ker(\za_{1}\zex...\zex\za_{r})$.
\bigskip

{\bf 5. Another equation.}

The aim of this paragraph is to establish another theorem on
some system defined by differential forms, which be needed later
on in the construction of versal models of Veronese webs. The
objects $M$, $\mathcal F$, etc... are as in the foregoing section
unless another thing is stated. Set
$J_{0}=\sum_{j=1}^{n-r}a_{j}{\frac{\zpar} {\zpar x_{j}}}\zte
dx_{j}+\sum_{k=1}^{r}a{\frac{\zpar} {\zpar y_{k}}}\zte dy_{k}$.
\bigskip

{\bf Theorem 5.1.} {\it In the standard case, given a germ at $p$
of maps $\zf_{kj}:S\zfl \mathbb K$, $k=1,...,r$, $j=1,...,n-r$,
such that every $\zf_{1j}(p)$, $j=1,...,n-r$, is a positive real
number, then there exists one and only one germ at $p$ on $M$ of
$1$-forms ${\tilde\za_{1}} =\zsu_{j=1}^{n-r}f_{1j}dx_{j}$,...,
${\tilde\za_{r}} =\zsu_{j=1}^{n-r}f_{rj}dx_{j}$ such that

$\,$

\noindent (4) \hskip .5truecm
 $\cases{ d{\tilde\za}_{k}\zex
dy_{1}\zex...\zex dy_{r}=0,\quad k=1,...,r\cr \, \cr\zgran \zpizq
d({\tilde\za}_{k}\zci J_{0})-\sum_{\zlma=1}^{r}
{\tilde\za}_{\zlma}\zex {\frac{\zpar {\tilde\za}_{k}} {\zpar
y_{\zlma}}}\zpder\zex dy_{1}\zex...\zex dy_{r}=0, \quad
k=1,...,r\cr}$

$\,$

\noindent and that ${f_{kj}}_{\zbv S}=\zf_{kj}$, $k=1,...,r$,
$j=1,..., n-r$.}
\bigskip

We shall prove this theorem by induction on $n$. For $n =r,r+1$
the result is obvious since $S =M$. Now assume that the theorem
holds up to $n-1$ (whichever $m$ and $a_{1},..., a_{n-r},a$ are).

Consider 1-forms ${\tilde\za}_{k} =\zsu_{j=1}^{n-r}f_{kj}dx_{j}$,
$k=1,...,r$, such that $d{\tilde\za}_{k}\zex dy_{1}\zex...\zex
dy_{r}=0$.

    By sake of convenience {\it we will suppose} $a_1 =0$ by
replacing $J_0$ by $J_{0}-a_{1}I$ (the main equation of theorem 1
does not change). Then:

$\,$

\noindent  $\zgran \zpizq d({\tilde\za}_{k}\zci J_{0})  -
\sum_{\zlma=1}^{r}{\tilde\za}_{\zlma}\zex
{{\partial}{\tilde\za}_{k}\over {{\partial}y}_{\zlma}}\zpder\zex
dy_{1}\zex...\zex dy_{r}$

$\zgran= dx_{1}\zex\zsu_{j=2}^{n-r}\zpizq
a_{j}{{\partial}f_{kj}\over {\partial}x_{1}} +
\sum_{\zlma=1}^{r}(f_{\zlma j}{{\partial}f_{k1}\over
{\partial}y_{\zlma}} - f_{\zlma 1}{{\partial}f_{kj}\over
{\partial}y_{\zlma}})\zpder dx_{j}\zex dy_{1}\zex...\zex dy_{r}
\zgran + \zsu_{2\leq 1<j\leq n-r}{\tilde h}_{ij}dx_{i}\zex
dx_{j}\zex dy_{1}\zex...\zex dy_{r}$.

$\,$

    Therefore the part of each $d({\tilde\za}_{k}\zci J_{0})  -
\sum_{\zlma=1}^{r}{\tilde\za}_{\zlma}\zex
{{\partial}{\tilde\za}_{k}\over {{\partial}y}_{\zlma}}$ that is
divisible by $dx_1$ modulo $dy_{1},...,dy_{r}$ vanishes if and
only if the following system holds:

$\,$

\noindent (5) \hskip 1truecm $\zgran a_{j}{{\partial}f_{kj}\over
{\partial}x_{1}} +\sum_{\zlma=1}^{r}( f_{\zlma
j}{{\partial}f_{k1}\over {\partial}y_{\zlma}} - f_{\zlma
1}{{\partial}f_{kj}\over {\partial}y_{\zlma}})=0$, $j = 2,...,
n-r$, $k=1,...,r$.

$\,$

On $S'$ endowed with coordinates $(z_{1},...,z_{r+2})$ system (5)
becomes:

$\,$

\noindent (6) \hskip 1truecm $\zgran a_{j}{{\partial}f_{kj}\over
{\partial}z_{1}} +\sum_{\zlma=1}^{r}( f_{\zlma
j}{{\partial}f_{k1}\over {\partial}z_{\zlma+2}} - f_{\zlma
1}{{\partial}f_{kj}\over {\partial}z_{\zlma+2}})=0$, $j = 2,...,
n-r$, $k=1,...,r$.

$\,$

    The restriction of each ${\tilde\za}_{k}\zex dy_{1}\zex...\zex dy_{r}$
to ${\mathcal F}\zin S'$ whose expression is
\bigskip

\centerline {$\zgran \zpizq \zcizq{f_{k1}dz_{1}} +
(\zsu_{j=2}^{n-r}f_{kj}){dz_{2}}\zcder\zex {dz_{3}}\zex...\zex
dz_{r+2}\zpder_{\zbv {\mathcal F}\zin S'}$}
\bigskip

\noindent is a closed 2-form. Hence

\centerline {$\zgran{{\partial}f_{k1}\over {\partial}z_{2}} -
\zsu_{j=2}^{n-r}{{\partial}f_{kj}\over {\partial}z_{1}} =0$,
$k=1,...,r$.}

    Now on $S'$ we can consider the system:

$\,$

\noindent (7) \hskip 1truecm
$\cases{\cases{\zgran{{\partial}f_{k1}\over {\partial}z_{2}} -
\zsu_{j=2}^{n-r}{{\partial}f_{kj}\over {\partial}z_{1}} =0
\cr\zgran a_{j}{{\partial}f_{kj}\over {\partial}z_{1}}
+\sum_{\zlma=1}^{r}( f_{\zlma j}{{\partial}f_{k1}\over
{\partial}z_{\zlma+2}} - f_{\zlma 1}{{\partial}f_{kj}\over
{\partial}z_{\zlma+2}})=0,\quad j = 2,..., n-r, \cr}\cr \,
\cr\quad k=1,...,r.\cr}$

$\,$
\bigskip

{\bf Lemma 5.1.} {\it In the standard case, given a germ at $p$ on
$S$ of functions $\zf_{kj}$, $k=1,...,r$, $j=1,...,n-r$, such that
every $\zf_{1j}(p)$, $j=1,..., n-r$, is a positive real number,
then there exists one and only one germ, at $p$ on $S'$, of
functions $f_{kj}$, $k=1,...,r$, $j=1,...,n-r$, which is a
solution to (7) and such that ${f_{kj}}_{\zbv S} = \zf_{kj}$,
$k=1,...,r$, $j = 1,..., n-r$.}
\bigskip

{\bf Proof.} On a neighbourhood of $p$ on $S'$ consider
functions $u_{1},..., u_{m-n}$ basic for ${\mathcal
F}\zin S'$ and such that $(z_{1},..., z_{r+2}, u_{1},...,
u_{m-n})$ is a system of coordinates. Since $ u_{1},..., u_{m-n}$
are basic for ${\mathcal F}\zin S'$ vector fields
${\partial}/{\partial}z_{1}, ...,{\partial}/{\partial}z_{r+2}$
defined above equal to partial derivative vector fields, with the
same name, which are associated to coordinates $(z_{1},...,
z_{r+2}, u_{1},..., u_{m-n})$.

Therefore (7) can be regarded like a system on an open set of
${\mathbb K}^{m+r-n+2}$ with coordinates $(z_{1},..., z_{r+2},
u_{1},..., u_{m-n})$, while $S$ is identify to the hypersurface
defined by $z_1 - z_2 = z_{1}(p) - z_{2}(p)$. In particular
${\partial}/{\partial}z_{1} - {\partial}/{\partial}z_{2}$ is
normal to $S$.

In this system the matrix associated to
${\partial}/{\partial}z_{1} - {\partial}/{\partial}z_{2}$ is
invertible. Indeed, it consists of $r$ blocks $(n-r)\zpor(n-r)$
along the diagonal corresponding to the different values of $k$
and zero outside of them, and every block is triangular with
entries $-1,a_{2},...,a_{n-r}$ on the diagonal. So in the complex
case and in the real analytic one it suffices to apply the
Cauchy-Kowalewsky theorem.

On the other hand if $r=1$, systems (3) and (7) have very similar
symbols and, for the $C^{\zinf}$ case, it is enough to reason as in
the proof of lemma 4.1. $\square$

Let us come back to the proof of the theorem 5.1.
\medskip

{\bf Uniqueness.} Let ${\tilde\za}_{k} =
\zsu_{j=1}^{n-r}f_{kj}dx_{j}$ and $\zg_{k} =
\zsu_{j=1}^{n-r}g_{kj}dx_{j}$, $k=1,...,r$, be two solutions to
(4) such that ${f_{kj}}_{\zbv S} = {g_{kj}}_{\zbv S}= \zf_{kj}$,
$k=1,...,r$, $j= 1,..., n-r$. On $S'$ functions $f_{kj}$ and
$g_{kj}$ are solutions to (7) which agree on $S$, so by lemma 5.1
we have $f_{kj} = g_{kj}$, $k=1,...,r$, $j = 1,..., n-r$, as germs
at $p$ on $S'$.

Now, like in the proof of theorem 4.1, we consider $x_1$ as a
new parameter. Let $J'_{0}$ be the restriction of $J_{0}$ to
${\mathcal F}'$ (recall that $dx_{1}\zci J_{0}=0$); then
$d'x_{j}\zci J'_{0} = a_{j}d'x_{j}$, $j=2,..., n-r$, $d'y_{k}\zci
J'_{0}=ad'y_{k}$, $k=1,...,r$.

Since $S'$ plays the same role with respect to $(x_{2},...,
x_{n-r}, y_{1},...,y_{r})$ as $S$ does with respect to
$(x_{1},..., x_{n-r}, y_{1},...,y_{r})$,
${\tilde\za}'_{k}=\zsu_{j=2}^{n-r}f_{kj}d'x_{j}$ and
$\zg'_{k}=\zsu_{j=2}^{n-r}g_{kj}d'x_{j}$ satisfy to system (4) of
theorem 5.1 for ${\mathcal F'}$ and $J'_{0}$, and ${\tilde\za}'
=\zg'$ on $S'$, from the induction hypothesis follows that $f_{kj}
= g_{kj}$, $k=1,...,r$,  $j = 2,..., n-r$, like germs at $p$ on
$M$.

Finally, as each $({\tilde\za}_{k} -\zg_{k})\zex dy_{1}\zex...\zex
dy_{r} = (f_{k1}-g_{k1})dx_1\zex dy_{1}\zex...\zex dy_{r}$ is
closed, function $f_{k1}-g_{k1}$ is constant along the leaves of
the foliation $Ker(d'y_{1}\zex...\zex d'y_{r}) \zco {\mathcal
F'}$. But $S$ is transverse to this foliation  and
$(f_{k1}-g_{k1})_{\zbv S}=0$ then $f_{k1}=g_{k1}$ and
${\tilde\za}_{k}=\zg_{k}$, $k=1,...,r$, as germs at $p$ on $M$.

For the existence we will need the following result.
\bigskip

{\bf Lemma 5.2.} {\it Consider $1$-forms $\zb_{1},...,\zb_{r}$
functional combination of $dx_{1},...,dx_{n-r}$. Let $G$ be the
$(1,1)$ tensor field along $\mathcal F$ defined by $dx_{j}\zci
G=a_{j}dx_{j}$, $j=1,...,n-r$, $dy_{k}\zci G=\zb_{k}+ady_{k}$,
$k=1,...,r$. Assume that $d\zb_{k}\zex dy_{1}\zex...\zex
dy_{r}=0$, $k=1,...,r$. Then $N_{G}=0$ if and only if

\centerline{$\zpizq d(\zb_{k}\zci J_{0})-\sum_{\zlma=1}^{r}
{\zb}_{\zlma}\zex {\frac{\zpar {\zb}_{k}} {\zpar
y_{\zlma}}}\zpder\zex dy_{1}\zex...\zex dy_{r}=0$, $k=1,...,r$.}}
\bigskip

{\bf Proof.} By lemma 2.1 one has $dx_{j}\zci N_{G}=0$ and

\noindent $dy_{k}\zci N_{G}=(d\zb_{k})_{G}-d(\zb_{k}\zci
J_{0}+a\zb_{k})=\zpizq \sum_{\zlma=1}^{r} dy_{\zlma}\zex
{\frac{\zpar {\zb}_{k}} {\zpar
y_{\zlma}}}\zpder_{G}-d_{x}(\zb_{k}\zci J_{0})-\sum_{\zlma=1}^{r}
dy_{\zlma}\zex ({\frac{\zpar {\zb}_{k}} {\zpar y_{\zlma}}}\zci
J_{0})-\sum_{\zlma=1}^{r}ady_{\zlma}\zex {\frac{\zpar {\zb}_{k}}
{\zpar y_{\zlma}}}=\sum_{\zlma=1}^{r} {\zb}_{\zlma}\zex
{\frac{\zpar {\zb}_{k}} {\zpar y_{\zlma}}}-d_{x}(\zb_{k}\zci
J_{0})$

\noindent where $d_{x}$ denotes the exterior derivative in
variables $(x_{1},...,x_{n-r})$ only. $\square$
\medskip

{\bf Existence.}  Given functions $\zf_{kj}$, $k=1,...,r$,
$j=1,...,n-r$, on $S$ such that every $\zf_{1j}(p)$ is a positive
real number, by means of system (7) we extend them to $S'$, around
$p$, with the same name.

    If we consider ${\mathcal F'}$ and $J'_{0}$, the induction
hypothesis allows us to find functions $f_{kj}$, $k=1,...,r$,
$j=2,...,n-r$, defined on an open neighbourhood of $p$ on $M$, in
such a way that $ d'{\tilde\za'}_{k}\zex d'y_{1}\zex...\zex
d'y_{r}=0$, $k=1,...,r$,

$\,$

\centerline{$\zgran\zpizq d'({\tilde\za'}_{k}\zci
J'_{0})-\sum_{\zlma=1}^{r} {\tilde\za'}_{\zlma}\zex {\frac{\zpar
{\tilde\za'}_{k}} {\zpar y_{\zlma}}}\zpder\zex d'y_{1}\zex...\zex
d'y_{r}=0, \quad k=1,...,r$,}

$\,$

\noindent and ${f_{kj}}=\zf_{kj}$ on $S'$, $k=1,...,r$, $j=2,...,
n-r$, where ${\tilde\za}'_{k} =\zsu_{j=2}^{n-r}f_{kj}d'x_{j}$
(note that ${\partial} /{\partial}x_{2},...,
{\partial}/{\partial}x_{n-r},
{\partial}/{\partial}y_{1},...,{\partial}/{\partial}y_{r}$ is the
dual basis of $d'x_{2},...,d'x_{n-r}$, $d'y_{1},...,d'y_{r}$ as
well).

Since each $d'{\tilde\za}'_{k}\zex d'y_{1}\zex...\zex d'y_{r}=0$
there exist functions $f_{k}$ such that $\zr_{k}\zex
dy_{1}\zex...\zex dy_{r}$ is closed where $\zr_{k}
=f_{k}dx_{1}+f_{k2}dx_{2}+...+f_{kn-r}dx_{n-r}$. On the other hand
the first equations of (7) means that every
$((\zf_{k1}dx_{1}+...+\zf_{kn-r}dx_{n-r})\zex dy_{1}\zex...\zex
dy_{r} )_{\zbv {\mathcal F}\zin S'}$ is closed. Therefore
$(\zr_{k}\zex dy_{1}\zex...\zex dy_{r})_{\zbv {\mathcal F}\zin S'}
-((\zf_{k1}dx_{1}+...+\zf_{kn-r}dx_{n-r})\zex dy_{1}\zex...\zex
dy_{r} )_{\zbv {\mathcal F}\zin S'} = ((f_{k}
-\zf_{k1}){dx_{1}\zex dy_{1}\zex...\zex dy_{r}})_{\zbv {\mathcal
F}\zin S'}$ has to be closed.

Consider coordinates $(x_{1},..., x_{n-r},y_{1},...,y_{r},
v_{1},..., v_{m-n})$, around $p$ on $M$, where $v_{1},..., v_{m-n}$
are basic functions for ${\mathcal F}$. Then, always
around $p$ on $M$, there exist functions ${\bar h}_{k}(x_{1}, y_{1},...,y_{r},
v_{1},..., v_{m-n})$ such that $f_{k}-\zf_{k1}={\bar h}_{k}$ on
$S'$. Now by setting $f_{k1} =f_{k}-{\bar h}_{k}$, we construct
1-form ${\tilde\za}_{k}=\zsu_{j=1}^{n-r}f_{kj}dx_{j}$,
$k=1,...,r$, along ${\mathcal F}$ such that ${f_{kj}}_{\zbv
S'}=\zf_{kj}$, $j=1,..., n-r$, $d{\tilde\za}_{k}\zex
dy_{1}\zex...\zex dy_{r}=0$ and

$\,$

\centerline{$\zgran\zpizq d({\tilde\za}_{k}\zci
J_{0})-\sum_{\zlma=1}^{r} {\tilde\za}_{\zlma}\zex {\frac{\zpar
{\tilde\za}_{k}} {\zpar y_{\zlma}}}\zpder\zex dx_{1}\zex
dy_{1}\zex...\zex dy_{r}=0$.}

$\,$

Therefore we can find 1-forms $\zg_{k},\zg_{k1},...\zg_{kr}$,
$k=1,...,r$, along ${\mathcal F}$, where each $\zg_{k}$ is closed
because $\zpizq d({\tilde\za}_{k}\zci J_{0})-\sum_{\zlma=1}^{r}
{\tilde\za}_{\zlma}\zex {\frac{\zpar {\tilde\za}_{k}} {\zpar
y_{\zlma}}}\zpder\zex dy_{1}\zex...\zex dy_{r}$ is closed since
$d{\tilde\za}_{k}\zex dy_{1}\zex...\zex dy_{r}=0$, such that

$\,$

\centerline{$\zgran\zpizq d({\tilde\za}_{k}\zci
J_{0})-\sum_{\zlma=1}^{r} {\tilde\za}_{\zlma}\zex {\frac{\zpar
{\tilde\za}_{k}} {\zpar
y_{\zlma}}}\zpder=dx_{1}\zex\zg_{k}+\zg_{k1}\zex
dy_{1}+...+\zg_{kr}\zex dy_{r}$.}

$\,$

Hence

$\,$

\noindent (8) \hskip .5truecm
 $\cases{\zgran d({\tilde\za}_{k}\zci
J_{0})=dx_{1}\zex\zg_{k}+\zg_{k1}\zex dy_{1}+...+\zg_{kr}\zex
dy_{r}+\sum_{\zlma=1}^{r} {\tilde\za}_{\zlma}\zex {\frac{\zpar
{\tilde\za}_{k}} {\zpar y_{\zlma}}}\cr k=1,...,r\cr}$

$\,$

Set $\zg_{k} = \zsu_{j=1}^{n-r}g_{kj}dx_{j} +
\zsu_{\zlma=1}^{r}g_{kn-r+\zlma}dy_{\zlma}$. Then

$\,$

\centerline {$\zgran g_{kj}=a_{j}{{\partial}f_{kj}\over
{\partial}x_{1}} +\sum_{\zlma=1}^{r}( f_{\zlma
j}{{\partial}f_{k1}\over {\partial}y_{\zlma}} - f_{\zlma
1}{{\partial}f_{kj}\over {\partial}y_{\zlma}})=0$, $j = 2,...,
n-r$, $k=1,...,r$. }

$\,$

\noindent (recall the construction of system (5)). So each
$g_{kj}$, $k=1,...,r$, $j = 2,..., n-r$, vanishes on $S'$ because
functions ${f_{kj}}_{\zbv S'}=\zf_{kj}$ satisfy to system (7).

Deriving (8) with respect to $y_{s}$ yields

$\,$

\noindent (9) \hskip .5truecm
 $\zgran d\zpizq{\frac{\zpar{\tilde\za}_{k}} {\zpar y_{s}}}\zci
J_{0}\zpder =dx_{1}\zex{\frac{\zpar\zg_{k}} {\zpar y_{s}}}
+\sum_{\zlma=1}^{r}{\frac{\zpar\zg_{k\zlma}} {\zpar y_{s}}}\zex
dy_{\zlma} +\sum_{\zlma=1}^{r}\zpizq {\frac{\zpar
{\tilde\za}_{\zlma}} {\zpar y_{s}}}\zex {\frac{\zpar
{\tilde\za}_{k}} {\zpar y_{\zlma}}}+ {\tilde\za}_{\zlma}\zex
{\frac{\zpar^{2} {\tilde\za}_{k}} {\zpar y_{s} \zpar
y_{\zlma}}}\zpder$

$\,$

On the other hand

 $\,$

\noindent (10) \hskip .5truecm
 $\zgran (d({\tilde\za}_{k}\zci
J_{0}))_{J_{0}}=dx_{1}\zex(\zg_{k}\zci J_{0})
+\sum_{\zlma=1}^{r}(\zg_{k\zlma}\zci J_{0}+a\zg_{k\zlma})\zex
dy_{\zlma}$

\hskip 4truecm$\zgran+\sum_{\zlma=1}^{r}
\zpizq({\tilde\za}_{\zlma}\zci J_{0})\zex {\frac{\zpar
{\tilde\za}_{k}} {\zpar y_{\zlma}}} +{\tilde\za}_{\zlma}\zex
({\frac{\zpar {\tilde\za}_{k}} {\zpar y_{\zlma}}}\zci
J_{0})\zpder$

$\,$

By lemma 2.1 applied along the leaves of ${\mathcal F}$ we have
$d((d({\tilde\za}_{k}\zci J_{0}))_{J_{0}})\zex dy_{1}\zex...\zex
dy_{r}=0$, whence by calculating $d((d({\tilde\za}_{k}\zci
J_{0}))_{J_{0}})$ from (10) and taking into account (8) and (9)
follows

$\,$

\noindent (11) \hskip .5truecm $\cases{\zgran \zpizq d(\zg_{k}\zci
J_{0})+\zsu_{\zlma=1}^{r}{{\partial}\zg_{k}\over{{\partial}y_{\zlma}}}
\zex{\tilde\za}_{\zlma} -\zsu_{\zlma=1}^{r} \zg_{\zlma}\zex
{{{\partial}{\tilde\za}_{k}}\over{{\partial}y_{\zlma}}}\zpder\zex
dx_{1}\zex dy_{1}\zex...\zex dy_{r}=0\cr k=1,...,r.\cr}$

$\,$

Set $\zg'_{k}=\zsu_{j=2}^{n-r}g_{kj}d'x_{j} +
\zsu_{\zlma=1}^{r}g_{kn-r+\zlma}d'y_{\zlma}$. Obviously
$d'\zg'_{k}=0$ because $d\zg_{k}=0$. Consider the $(1,1)$ tensor
field $J'$ on ${\mathcal F'}$ defined by $d'x_{j}\zci J'
=a_{j}d'x_{j}$, $j=2,..., n-r$, and $d'y_{\zlma}\zci J'=
{\tilde\za}'_{\zlma}+ad'y_{\zlma}$, $\zlma=1,...,r$ (recall that
${\tilde\za}'_{\zlma} =\zsu_{j=2}^{n-r}f_{\zlma j}d'x_{j}$).

Since

$\,$

\centerline{$\zgran\zpizq d'({\tilde\za'}_{k}\zci
J'_{0})-\sum_{\zlma=1}^{r} {\tilde\za'}_{\zlma}\zex {\frac{\zpar
{\tilde\za'}_{k}} {\zpar y_{\zlma}}}\zpder\zex d'y_{1}\zex...\zex
d'y_{r}=0, \quad k=1,...,r$,}

$\,$

\noindent by lemma 5.2, applied to ${\mathcal F'}$ and $J'$, the
Nijenhuis torsion of $J'$ vanishes. Set $\zgran\zr_{\zlma
k}=-{{{\partial}{\tilde\za}'_{k}}\over {{\partial}y_{\zlma}}}$.
Now system (11) becomes (note that $\zgran
d'g_{kn-r+\zlma}={{\partial}\zg'_{k}\over{{\partial}y_{\zlma}}}$
because $\zg'_{k}$ is closed)

$\,$

\noindent (12) \hskip .5truecm $\zgran\zpizq d'(\zg'_{k}\zci
J')+\sum_{\zlma=1}^{r} \zg'_{\zlma}\zex \zr_{\zlma k}\zpder\zex
d'y_{1}\zex...\zex d'y_{r}=0, \quad k=1,...,r$.

$\,$

On $S'$, $\zg'_{k}\zex dy_{1}\zex...\zex dy_{r}=0$ as
${g_{kj}}_{\zbv S'}=0$, $k=1,...,r$, $j=2,...,n-r$. Since the
restriction of $\zg'_{k}$ to ${\mathcal F'}\zin S'$ is closed,
around $p$ on $M$ there exist functions
$\zfi_{k\zlma}(x_{1},y_{1},...,y_{r}, v_{1},..., v_{m-n})$,
 $k,\zlma=1,...,r$, such that every $\zg'_{k}=
\zsu_{\zlma=1}^{r}\zfi_{k\zlma}d'y_{\zlma}$ on $S'$.

Set $\zl_{k}=\zsu_{\zlma=1}^{r}\zfi_{k\zlma}d'y_{\zlma}$. Then
each $\zl_{k}$ is a closed 1-form along ${\mathcal F'}$ defined on
an open neighbourhood of $p$ on $M$ and $\zpizq d'(\zl_{k}\zci
J')+\sum_{\zlma=1}^{r} \zl_{\zlma}\zex \zr_{\zlma k}\zpder\zex
d'y_{1}\zex...\zex d'y_{r}=0, \quad k=1,...,r$. Now lemma 4.2
applied to ${\mathcal F'}$ and $J'$ implies that
$\zg'_{k}=\zl_{k}$, $k=1,...,r$. In other words every $\zg_k$ is a
functional combination of $dx_{1},dy_{1},...,dy_{r}$. Therefore

$\,$

\centerline{$\zgran \zpizq d({\tilde\za}_{k}\zci
J_{0})-\sum_{\zlma=1}^{r} {\tilde\za}_{\zlma}\zex {\frac{\zpar
{\tilde\za}_{k}} {\zpar y_{\zlma}}}\zpder\zex dy_{1}\zex...\zex
dy_{r}=0$, $k=1,...,r$,}

$\,$

\noindent and {\it the proof of theorem 5.1 is finished}.
\bigskip

{\bf 6. Local classification of codimension one Veronese webs.}

On a real or complex manifold $N$ of dimension $n$ consider a
Veronese web $w$ of codimension $r\zmai 1$. Given non-equal
scalars $a_{1},...,a_{n-r},a$ and any point $p\zpe N$, let $J$ be
a $(1,1)$ tensor field like in part $(1)$ of theorem 2.1 and let
$(x_{1},...,x_{n-r},y_{1},...,y_{r})$ be a system of coordinates,
around $p$, such that $dx_{j}\zci J=a_{j}dx_{j}$, $j=1,...,n-r$,
and $Ker(dy_{1}\zex...\zex dy_{r})=w(\zinf)$. Then $dy_{k}\zci
J=ady_{k}+{\tilde\za}_{k}$, $k=1,...,r$, where each
${\tilde\za}_{k}=\zsu_{j=1}^{n-r}f_{kj}dx_{j}$. As
$(w(\zinf)',J^{*})$ spans the cotangent bundle around $p$, by
linearly recombining functions $y_{1},..., y_{r}$ and considering
$b_{j}x_{j}$ instead $x_j$ for a suitable $b_{j}\zpe{\mathbb
K}-\{0\}$, we assume that each $f_{1j}(p)$, $j = 1,..., n-r$, is a
positive real number (see the beginning of section 4).

On the other hand $d(dy_{k}\zci J)\zex dy_{1}\zex...\zex dy_{r}=0$
and $N_{J}=0$; by lemma 5.2 these last two conditions are
equivalent to system

$\,$

\noindent (13) \hskip .5truecm
 $\cases{ d{\tilde\za}_{k}\zex
dy_{1}\zex...\zex dy_{r}=0,\quad k=1,...,r\cr \, \cr\zgran \zpizq
d({\tilde\za}_{k}\zci J_{0})-\sum_{\zlma=1}^{r}
{\tilde\za}_{\zlma}\zex {\frac{\zpar {\tilde\za}_{k}} {\zpar
y_{\zlma}}}\zpder\zex dy_{1}\zex...\zex dy_{r}=0, \quad
k=1,...,r\cr}$

$\,$

\noindent where $J_{0}=\zsu_{j=1}^{n-r}a_{j}{\frac{\zpar} {\zpar
x_{j}}}\zte dx_{j}
+\zsu_{\zlma=1}^{r}a{\frac{\zpar} {\zpar
y_{\zlma}}}\zte dy_{\zlma}$.

Moreover
$\zg(t)=(\zpr_{j=1}^{n-r}(t+a_{j}))(t+a)^{r}((J+tI)^{-1})^{*}
(dy_{1}\zex...\zex dy_{r})$ represents $w$.

Therefore, in view of $(3)$ of theorem 2.1, locally Veronese webs
correspond to those solutions of system $(13)$ such that
$f_{11}(p),...,f_{1n-r}(p)\zpe {\mathbb R}^{+}$ (this last
assumption implies that $(dy_{1},...,dy_{r},J^{*})$ spans the
cotangent bundle near $p$). In turn, for the standard case, this
kind of solutions to $(13)$ are given by theorem 5.1 by setting
$M=N$ and ${\mathcal F}=TN$, which means that now $S$ is the
submanifold defined by $x_{j}-x_{n-r}=x_{j}(p)-x_{n-r}(p)$,
$j=1,...,n-r-1$.

When $r\zmai 2$ the tensor field $J$ is not unique and
consequently we may associate more than one model to a same
Veronese web; {\it thus our model of every Veronese web is
versal}.

To remark that a classification in codimension$\zmai 2$ seems
rather difficult as the following example shows. Consider a field
of $2$-planes and a local basis of it $\{ X,Y\}$. Let ${\tilde
w}(t)$, $t\zpe \mathbb K$, be the $1$-foliation defined by $X+tY$.
Then to classify the $1$-dimensional (local) Veronese web ${\tilde
w}=\{ {\tilde w}(t)\zbv t\zpe {\mathbb K}\}$, roughly speaking, is
like locally classifying the fields of $2$-planes in any
dimension; but it is well known the difficult of this problem
(first dealt with by \'Elie Cartan in ''Les syst\`emes de Pfaff \`a
cinq variables'' and later on by several authors).

Now let us examine the remainder case. {\it Assume $r=1$ until the
end of this section.} Then $a_{1},...,a_{n-1},a$ completely
determines $J$ since $Ker(J^{*}-a_{j}I)$, $j=1,...,n-1$, is the
annihilator of $w(-a_{j})$ and $Ker(J^{*}-aI)$ that of $w(-a)$.
The next step will be to construct an intrinsic surface $S$. By
technical reasons {\it one will suppose that $a_{1},...,a_{n-1},a$
are non-equal real numbers}.

The polynomial $\zsu_{j=1}^{n-1}\zpr_{k=1;k\znoi
j}^{n-1}(t+a_{k})$ has $n-2$ different roots $b_{1},...,b_{n-2}$
since it is the derivative of $\zpr_{k=1}^{n-1}(t+a_{k})$, whose
roots are $-a_{1},...,-a_{n-1}$; moreover $b_{\zlma}\znoi -a_{j}$,
$\zlma=1,...,n-2$, $j=1,...,n-1$ (warning this property is not
true when a polynomial, even real, has some complex root, for
example $t^{3}-1$ and $3t^{2}$; by this reason one chooses real
numbers $a_{1},...,a_{n-1},a$).

Let $R$ be the germ at $p$ of the leaf of the $1$-foliation
$w(b_{1})\zin ...\zin w(b_{n-2})\zin w(\zinf)$ passing through
this point, and let $S_{0}$ be the germ at $p$ of the surface
containing $R$ and to which the $1$-foliation $w(-a_{1})\zin
...\zin w(-a_{n-1})$ is tangent. By construction $S_{0}$ is
intrinsic.

Since $R$ is transverse to every $w(-a_{j})$, $j=1,...,n-1$, one
may take coordinates $(x_{1},...,x_{n-1},y)$ constructed before,
with two additional properties: $R$ is defined by the equations
$x_{1}=...=x_{n-1}, y=0$, and $x_{1}(p)=...=x_{n-1}(p)=y(p)=0$; of
course we write $y$ and ${\tilde\za}=\zsu_{j=1}^{n-r}f_{j}dx_{j}$
instead $y_{1}$ and
${\tilde\za}_{1}=\zsu_{j=1}^{n-r}f_{1j}dx_{j}$. In these
coordinates $S_{0}$ is defined by the equations
$x_{1}=...=x_{n-1}$. Moreover

$\,$

\centerline{$\zgran\zg(t)=-\zsu_{j=1}^{n-1}\zpizq\zpr_{k=1;k\neq
j}^{n-1}(t+a_{k})f_{j}\zpder dx_{j} +
\zpr_{k=1}^{n-1}(t+a_{k})dy$}

$\,$

\noindent because a straightforward calculation shows that

$\,$

 \centerline{$\zgran\zpizq -\zsu_{j=1}^{n-1}(\zpr_{k=1;k\neq
j}^{n-1}(t+a_{k})f_{j})dx_{j} + \zpr_{k=1}^{n-1}(t+a_{k})dy \zpder
\zci(J+tI)= \zpizq\zpr_{k=1}^{n-1}(t+a_{k})(t+a)\zpder dy$.}

$\,$

On the other hand $\zg(b_{\zlma})(q)((\zpar /\zpar
x_{1})+...+(\zpar /\zpar x_{n-1}))=0$, $\zlma=1,...,n-2$, for
every $q\zpe R$ because $(\zpar /\zpar x_{1})+...+(\zpar /\zpar
x_{n-1})$ is tangent to $R$ and $T_{q}R=(w(b_{1})\zin ...\zin
w(b_{n-2})\zin w(\zinf))(q)$. Therefore $b_{1},...,b_{n-2}$ are
the roots of $\zsu_{j=1}^{n-1}\zpr_{k=1;k\neq
j}^{n-1}(t+a_{k})f_{j}(q)$ when $q\zpe R$; so $f_{1}=...=f_{n-1}$
on $R$ since $b_{1},...,b_{n-2}$ are the roots of
$\zsu_{j=1}^{n-1}\zpr_{k=1;k\neq j}^{n-1}(t+a_{k})$ too, which
implies that both polynomials are equal up to multiplicative
factor (conversely, if  $f_{1}=...=f_{n-1}$ on $R$ then $(\zpar
/\zpar x_{1})+...+(\zpar /\zpar x_{n-1})$ is tangent to this curve
and $R$ is defined by $x_{1}=...=x_{n-1}, y=0$).

The change of coordinates between two of such system can be
regarded as a diffeomorphism $(x_{1},...,x_{n-1},y)\zfl
G(x_{1},...,x_{n-1},y)$. But $G$ has to preserve $R$, $S_{0}$, the
foliations of dimension $n-1$ defined by $dx_{1}$,..., $dx_{n-1}$
and $dy$ respectively (that is to say $w(-a_{1})$,...,
$w(-a_{n-1})$ and $w(\zinf)$), and the origin. Therefore
$G(x_{1},...,x_{n-1},y)=(h_{1}(x_{1}),...,h_{1}(x_{n-1}),h_{2}(y))$
where $h_{1},h_{2}$ are one variable functions such that
$h_{1}(0)=h_{2}(0)=0$ and $h'_{1}(0)\znoi 0$, $h'_{2}(0)\znoi 0$.

Denote by $J'$ the pull-back of $J$ by the diffeomorphism $G$.
Then $dx_{j}\zci J'=a_{j}dx_{j}$, $j=1,...,n-1$, and $dy\zci
J'=ady+{\tilde\za}'$ where
\medskip

\centerline{$\zgran{\tilde\za}'=\zsu_{j=1}^{n-1}h'_{1}(x_{j})(h'_{2}(y))^{-1}
f_{j}(h_{1}(x_{1}),...,h_{1}(x_{n-1}),h_{2}(y))dx_{j}$.}
\medskip

Now we may take $h_{1},h_{2}$ in such a way that

\centerline{$h'_{1}(x_{1})(h'_{2}(y))^{-1}f_{1}(h_{1}(x_{1}),...,h_{1}
(x_{n-1}),h_{2}(y))=1$}

\noindent on the curves $x_{1}=...=x_{n-1}$, $y=0$, and
$x_{1}=...=x_{n-1}=0$. Indeed, first consider the function $h_{2}$
defined by $(h'_{2}(t))^{-1}f_{1}(0,...,0,h_{2}(t))=1$,
$h_{2}(0)=0$, and then the function $h_{1}$ defined by
$h'_{1}(t)(h'_{2}(0))^{-1}f_{1}(h_{1}(t),...,h_{1}(t),0)=1$,
$h_{1}(0)=0$; note that $h'_{1}(0)=1$ since

\centerline{$h'_{1}(0)(h'_{2}(0))^{-1}f_{1}(0,...,0,0)=
(h'_{2}(0))^{-1}f_{1}(0,...,0,0)=1$.}

In other words, there exist
coordinates $(x_{1},...,x_{n-1},y)$ as before with a third
additional property: $f_{1}=...=f_{n-1}=1$ on the curve
$x_{1}=...=x_{n-1}$, $y=0$, and $f_{1}=1$ on the curve
$x_{1}=...=x_{n-1}=0$.

In turn, a change of coordinates between two system with this last
property is given by two functions $h_{1},h_{2}$ such that
$h'_{1}(x_{1})(h'_{2}(y))^{-1}=1$ on the curves
$x_{1}=...=x_{n-1}$, $y=0$, and $x_{1}=...=x_{n-1}=0$. Therefore
$h'_{1},h'_{2}$ are constant. In short, the only possible change
of coordinates is a homothety by some $b\zpe{\mathbb K}-\{0\}$,
and ${\tilde\za}'(x_{1},...,x_{n-1},y)
={\tilde\za}(bx_{1},...,bx_{n-1},by)$.

A germ at the origin of a map $\zfi=(\zf_{1},...,\zf_{n-1})$ from
$S_{0}$ to ${\mathbb K}^{n-1}$ will be called {\it admissible} if
$\zf_{1}=...=\zf_{n-1}=1$ on the curve $x_{1}=...=x_{n-1}$, $y=0$,
and $\zf_{1}=1$ on the curve $x_{1}=...=x_{n-1}=0$. Two admissible
germs $\zfi$ and $\bar\zfi$ will be named {\it equivalent} if
there exists $b\zpe{\mathbb K}-\{0\}$ such that
${\bar\zfi}(x_{1},...,x_{n-1},y)=\zfi(bx_{1},...,bx_{n-1},by)$.

From theorem 2.1, theorem 5.1 and system (13), applied to the last
kind of coordinates system, follows (remark that in this last step
the number $a$ does not play any role, which is due to the fact
that a Veronese web is determined by $w(\zinf)$ and $J_{\zbv
w(\zinf)}$):
\bigskip

{\bf Theorem 6.1.} {\it Consider non-equal real numbers
$a_{1},...,a_{n-1}$. One has:

\noindent (1) Given a Veronese web of codimension $1$ on a real or
complex $n$-manifold $N$ and any point $\zpe N$, there exist
coordinates $(x_{1},...,x_{n-1},y)$ around $p$ such that
$x_{1}(p)=...=x_{n-1}(p)=y(p)=0$ and the Veronese web is
represented by

$\,$

\centerline{$\zgran\zg(t)=-\zsu_{j=1}^{n-1}\zpizq\zpr_{k=1;k\neq
j}^{n-1}(t+a_{k})f_{j}\zpder dx_{j} +
\zpr_{k=1}^{n-1}(t+a_{k})dy$,}

$\,$

\noindent where ${\tilde\za}=\zsu_{j=1}^{n-r}f_{j}dx_{j}$
satisfies to the system

$\,$

\noindent  \hskip .5truecm
 $\cases{ d{\tilde\za}\zex dy=0\cr \, \cr
 \zgran \zpizq d\zpizq\sum_{j=1}^{n-1}a_{j}f_{j}dx_{j}\zpder
-{\tilde\za}\zex {\frac{\zpar {\tilde\za}} {\zpar y}}\zpder\zex
dy=0,\cr}$

$\,$

\noindent $f_{1}=...=f_{n-1}=1$ on the curve $x_{1}=...=x_{n-1}$,
$y=0$, and $f_{1}=1$ on the curve $x_{1}=...=x_{n-1}=0$.

\noindent (2) Let $S_{0}$ be the surface of equation
$x_{1}=...=x_{n-1}$ and let $\zfi=(\zf_{1},...,\zf_{n-1})$ be a
germ at the origin of a map from $S_{0}$ to ${\mathbb K}^{n-1}$.
Assume $\zfi$ admissible. Then there exists one and only one germ
at the origin of $1$-form
${\tilde\za}=\zsu_{j=1}^{n-r}f_{j}dx_{j}$, which satisfies to the
system of part (1) and such that ${f_{j}}_{\zbv S_{0}}=\zf_{j}$,
$j=1,...,n-1$.

Moreover

$\,$

\centerline{$\zgran\zg(t)=-\zsu_{j=1}^{n-1}\zpizq\zpr_{k=1;k\neq
j}^{n-1}(t+a_{k})f_{j}\zpder dx_{j} +
\zpr_{k=1}^{n-1}(t+a_{k})dy$,}

$\,$

\noindent defines a Veronese web of codimension $1$ around the
origin.

\noindent (3) Finally given two admissible germs at the origin
$\zfi$ and $\bar\zfi$ of maps from $S_{0}$ to ${\mathbb
K}^{n-1}$, the germs of $1$-codimensional Veronese webs associated
to them by virtue of part (2) are equivalent, by diffeomorphism,
if and only if $\zfi$ and $\bar\zfi$ are equivalent as admissible
germs.}
\bigskip

The local classification of Veronese webs of codimension $1$ is
due to Turiel (see \cite{TUB} whose exposition is closely followed here).

\end{document}